\documentclass[12]{ws-jktr}

\usepackage{amsfonts,amssymb,wasysym,wrapfig}
\usepackage{graphicx}

\newcommand{\R}{\mathbb{R}}

\pdfoutput=1

\title{Non-orientable surfaces in $4$-dimensional space}

\author{Yongju Bae\\J. Scott Carter \\Seonmi Choi \\Sera Kim  }

\begin{document}
\maketitle

\begin{abstract} This article is a survey article that gives detailed constructions and illustrations of some of the standard examples of non-orientable surfaces that are embedded and immersed in $4$-dimensional space. The illustrations depend upon their $3$-dimensional projections, and indeed the illustrations here depend upon a further projection into the plane of the page. The concepts used to develop the illustrations will be developed herein.
\end{abstract}

This article falls into the category of ``descriptive topology." It consists of detailed constructions of  some of the standard examples of non-orientable surfaces that are embedded (and immersed) in ordinary $4$-dimensional space. Many of the examples here are commonly known throughout the literature. For example, George Francis's book \cite{GF} contains elaborate pen and ink drawings of Boy's surface, the cross-cap, the Roman surface, and the Klein bottle. The surface that has been dubbed ``girl's surface" appears in Ap\'{e}ry \cite{Ap} and in a more recent paper by Goodman and Kossowski \cite{GW}. Other surfaces, such as the connected sum of the cross-cap and the $3$-twist-spun trefoil, are lesser known and are the subject of some intensive research agendas. Our hope in this survey is to present as detailed descriptions as possible using print media and $2$-dimensional drawings. The original illustrations are in full color as can be seen in the internet based versions of the article or  the prepublication version that is available on the arXiv.

 The embeddings and immersions are always smooth (of class ${\mathcal C}^\infty$). The descriptions and constructions are given by means of movies, charts, diagrams, or decker sets. An informal description of these ideas is given in \cite{R2B}, and it is discussed more formally in \cite{CS:Book}. In addition to simply describing such surfaces, we also will demonstrate that some of the embeddings are ambiently isotopic while others are not. 
 
 Section~\ref{prelim} gives preliminary results on movie and chart descriptions. A full restatement of the movie-move theorem \cite{CS:MM,CRS} would take us too far afield; so in the subsequent sections when a movie-move or its associated chart move is applied, we will refer to the standard sources \cite{CS:Book} and indicate the page number in which the move is described. 
The discussion of non-orientable surfaces begins in earnest in Section~\ref{peetwo} with examples of the projective plane. These include the two different cross-cap embeddings whose normal Euler classes are $\pm2$. The non-trivial values of the normal class indicate that neither can projected into $3$-space without having branch points. In contrast, Boy's surface is an immersion of the projective plane into $3$-space that is not the projection of an embedding in $4$-space. It has a sister that is informally called {\it girl's surface}; it has a distinct double decker structure from Boy's surface. The {\it double decker set} of a generic map of a surface in $3$-space is the lift of the double points, branch points, and triple points to the surface that is being mapped. When the surface has been further projected to the plane, it has folds and cusps. The fold set also can be lifted to the ambient surface, and we will include the pre-images of the folds and cusps in the {\it decker set}. The double decker set, then, is a subset of the decker set.

The connected-sum of several cross-caps can have normal Euler classes that range over a specific set that depends on the number of connect-summands. The number of connect-summands is called the {\it non-orientable genus}. Whitney conjectured, and Massey proved, that any embedded projective plane in $4$-space of non-orientable genus $g$ has a normal Euler class that takes values in the set $\{-2g, -2g+4, \ldots , 2g-4, 2g \}$. Massey proved this result using the Atyiah-Singer Index Theorem. A more elementary proof is found in \cite{Kamada}. The normal Euler class is discussed in subsection~\ref{MTWC}.

Section~\ref{Klein} demonstrates that the connected sum of a positive and a negative cross-cap is ambiently isotopic to a standard Klein bottle which is ambiently isotopic to the twisted figure-8 Klein bottle. 
In Section~\ref{knotted}, we give examples of projective planes that are knotted. We demonstrate knottedness by using coloring invariants. These examples are constructed explicitly as the connected sum of cross-caps and knotted spheres. A fundamental problem in $4$-dimensional knot theory is to determine if other examples exist. That is, ``Does every knotted projective plane arise as the connected-sum of a standard projective plane and  a knotted $2$-sphere?" This is conjectures to be true and is known as ``Kinoshita's conjecture."

\subsection*{Acknowledgements} Some of the research for this talk was supported by the Ministry of Education Science and Technology (MEST) and the Korean Federation of Science and Technology Societies (KOFST). Additional support came from NIMS. We gratefully acknowledge a number of useful conversations with Seiichi Kamada and Masahico Saito. The members of the TAPU-KOOK seminar have also been quite helpful. The referee's suggestions due to the careful reading of the manuscripts improved the paper greatly.

\section{Preliminaries}
\label{prelim}

The technique of movie descriptions of surfaces that are embedded in $4$-space are at least as old as Fox's seminal paper {\it A quick trip through knot theory} \cite{Fox}. When a surface is smoothly embedded in $4$-space, the generic intersection with a $3$-dimensional cross-sectional hyperplane is an embedded closed curve in space. We can project such a cross-section further to the plane, and generically such a projection will have no cusps, no self-tangencies, and no transverse multiple points of multiplicity greater than two. An additional projection to a ``vertical axis" gives a Morse function that is called the {\it height function for the curve}. At this final level, we can assume that all the critical points are non-degenerate ($f^{\prime\prime}(c)\ne0$), and that double points and critical points all occur at distinct levels. 

\begin{wrapfigure}[19]{l}{3.32in}\vspace{-.32in}\includegraphics[width=3.3in]{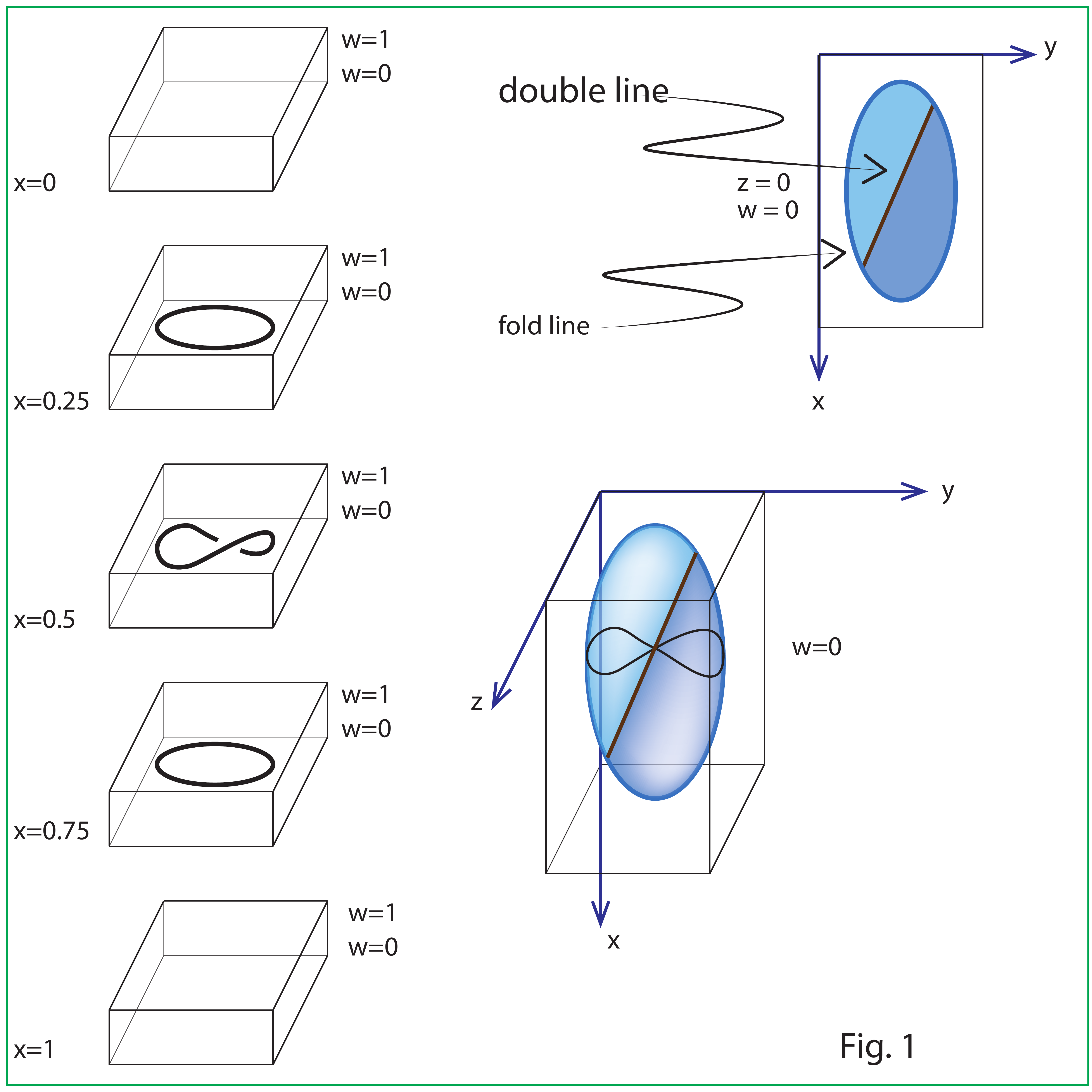}\end{wrapfigure}
\begin{sloppypar}
The situations that are excluded in the previous paragraph correspond to certain co\-dim\-en\-sion-$1$ singularities in the space of all maps from a curve to the plane. As such, they each have exactly one pair of distinct resolutions. The resolution of the cusp singularity corresponds to a type-I Reidemeister move; the resolution of a self-tangency corresponds to a type-II move; and the resolution of a triple point corresponds to a type-III move. Similarly, resolutions of the height function correspond to moves of knot diagrams. Specific details will follow shortly.\end{sloppypar}

A surface embedded in $4$-space, then, can be decomposed as a sequence of curves that are intersections with $3$-dimensional cross-sections. Successive generic cross-sections differ by codimension-$1$ singularities (that were given above) and critical points (Maxima, minima, and saddle points) for surfaces. Now we turn to a more technical description of the situation.

Let $F$ denote a closed surface (compact without boundary). Suppose that $K:F \hookrightarrow \R^4$ is a smooth embedding. By compactness and by rescaling if necessary, we assume that the image $K(F)$ lies in the unit hypercube $[0,1]\times [0,1] \times [0,1] \times [0,1] = \{ (x,y,z,w): \ 0 \le x,y,z,w \le 1 \}$. The embedding can be perturbed slightly (again if this is necessary) so that the composition $p_w \circ F$ is a generic surface in the cube $[0,1]\times [0,1] \times [0,1] \times \{0 \}$. Here $p_w: [0,1]^4 \rightarrow [0,1]^3\times \{0\}$ is the orthogonal projection of the hypercube to the unit $3$-dimensional cube.

A generic surface in $3$-space has isolated branch points, arcs of transverse double points, and isolated triple points. A {\it branch point} is a point on the surface for which any sufficiently small neighborhood intersects the surface in a cone on a figure-$8$ curve. In a neighborhood of a double point the surface is homeomorphic to the image $\{ (x,y,z): \   xy=0  \}$. In a neighborhood of a triple point the surface is homeomorphic to the image $\{(x,y,z): \ xyz=0 \}$.
The composition, $p_w\circ K$, is a generic surface in $3$-space that possibly has branch points, double point arcs, and triple points.

 \begin{wrapfigure}[19]{r}{3.4in}\vspace{-0.32in}\includegraphics[width=3.3in]{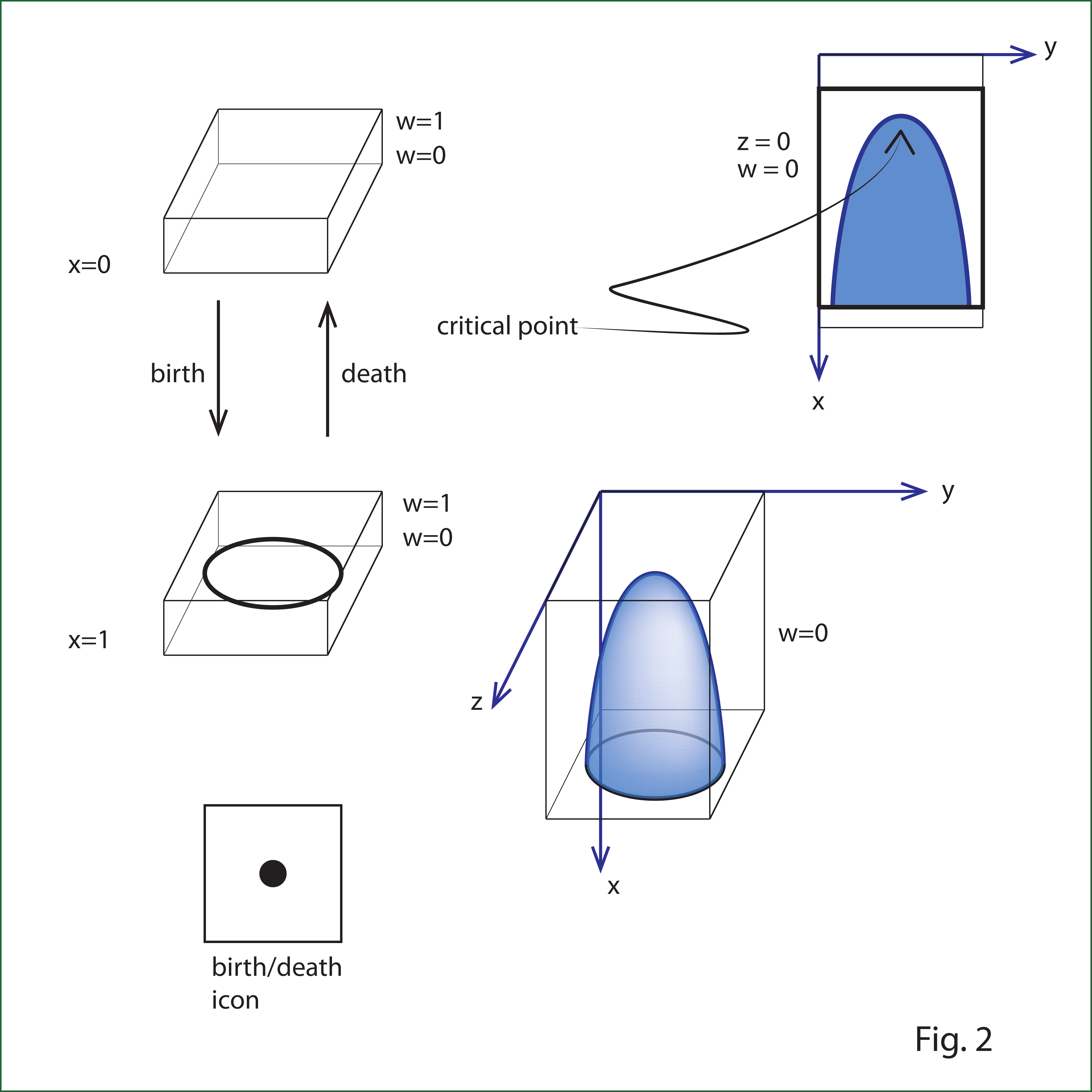}\end{wrapfigure} 
Let $p_{wz}: [0,1]^4 \rightarrow [0,1]^2 \times \{0\}\times \{0\}$ denote the orthogonal projection
$(x,y,z,w)\mapsto (x,y,0,0)$ to the unit $(x,y)$-square. Again the embedding $K$ may be perturbed slightly so that the image of $F$ projects with cusp and fold singularities. Meanwhile, for all but a finite number of values of $x \in [0,1]$, the intersection $K(F) \cap (\{x\} \times [0,1]^3)$ is an embedded curve (knot or link) in $3$-space.   Such a generic curve is called {\it a still of  a movie}. 

 \begin{wrapfigure}[19]{l}{3.4in}\vspace{-0.32in} \includegraphics[width=3.3in]{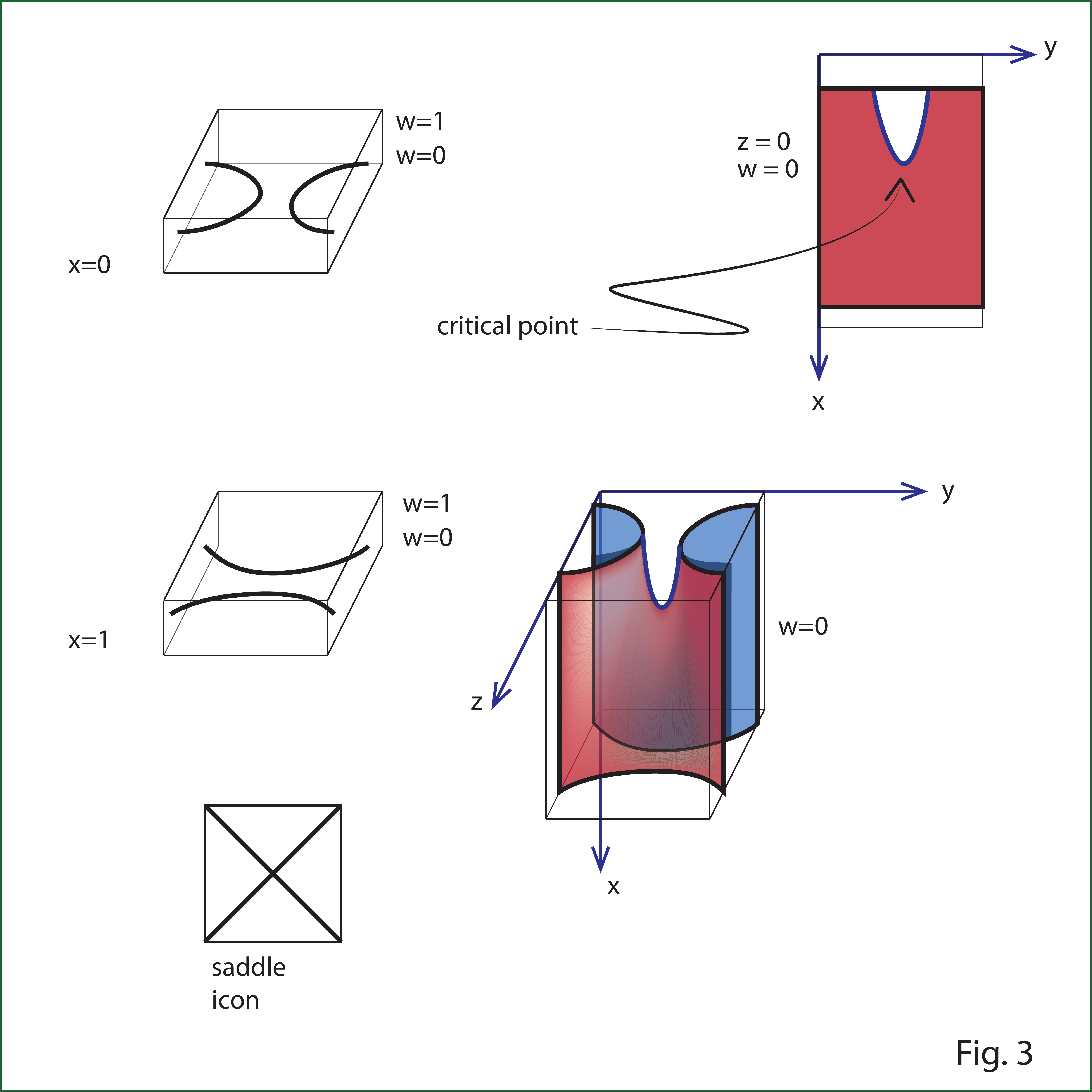}\end{wrapfigure} The {\it diagram of a still} is obtained by projecting to the ($x=$ constant, $w=0$)-plane and by denoting crossing information by means of the standard method.  Fixing $x=c,$ consider the projection of the still to the interval $\{ c\} \times [0,1] \times \{0\} \times \{0 \}$.  The surface can be put into general position so that this projection is a Morse function for the curve. More specifically, the double points usually project to distinct levels as do the non-degenerate critical points. And double points and critical points usually do not coincide.

  Above in Fig. 1, we illustrate a particularly  simple sphere  that is embedded \begin{wrapfigure}[19]{l}{3.4in}\vspace{-0.32in} \includegraphics[width=3.3in]{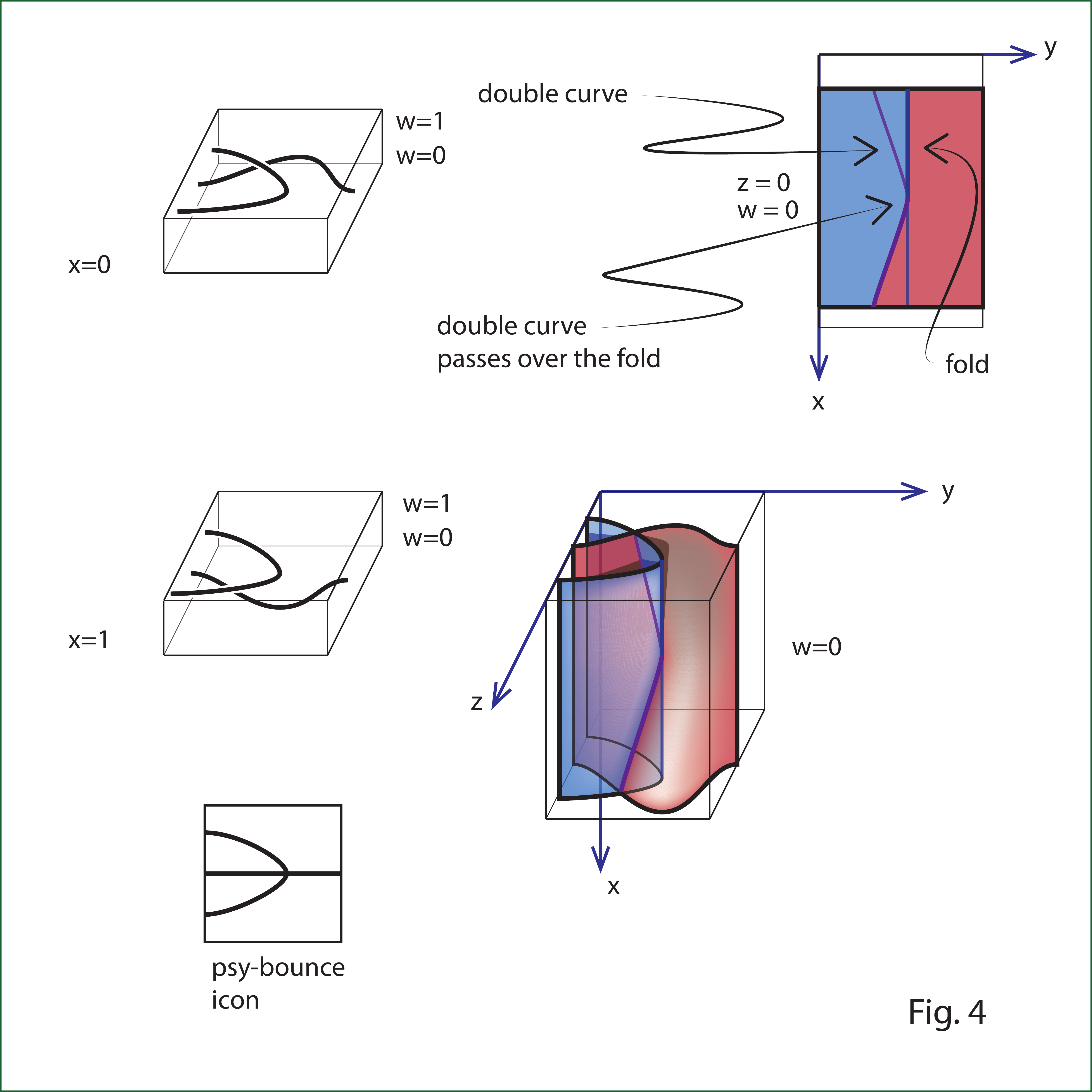}\end{wrapfigure} in $4$-space by means of its {\it movie} --- sequence of generic stills that differ by at most a single critical event (on the left of this illustration), in the center of the figure the projection is illustrated, and on the upper right its  chart is illustrated.   The chart is coincident with the projection of the surface onto the $(x,y,0,0)$-plane. In more generality, the chart is a labeled graph with two types of edges (double points and folds) and several types of vertices; these correspond to the singularities that have been discussed above.

The images of these singularities under projection are depicted. In Fig. 2,  
a {\it birth or death} of a simple closed curve is depicted as a ($2$-dimensional) cap or cup. Such a critical point yields an optimum on the fold set of the projection into the $(x,y)$-plane. 

Fig. 3 indicates a {\it saddle} which also yields an optimum on the fold set. Note that the critical points in the $y$-direction form the folds. In subsequent figures, a minimal point in the $y$-direction (a left pointing optimum) will be denoted by a subscripted lower case $m$ while a {Maximal{\footnote {We adopt the affectation that the word ``Maximal" always begins with an upper case M while ``minimal" begins with a lower-case m.}}}  point (a right pointing optimum) is denoted by a subscripted $M$.

The {\it $\psi$-bounce} or simply {\it bounce} singularity occurs when a double curve passes over a fold (Fig. 4). In the projection, the image of the double curve and the fold curve become tangent --- thus the name ``bounce." The singular situation involves two intersecting sheets in which one passes directly through the vertex of a parabola in the same manner as the line $x=1/2$, $z=1/2, w=0$ passes through the parabola $x=1/2$, $w=0$, $y=1/2-(z-1/2)^2.$ 
 
 \begin{wrapfigure}[19]{r}{3.4in} \includegraphics[width=3.3in]{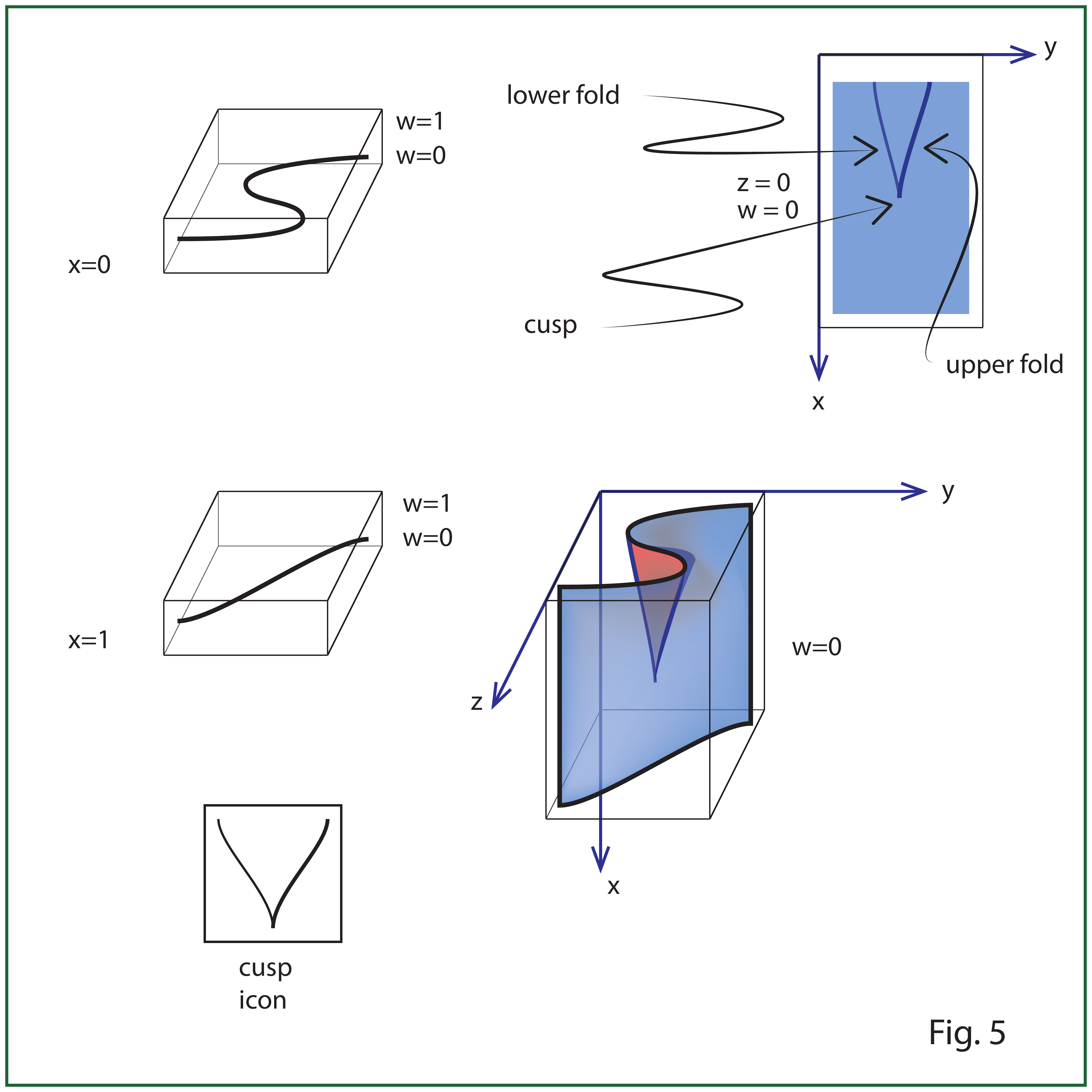}\end{wrapfigure} 
The {\it cusp} is a critical point of the fold set at which two folds converge (Fig. 5). One is closer to the observer (who resides in the positive $z$-direction), and the second fold is veiled (or hidden from the observer by means of the surface).  In the illustration, the hidden fold is drawn with a smaller line thickness, and it is indicated as slightly transparent. More generally, the cusp is represented by one of the types of vertices in the chart of the surface. It is a bivalent vertex, and the edges that converge both have labels that indicate the number of sheets in front or behind it. Thus the lower fold has a label of the form $F(n+1,m)$, and the upper fold has a label of the form $F(n, m+1)$. More details about this notation will follow shortly.

\begin{wrapfigure}[19]{l}{3.4in}\vspace{-0.32in} \includegraphics[width=3.3in]{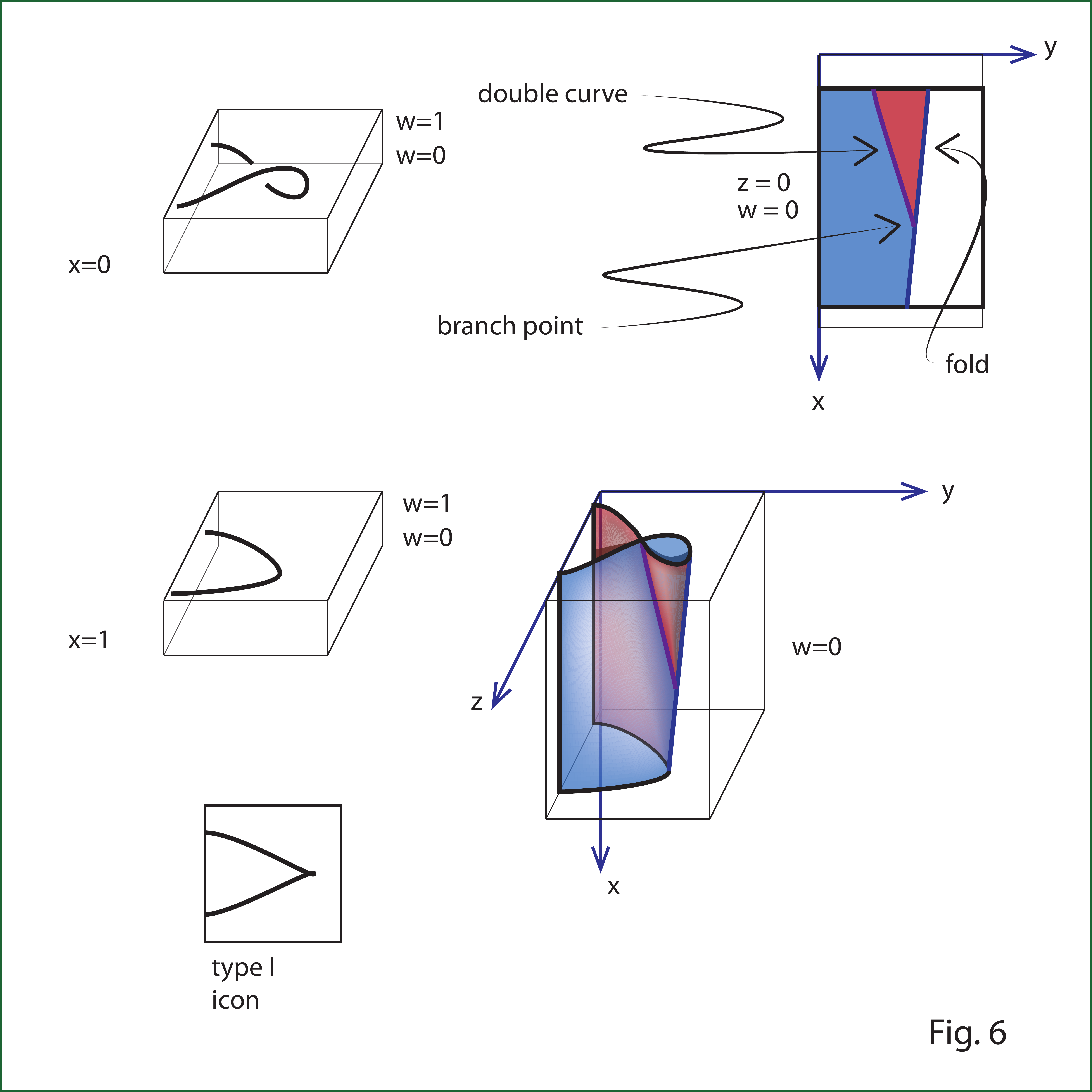}\end{wrapfigure}
 Each type of singularity has an icon that is associated to it. The icon is depicted in the lower right of the associated figure, and it is meant to suggest the nature of the singularity that is involved. 
For example, the icon associated to a birth is a solid dot that is meant to indicate the critical point. 
The icon associated to a saddle is an $x$; again this indicates the critical point but includes the level curves. The icon associated to a $\psi$-bounce is a horizontal $\psi$ that suggests the parametrization of the singularity that was described above. The icon associated to the cusp is a narrowed $v$ whose vertex is meant to suggest the cusp.

\begin{wrapfigure}[19]{r}{3.4in}\vspace{-0.32in}\includegraphics[width=3.3in]{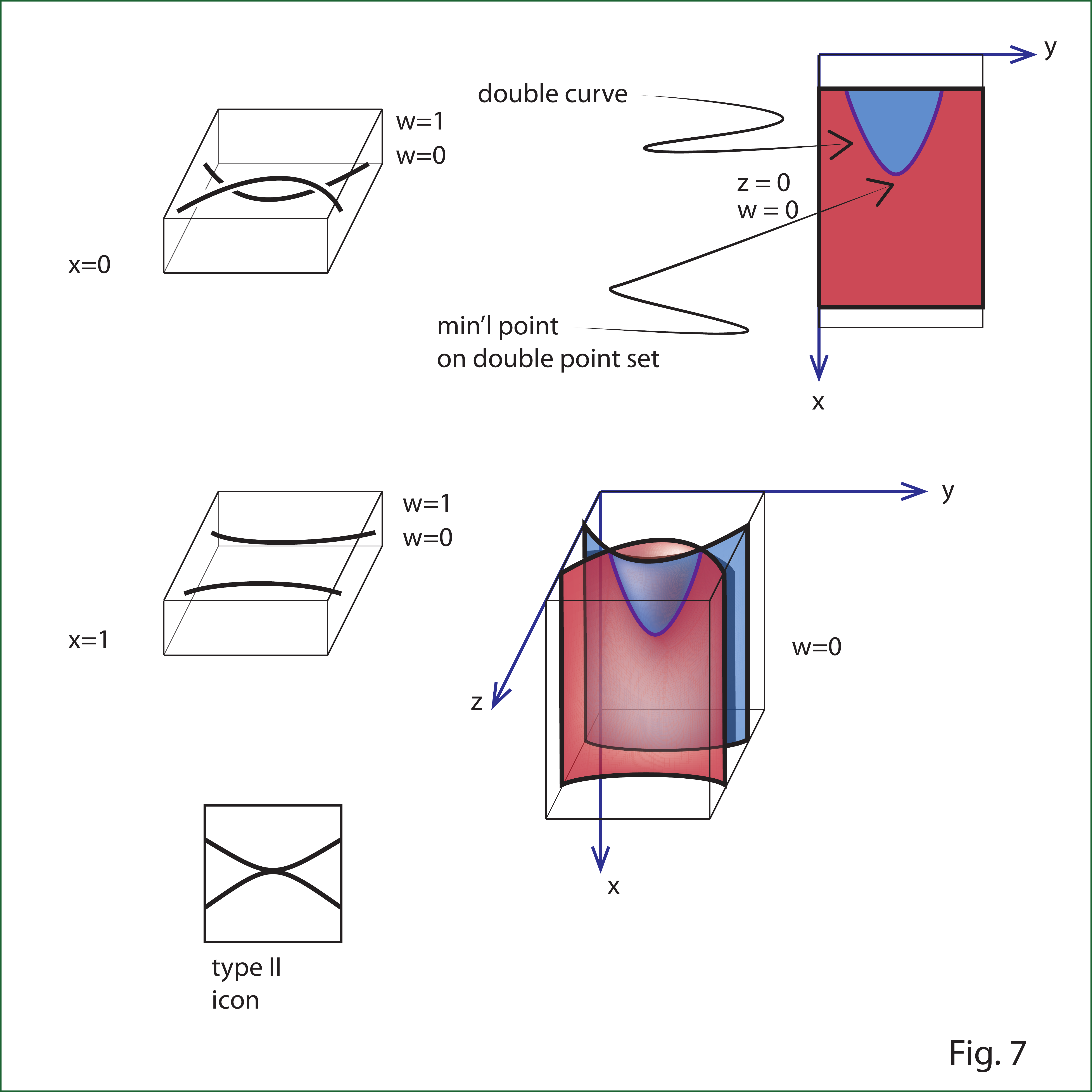}\end{wrapfigure}
 The {\it type-I move} is depicted in Fig. 6, for example, as the death of a curl on the right that has positive writhe. Of course, there are left/right and birth/death variations of this move. In addition, the writhe of the curl that is annihilated or created can be either positive or negative. The resulting branch point always occurs on a fold, but generically, it is distinct from a cusp or a saddle point{{\footnote{In a braid movie, following \cite{Kamada:book}, the branch point is simultaneously coincident with three cusps and a saddle point, but we may perturb this situation generically so that these singularities are isolated.}}}. The {\it sign of the branch point}  is said to be positive if the branch point is the birth of positive writhe or the death of negative writhe (with respect to the $x$-direction); the death of positive writhe and the birth of negative writhe yield a negative branch point.  The {\it writhe} of the crossing is determined by choosing an orientation for the curve and computing the sign of the crossing via the usual conventions. The icon for a type-I move is a sideways cusp.

To re-emphasize: The $w$-direction is indicative of over/under information. It is the bottom-to-top direction in any of the stills. By projecting to the $(w=0)$-hyperplane, the surface is generically immersed in a $3$-dimensional space (or {\it solid}). The subsequent projection to the $(x,y,0,0)$-plane is the chart of the surface.   The singular points for a surface that is mapped to a plane consists of folds that may end in cusps. In the chart, the $x$-direction points downwards, and the $y$-direction points to the right. The $y$-direction provides a height function for all of the stills, while the $x$ direction provides a height direction (down) for the surface. Thus generic critical points for the surface (births, deaths, and saddles) and critical points for the double point set (type-II moves), are with respect to the $x$-direction. We remark further that cusps always point up or down with respect to the $x$-direction, and branch points (type-I) moves always occur along folds.

 The {\it type-II move} is depicted in Fig. 7 in one of its variants. Note that both   of the crossings in the $x=0$ still could change and the knot isotopy could occur. \begin{wrapfigure}[19]{l}{3.25in}\vspace{-.32in}\includegraphics[width=3.2in]{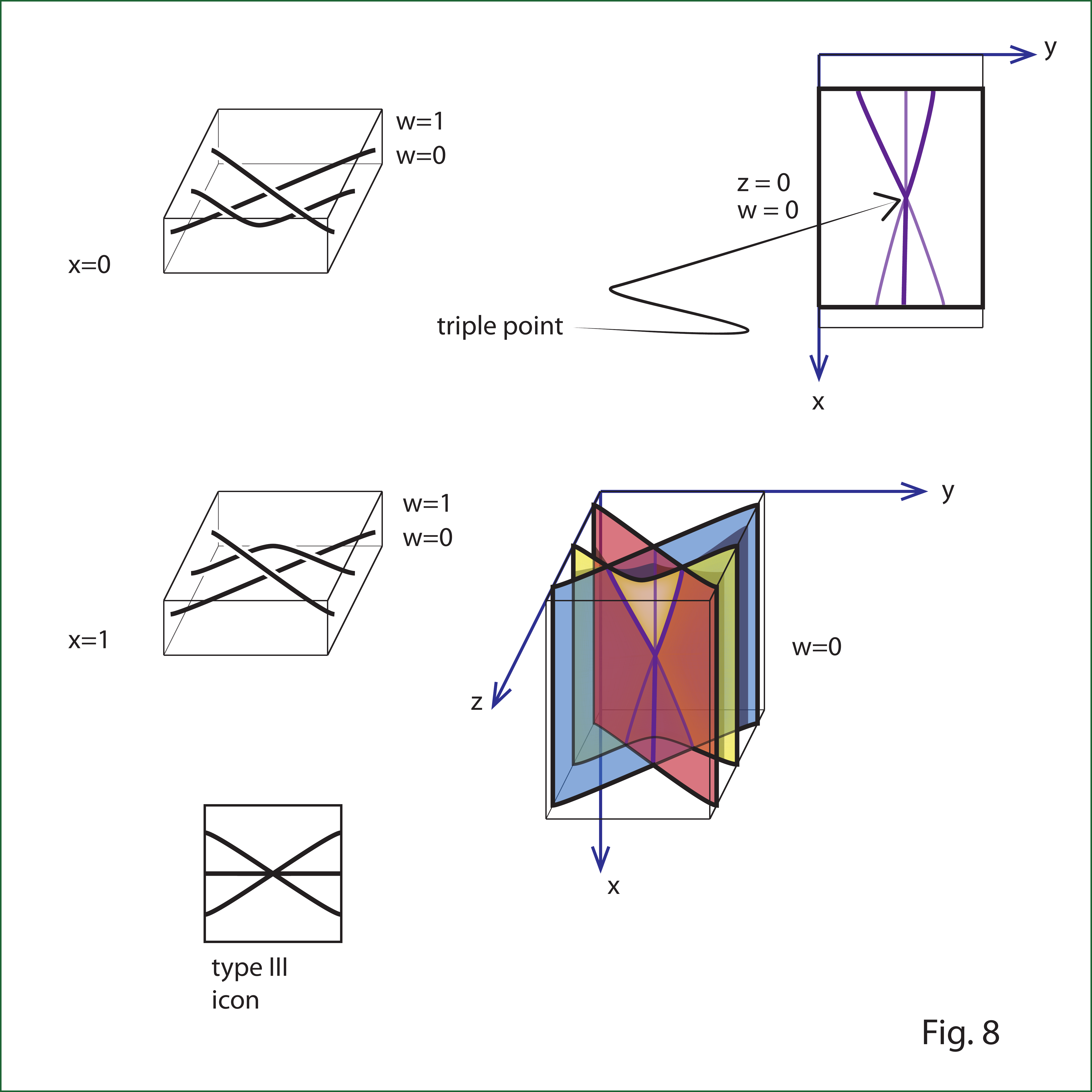}\end{wrapfigure}  Also the $x$-direction could be reversed. Thus there are four variations for the type-II move: Maximal/minimal points on the double curves reversed and crossings reversed. The type-II move represents an optimal point of the double point set. Its icon is a pair of tangent curves. The point of tangency occurs at the optimal point for the double points.

Here and in the illustrations for the type-I move the double point set is depicted in purple while the fold set is depicted in blue. Within the current illustrations, the lines of the fold and double point sets are depicted in their full thickness and opacities. In subsequent illustrations, the thicknesses and opacities will vary as the curves are veiled by other sheets of the surfaces. These drawing parameters can be reflected in labels of the corresponding arcs in the chart. For example, a double point arc  is indicated as a certain type (purple or D) of edge. A fold line is another type (blue or F). The {\it depth} of an edge is a pair of non-negative integers that indicate the number of sheets of the surface that lie in-front (towards the observer or in the increasing $z$-direction) or behind (towards the $z=0, w=0$ plane). Thus a double point arc is a labeled edge in a chart in which the label is of the form $D(m,n)$ and a fold is an edge of the chart of the form $F(m,n).$ For example, the branch point at a type-I move is a vertex in a chart that is coincident to an edge labeled $D(m,n)$ and a pair of edges labeled $F(m,n)$.

 The {\it type-III move} occurs when three lines coincide (say at $x=1/2$) in the projection to the $w=0$ hyperplane (Fig. 8).   The three sheets of surface are called the {\it top, middle,} and {\it bottom} sheets. These directions refer to the $w$-directions.   In the $(w=0)$-solid, a neighborhood of this {\it triple point} is homeomorphic to the intersection of the three coordinate planes in $3$-space. The projection of the triple point in the $(x,y)$-plane ($(x,y,0,0)$) is a valence $6$-vertex that is coincident to three edges of a type $D(m,n+1)$ and three edges of type $D(m+1,n)$.

 \begin{wrapfigure}{r}{3.5in}\includegraphics[width=3.4in]{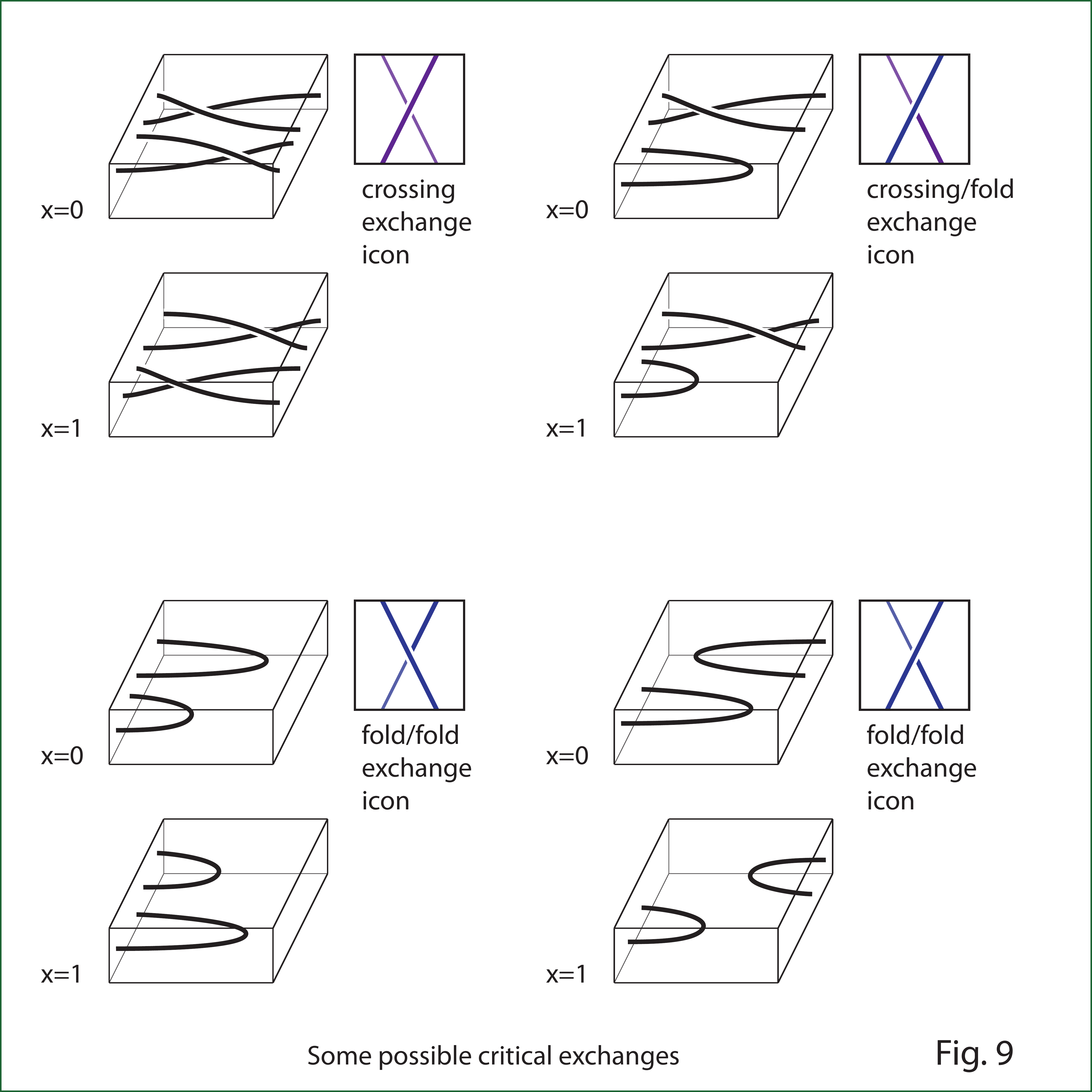} \end{wrapfigure}
Consider the  three classical crossings in the still $x=0$. There are eight possible configurations for such crossings, but only six of them will yield configurations of arcs in which there is a well-defined top, middle, and bottom arc. These six variations can occur as triple points. The icon for a type-III move is the coincidence of three line segments. 
\vspace{36pt}

{\it Critical exchanges}  occur as $4$-valent vertices in a chart when crossings or folds change heights with respect to the $y$ direction (Fig. 9). Four of the possible critical exchanges are depicted here. Either a pair of double points cross each other, a double point arc crosses a fold, or a pair of folds cross in one of two ways. The crossing point in the chart is a $4$-valent vertex. If a fold arc is incident then the labels on the other arcs will shift accordingly from top-to-bottom ($x$-direction) in the chart. 

Thus our discussion of the {\it atomic pieces} that constitute an embedded surface in $4$-dimensional space is complete. As we describe surfaces in $4$-space (or in $3$-space), we will decompose them in terms of their movie and their chart. The chart will also be used to construct the projection into $3$-space which, of course, is schematized as a drawing on the paper or screen that the reader is viewing.

Before we elaborate our examples, we have two more preliminary items to discuss. The first is the attachment of a band that is orientation reversing. The second is the Gauss-Morse code that can be used as an alternate (but incomplete) description  of the stills in a movie. Each surface that we construct can be further described via a sequence of Gauss-Morse codes. These elaborate the  {\it decker set} which is the lift of the double points, triple points, folds and their critical points to the ambient surface.

\subsection*{M\"{o}bius Bands and Saddles}
\label{peetwo}

Consider an oriented knot or link embedded in the $(x=0)$-solid face of the hypercube $[0,1]^4$. 
Each component of the link is oriented. Suppose that at $(x=1)$ a link is given that is related via a saddle point. \begin{wrapfigure}[14]{l}{2.4in}\vspace{-.25in}\includegraphics[width=2.4in]{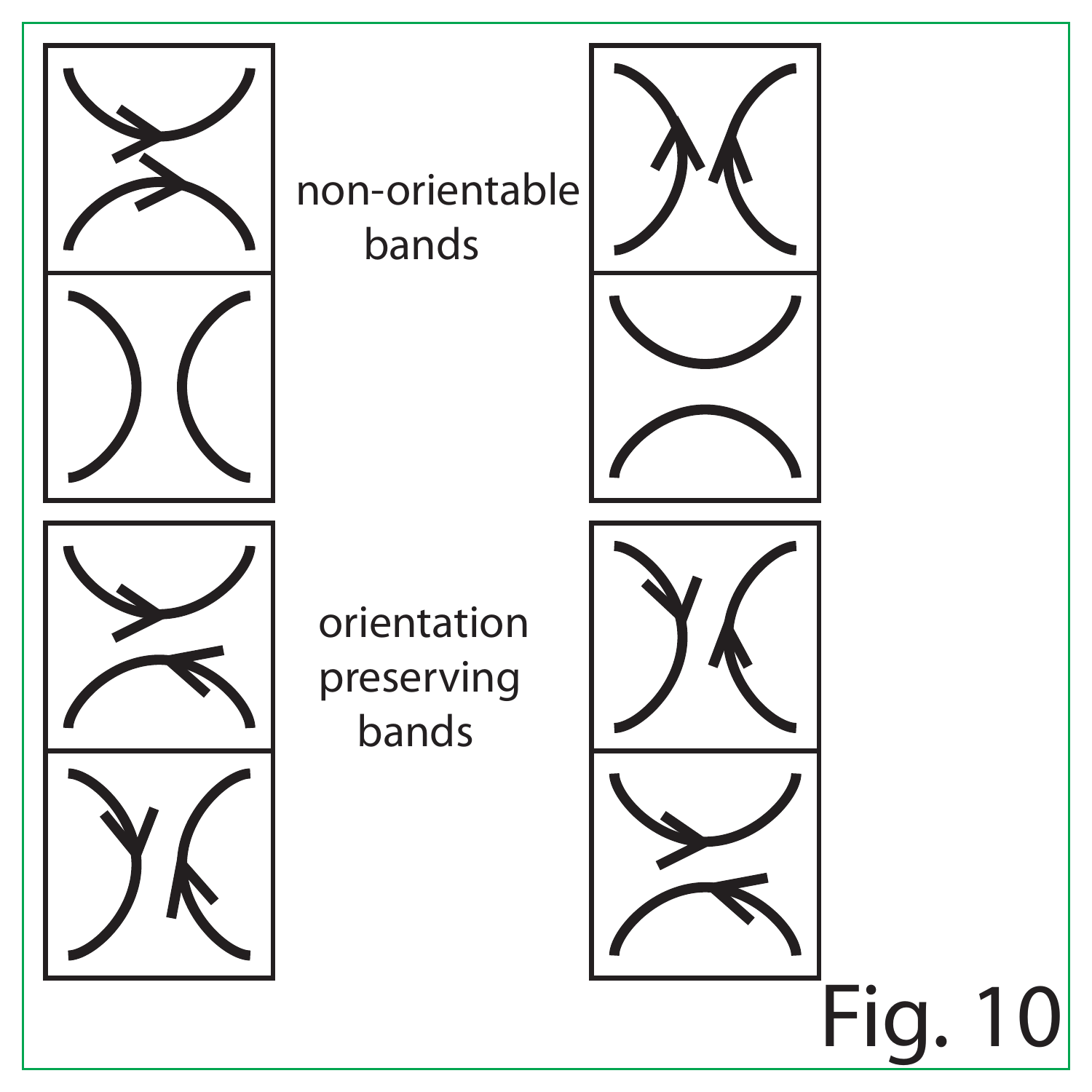}\end{wrapfigure} That is the links at $(x=0)$ and $(x=1)$ differ precisely by a single saddle point. We examine (Fig. 10) the methods to determine that the saddle is orientation reversing, or equivalently that the link cobordism in the $x$-direction is a non-orientable surface. 
When the orientations between the arcs at which a saddle is to be attached run parallel, the after the saddle attachment, any orientation on the result is inconsistent with the original orientation. This situation is indicated at the top of the figure to the left. On the other hand, the pair of arcs to be rejoined via a saddle are anti-parallel if and only if the saddle preserves the orientation. Note that the orientation problem of the saddle band is a problem of having consistent orientations throughout a movie of a surface. 

\begin{wrapfigure}[17]{l}{3.3in}\vspace{-0.32in}\includegraphics[width=3.2in]{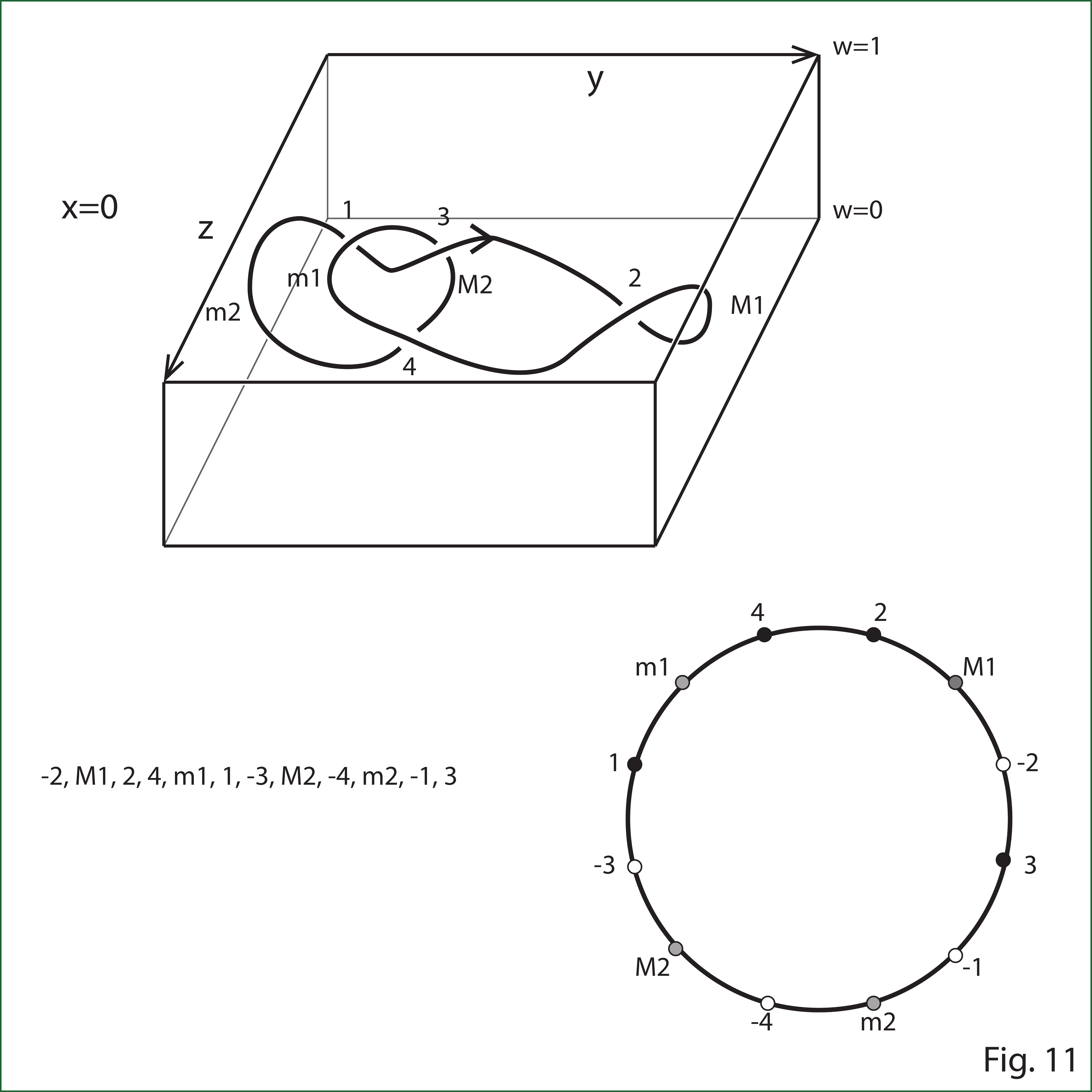}\end{wrapfigure} We remark further that if the number of components of a link are preserved under a saddle band attachment, then one component of the cobordism is non-orientable. This last condition  is sufficient, but not necessary, as the standard Klein bottle illustrates (see below).

\subsection*{Gauss-Morse Codes}\label{GM}

Consider the simple closed curve  depicted in Fig. 11 that is embedded in the $(y,z,w)$-cube $\{0\} \times [0,1] \times [0,1] \times [0,1]$. The critical points in the $y$-direction are indicated as minima: $\{m_1,m_2\}$ and Maxima: $\{M_1,M_2\}$.  
The crossing points are labeled in such a way that points of positive writhe are labeled by even integers and those of negative writhe are labeled by odd integers. \begin{wrapfigure}[20]{r}{3.5in}\vspace{-0.32in}\includegraphics[width=3.4in]{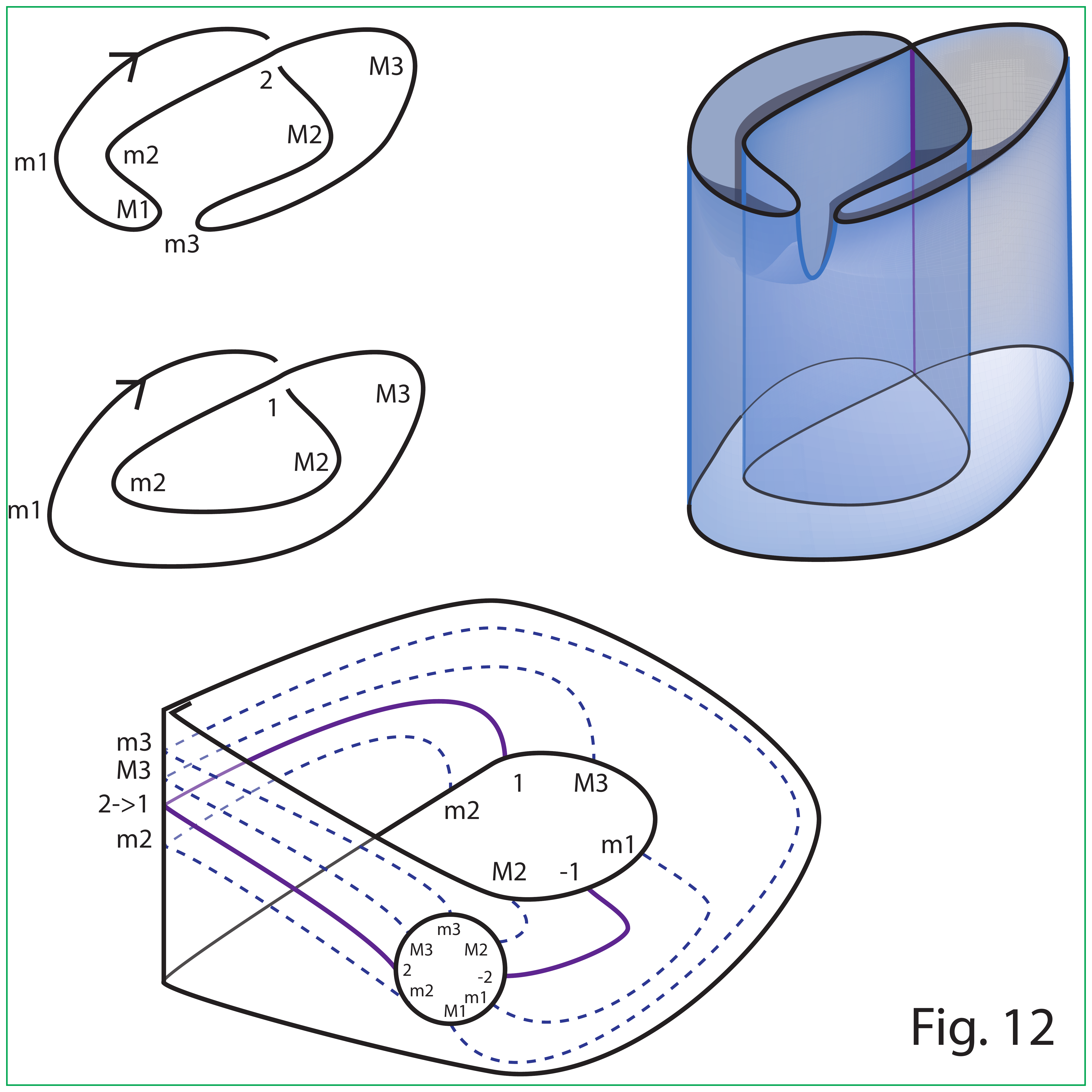}\end{wrapfigure} The base point and orientation of the curve are indicated by the small arrow at the top of the curve and to the right of the crossing point $3$. For any generic curve in $3$-space such a labeling is possible. The {\it Gauss-Morse} code of such a curve is obtained by reading around the curve following the orientation and recording each optimal point and crossing point. Under-crossings receive negative integers, and the corresponding over-crossings receive positive integer labels. The absolute value of these labels coincide for each crossing. Writhe is encoded by parity as the example suggests. The Gauss-Morse code is indicated as a sequence of letters in the labels for the crossings and optimal points, but in fact it is a cyclic word. Thus the circle at the bottom right of the illustration indicates the Gauss-Morse code on a circle. The under-crossing points are indicated as open dots on the curve while the over-crossing points are indicated as solid dots and the optimal points are grey.

Each of the critical events: birth, death, cusp, $\psi$-bounce, type-I, type-II, and type-III, have a well understood effect upon the Gauss-Morse codes of the underlying curves. Thus associated to a movie of an embedded surface, we can construct a sequence of Gauss-Morse codes associated to each of its stills. These are interconnected by crossings and critical points of the fold set and the double decker set in the ambient surface. As the surfaces are illustrated, often their double decker sets will also be indicated.

 The ideas of the last two subsections are illustrated in Fig. 12.  A pair of stills that differ by means of an orientation reversing saddle are indicated on the left. The orientation on the bottom curve cannot be induced from the orientation on the top curve. The Gauss-Morse code on the top curve is indicated on the circle that is in the interior of the M\"{o}bius band at the bottom of the figure. The Gauss-Morse code for the bottom curve is indicated on the other boundary component of the punctured M\"{o}bius band. The fold arcs for the interpolating surface are indicated as dotted arcs on punctured M\"{o}bius band, and the double decker arcs are indicated by solid (purple) lines. A feature of the movies of non-orientable surfaces is that the writhe of crossings may not be consistent from still to still particularly when an orientation reversing handle is attached.  \begin{wrapfigure}[12]{l}{2.4in}\vspace{-.32in}\includegraphics[width=2.3in]{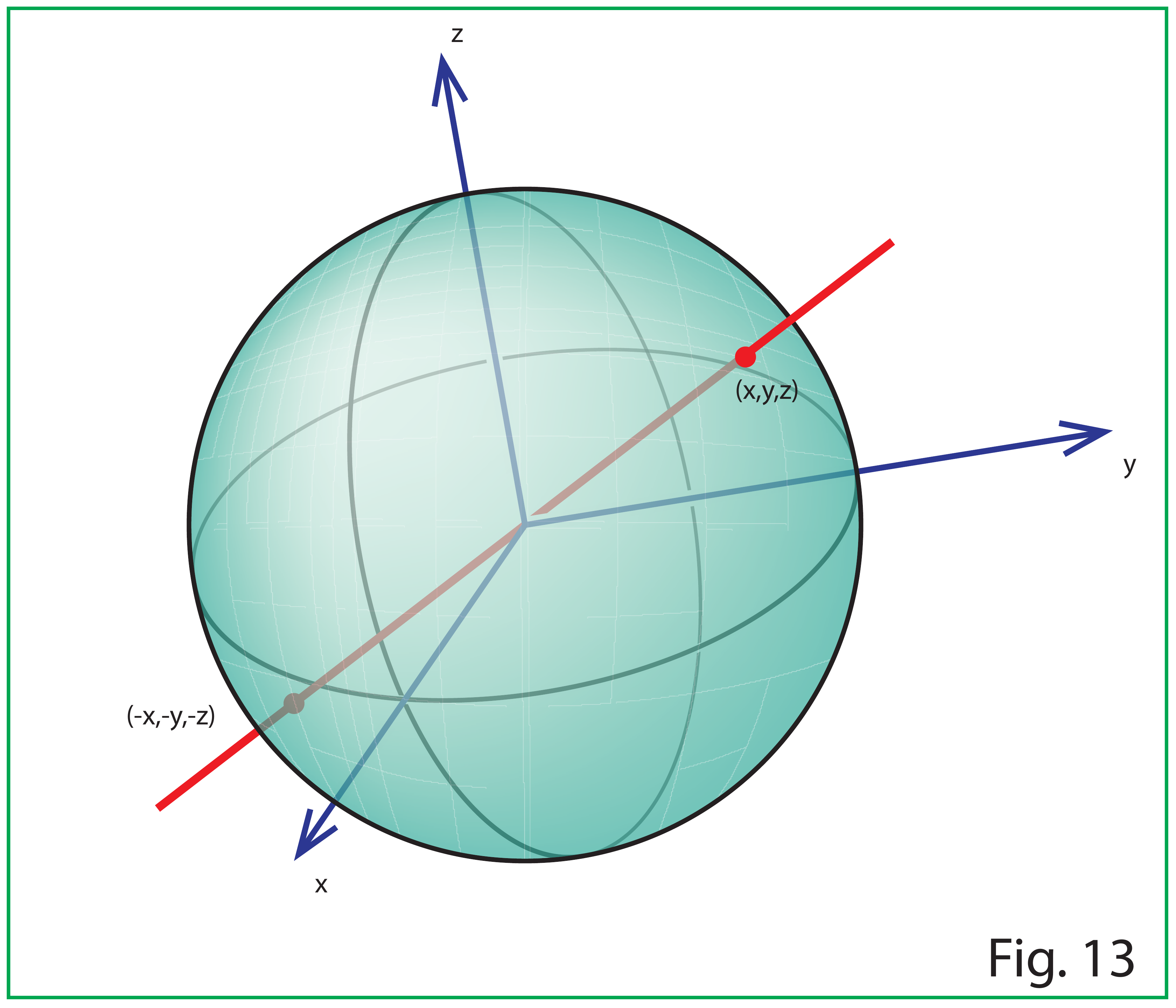}\end{wrapfigure} Thus the double point labeled $2$ at the top of the movie changes to a $1$ after the handle has been attached. Additional information can be encoded in the double decker set. For example the double decker curve associated to the lower deck that is indicated by $-1$ and $-2$ could be weighted differently than the upper decker curve that is indicated by a $1$ and $2$. Here that indication is suppressed since the figure is relatively simple, and the sign on the decker curve seems sufficient as a $w$-height indicator. 
 
 \section{The Cross-cap embeddings}

   \begin{wrapfigure}[12]{r}{2.5in}\vspace{-0.42in}\includegraphics[width=2.4in]{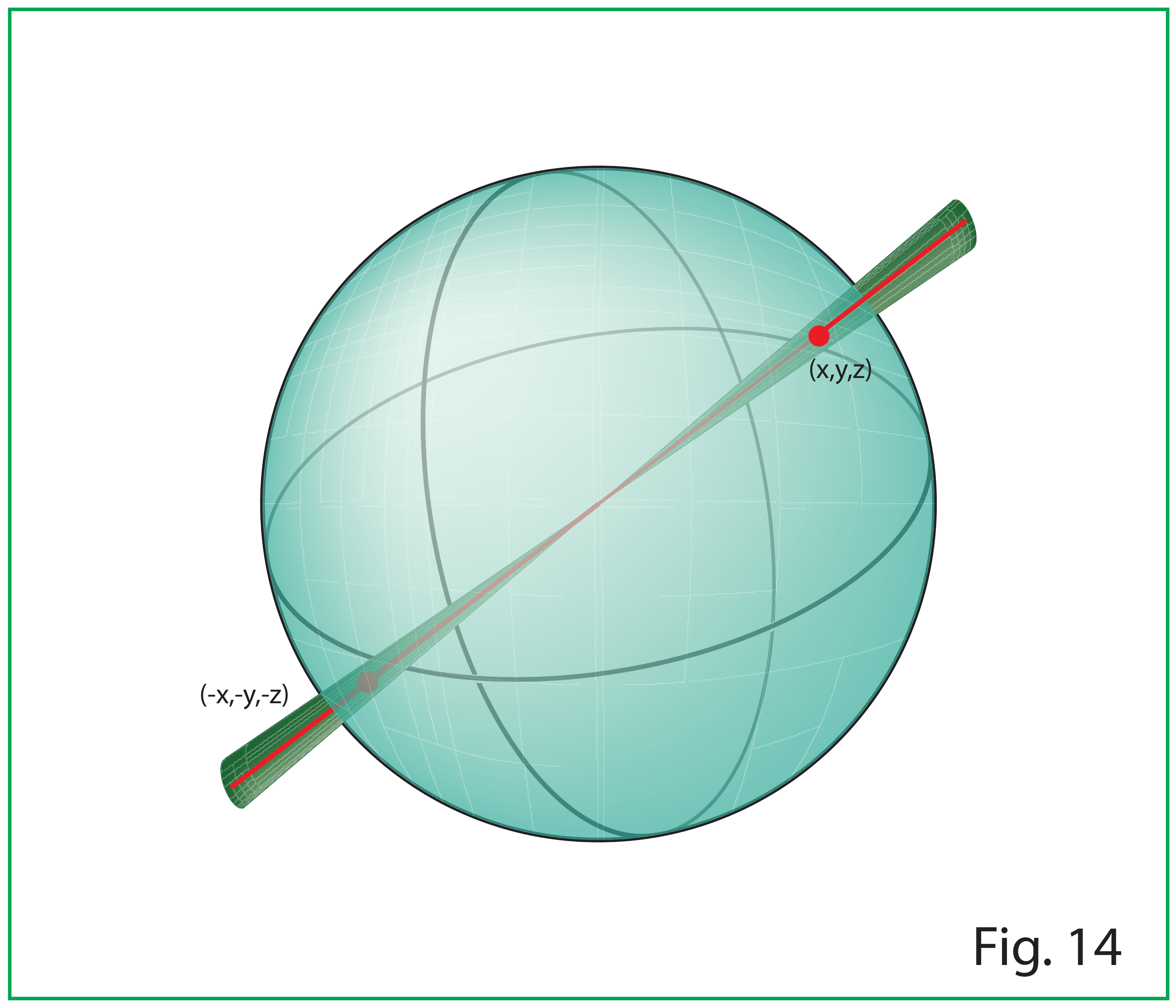}\end{wrapfigure}
 The projective plane is the most simple of the closed non-orientable surfaces. A {\it closed surface} is a compact surface that has no boundary. So each point in the surface has a neighborhood that is homeomorphic to an open subset of the plane $\R^2.$ Of course, our surfaces are Hausdorff and second countable. In the case of surfaces, the condition of non-orientability is equivalent to the existence of an embedded open subset that is homeomorphic to an open M\"{o}bius band. This prototypical non-orientable surface is the quotient space $[0,2\pi] \times (-1,1)/\sim$ where $(0,t)\sim (2\pi, -t).$ The projective plane $\R P^2$ is constructed in one of four equivalent ways.

 First, we may consider the $2$-sphere $S^2$ as the set $$S^2= \{ (x,y,z) \in \R^3: x^2+y^2+z^2=1 \}.$$ Define an equivalence relation $\sim$ on $S^2$ by $(x,y,z)\sim (-x,-y,-z)$; thus antipodal points on the sphere are identified. The {\it projective plane} is the quotient space $S^2/\sim$ (Fig. 13). Note that there is a representative point of each equivalence class on the closed upper hemi-sphere $\{(x,y,z): \ 0\le z \}$. Furthermore, diametrically opposite points on the equator ($\{(x,y,0): x^2 +y^2=1 \}$ are identified in the quotient.

\begin{wrapfigure}[12]{l}{2.4in}\includegraphics[width=2.3in]{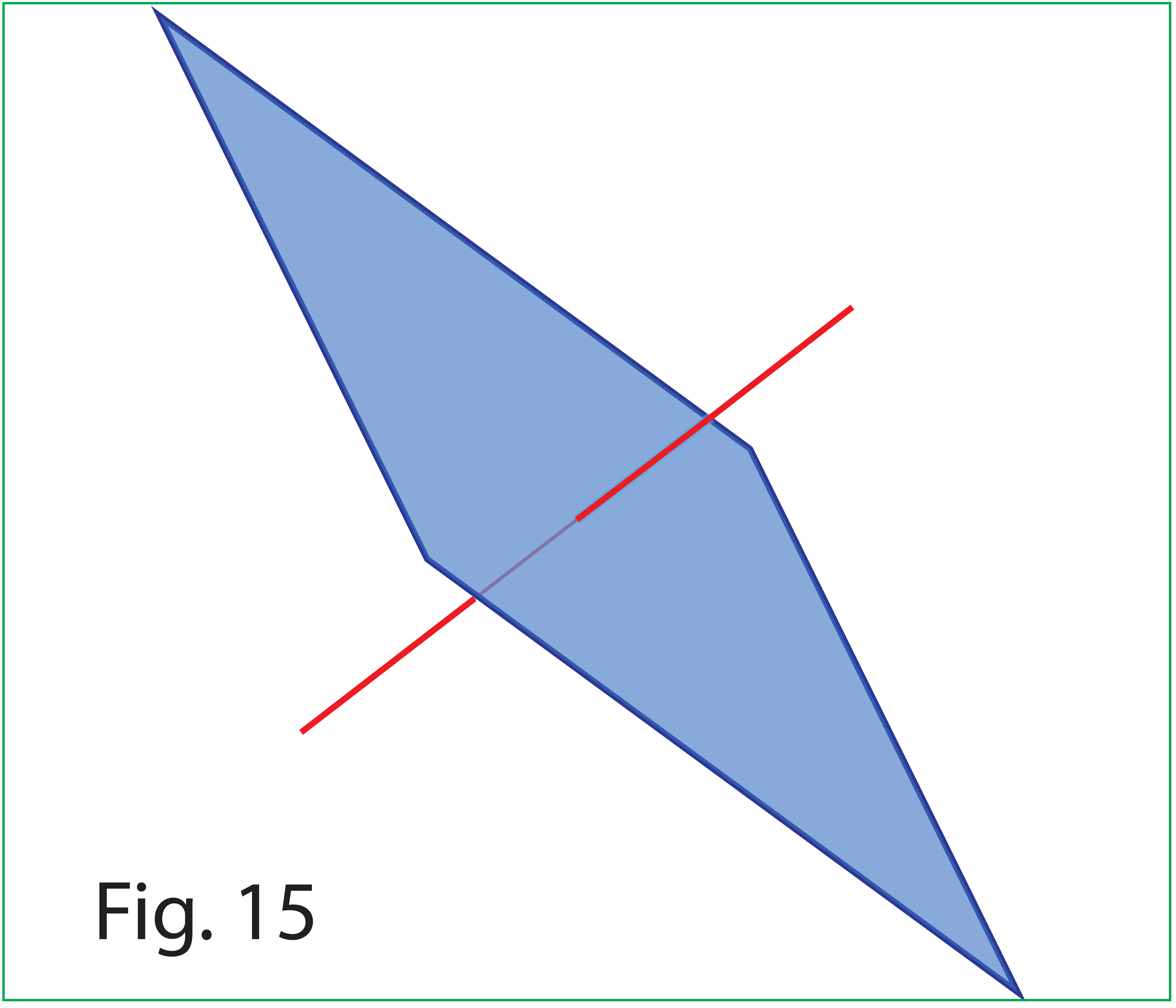}\end{wrapfigure}
 Alternatively, we may consider the projective plane $\R P^2$ as   the set of all lines through the origin in $\R^3$ . Two such lines are considered to be close if they lie within the same small double cone whose vertex is at the origin (Fig. 14). 
Any such line is determined by its direction vector $\vec{v}=(A,B,C)$. Say $\vec{\ell}(t)= t(A,B,C).$ Thus such lines are close if the angle, $\theta = \arccos{ \left( |\langle\vec{v_1}, \vec{v_2}\rangle /( || \vec{v_1}|| \ || \vec{v_2} ||) | \right)}$, between their direction vectors is small. By normalizing the direction vector $\vec{v}\mapsto \vec{u}= \vec{v}/||\vec{v}||$, we obtain a unique point $\vec{u}$ on the unit sphere. But the direction vector $\vec{v}$ determines a {\it directed} line. Thus both $\vec{v}$ and $-\vec{v}$ (and hence $\vec{u}$ and  $-\vec{u}$) determine the same line in space. Thus the space of all lines is equivalent to quotient space of the sphere under the antipodal map.

Similarly, the vector $\vec{v}$ determines a plane through the origin $\{ (x,y,z): \ Ax+By+Cz=0 \}$. So the projective plane can topologized as the set of planes containing the origin in $3$-space where two planes are close if their normal vectors make a small angle (Fig. 15). 

\begin{wrapfigure}[13]{r}{2.4in}\vspace{-0.48in}\includegraphics[width=2.3in]{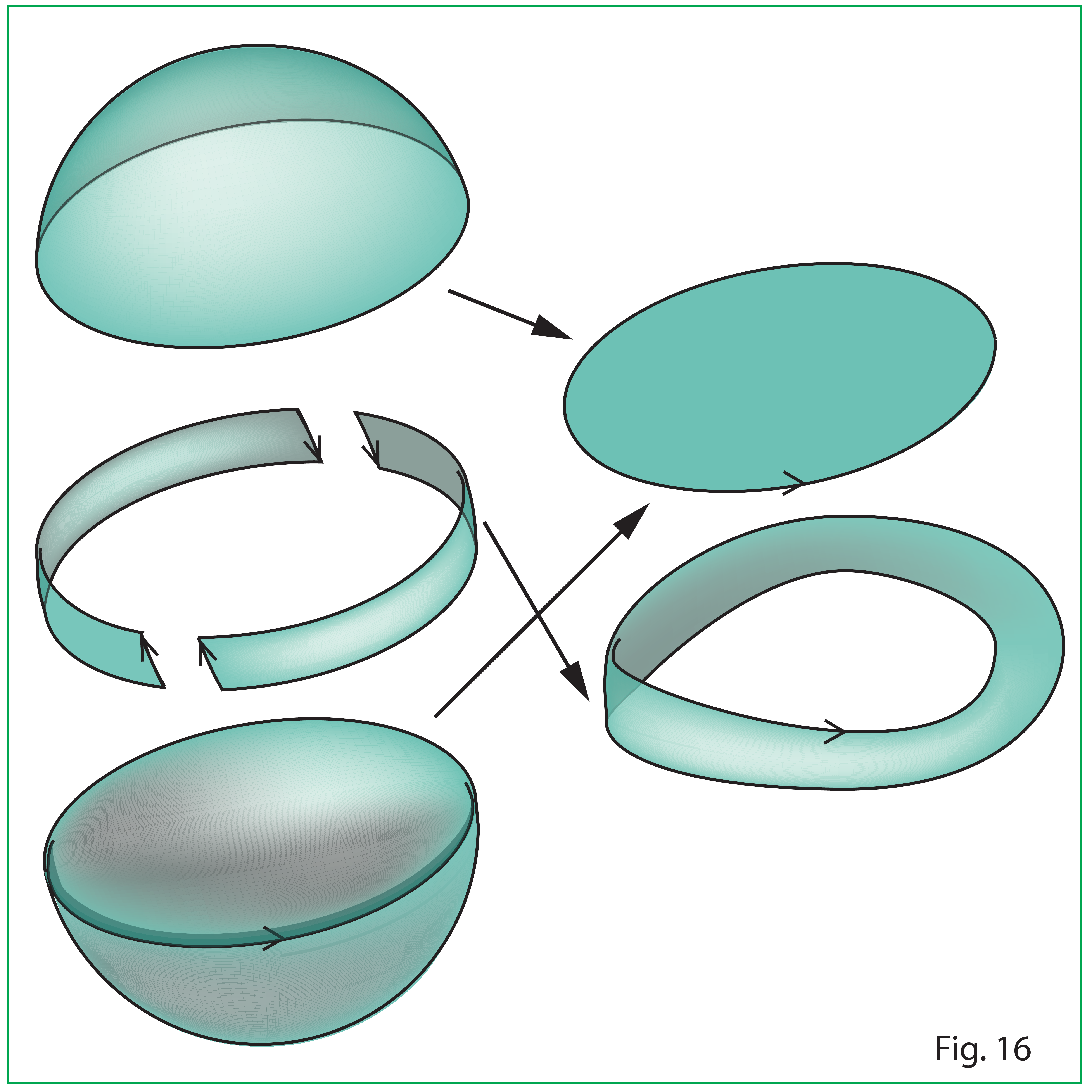}\end{wrapfigure}
It is also convenient to consider the $2$-sphere, $S^2$, as decomposed into four regions that we describe roughly as the northern polar region, the southern polar region, the eastern tropics, and the western tropics. In this loose description, let us say that the tropical latitudes are at $22^\circ$ north or south of the equator, that the  northern and southern polar regions have latitudes greater than or equal to $22^\circ$, and that the western and eastern regions are delineated by means of the prime meridian and the international dateline. The quotient of the $2$-sphere under the antipodal map, consists of two pieces:  a polar disk that is the quotient of the northern and southern polar regions, and a tropical M\"{o}bius band that is obtained from either tropical strip by identifying the northern latitudes on the prime meridian to the corresponding southern latitudes on the international dateline. In this model, the boundary of the M\"{o}bius band is identified to the boundary of the polar disk.

\begin{center}
\includegraphics[width=5in]{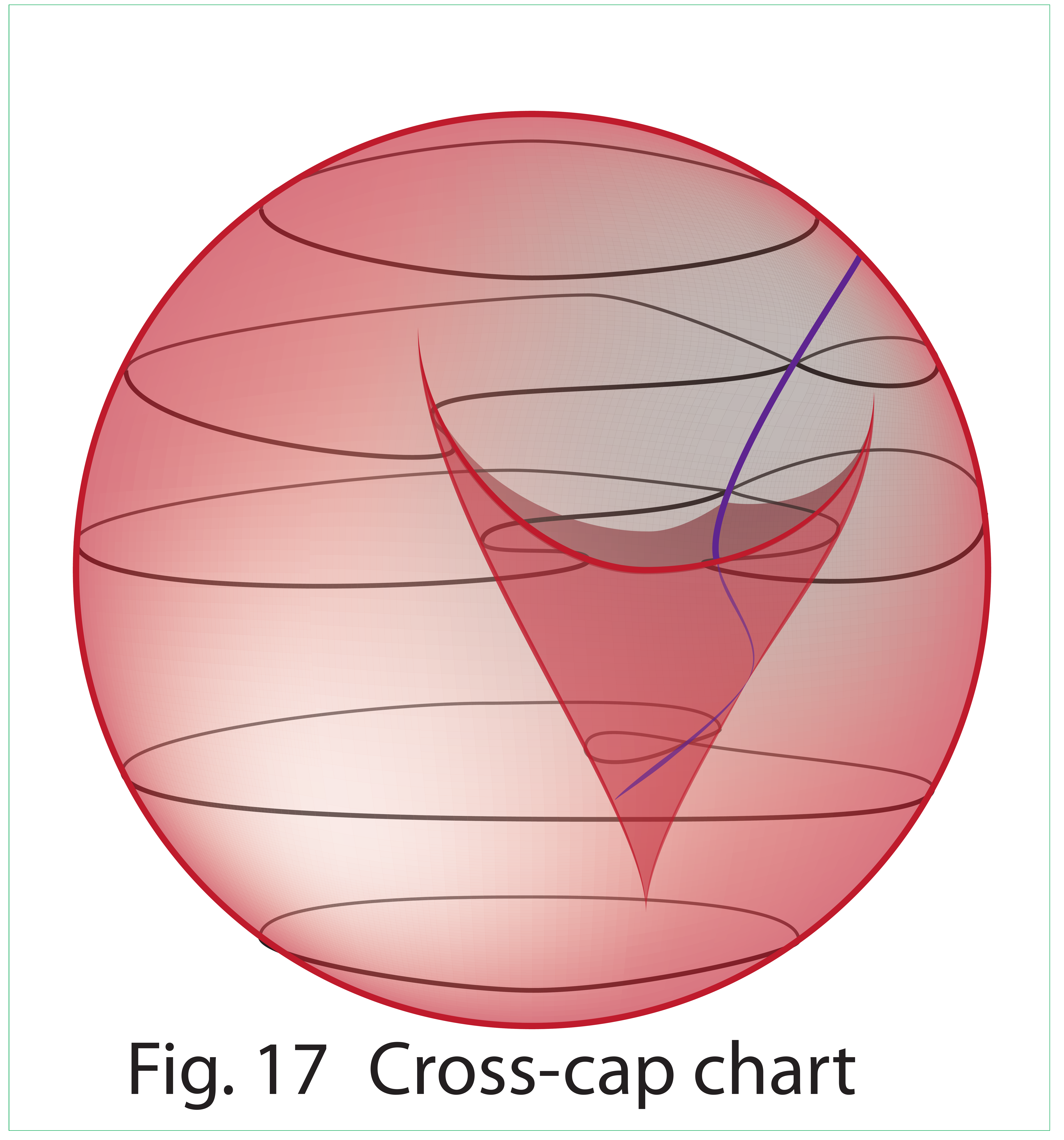} \\

\includegraphics[width=4.75in]{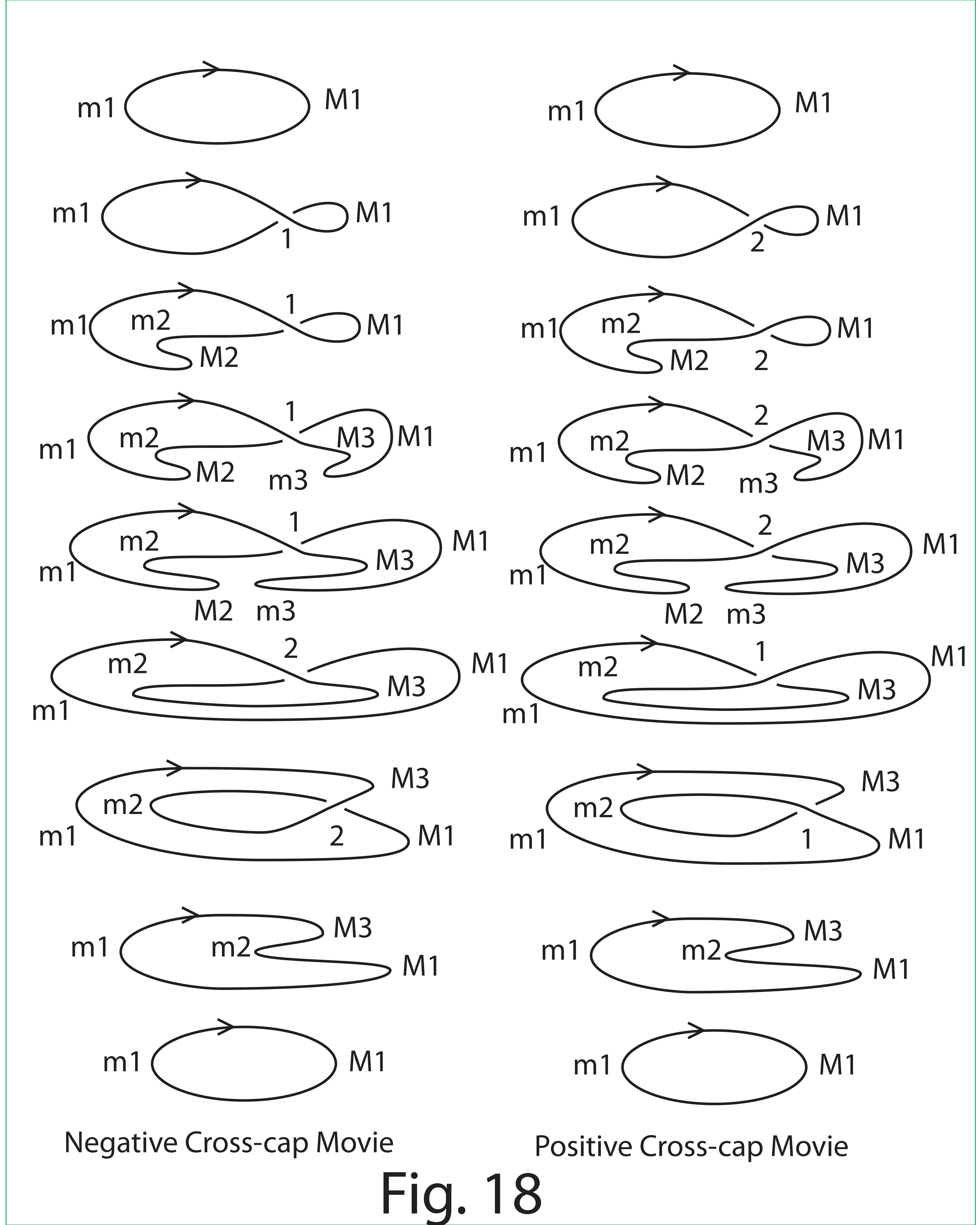} \\ 

\includegraphics[width=4.75in]{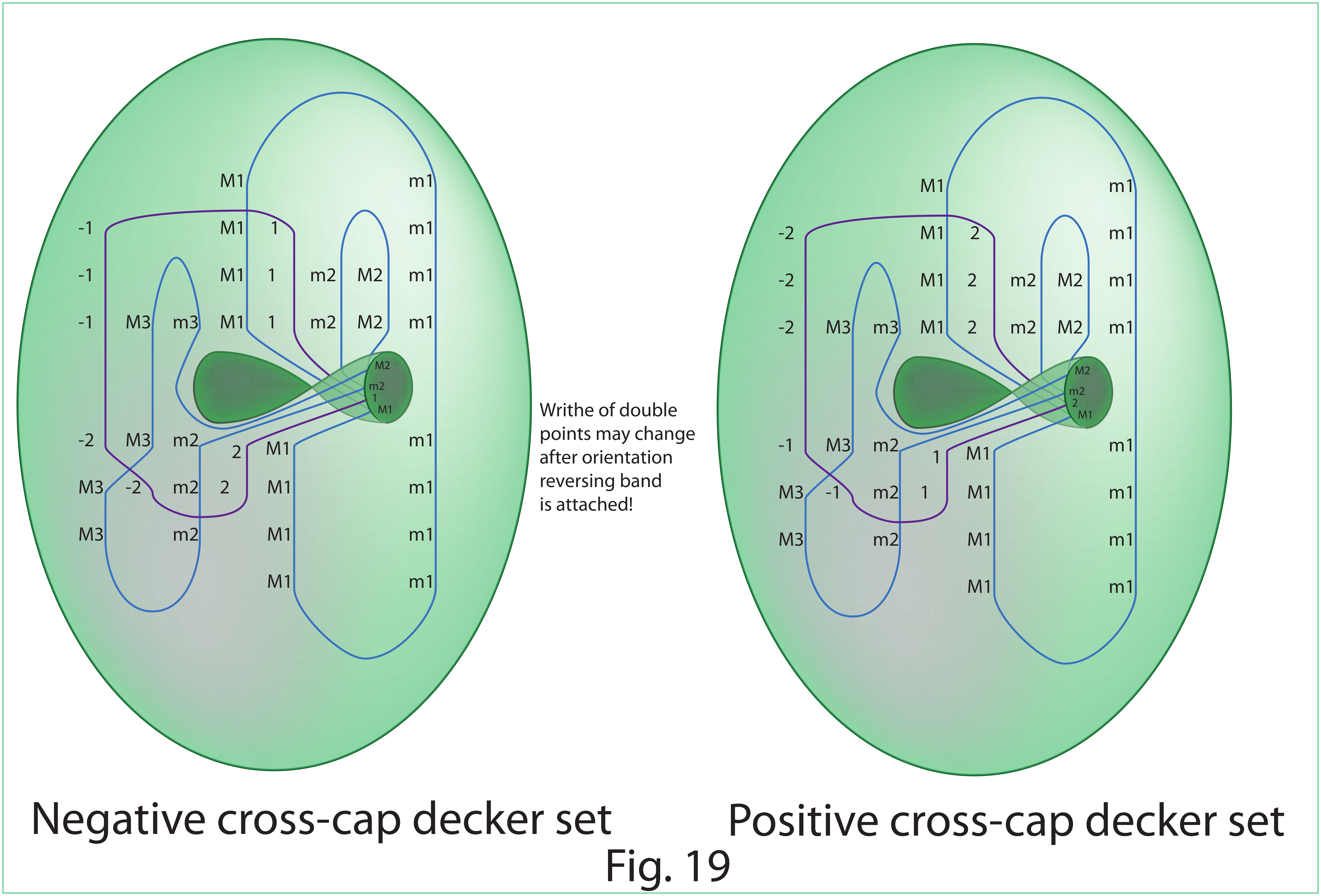}   \\

\end{center}

Up to ambient isotopy, there are two {\it standard cross-cap} embeddings into $4$-space of the projective plane. We denote these by $CC(\R P^2(\pm)).$ These are depicted here as movies (Fig. 18), in terms of their projections (Fig. 17), and with their decker sets (Fig. 19). The difference between them is contained within the nature of their branch points. The negative cross-cap $CC(\R P^2(-))$ has two negative branch points while the positive cross-cap $CC(\R P^2(+))$ has two positive branch points.

 In our movies of the cross-caps, we have assigned a base point and direction for each curve in the stills of the movie in order to reconstruct the decker set from a sequence of Gauss Morse codes. The sequence of base points can be thought of the left edge of the diagram of the decker set. One reads the Gauss-Morse codes and stacks them vertically when writing them. The book \cite{R2B} discusses such conventions extensively in the context of orientable surfaces.

 Massey's theorem (Whitney's conjecture) that gives limits upon the number and types of branch points of projections of embedded non\-or\-ientable surfaces in $4$-space. The branch points contribute locally to the normal Euler class of the embeddings. This characteristic class is a cohomology class (with local coefficients) that detects the non-triviality of the normal bundle for the embedding. A geometric description is obtained as follows. 

%


%

 \subsection{Massey's Theorem}\label{MTWC}

Let an embedded surface $F$ in $4$-space be given. Consider the intersection between the embedding and a generic push-off in the normal direction. If the normal bundle is trivial ({\it i.e.} homeomorphic to $F\times \R^2$), then there is a push-off that does not intersect the original surface. If not, we examine the neighborhoods of the intersection points. Generically, such an intersection has a neighborhood in which the intersection is homeomorphic to a pair of coordinate disks $D^2(x,y)= \{(x,y,0,0): x^2+y^2 \le 1\}$ and $D^2(z,w)= \{(0,0,z,w): z^2+w^2\le 1\}$. Even though we are thinking about non-orientable surfaces, we can give the neighborhood of the intersection point in $F$ a local orientation and give the push-off an induced orientation (Fig. 20). \begin{wrapfigure}{r}{2.1in}\vspace{-0.32in}\includegraphics[width=2in]{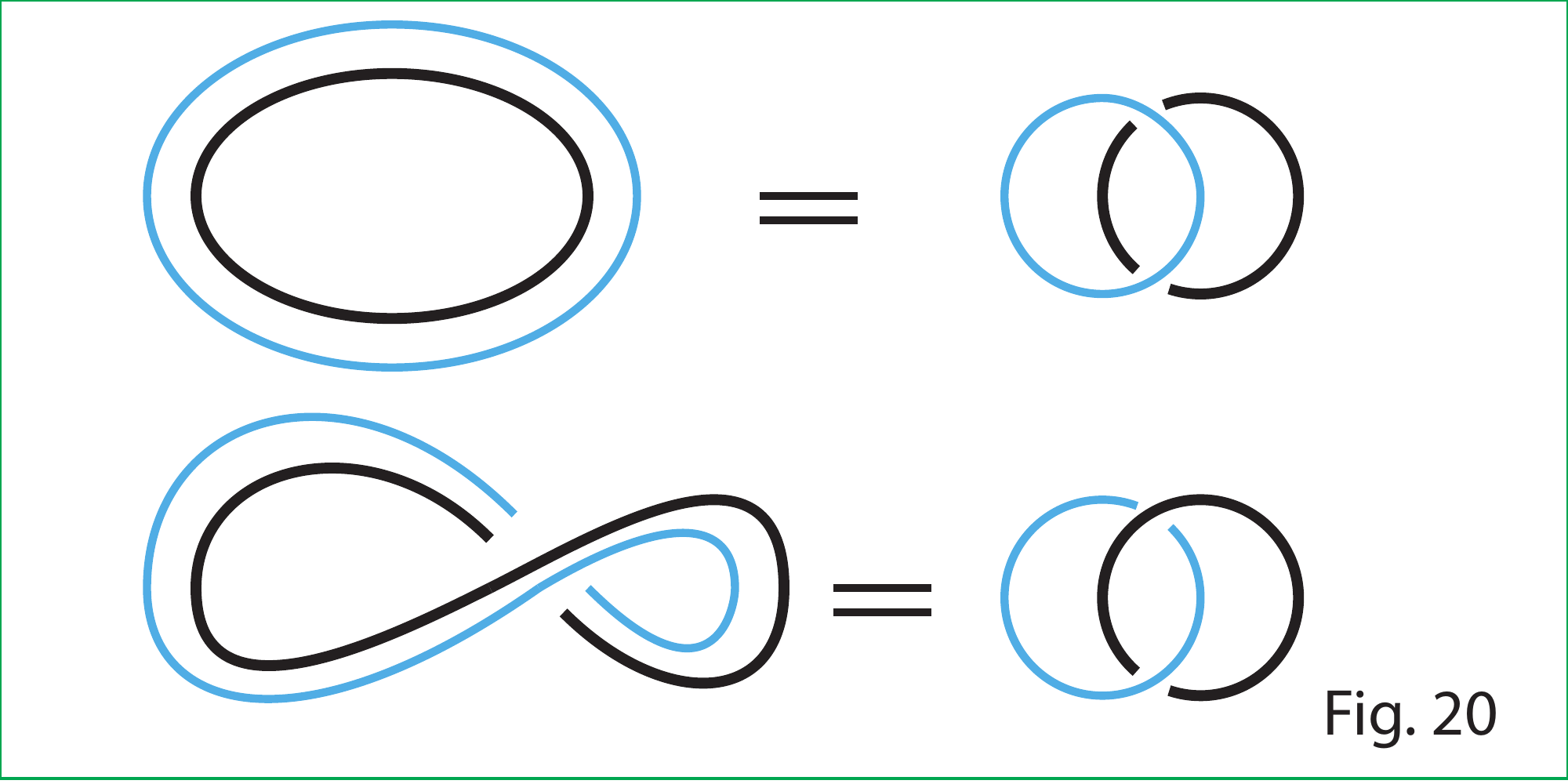}\end{wrapfigure}  
The {\it sign of the intersection point} is defined by comparing  the fixed (right-handed) orientation of $4$-space with the induced orientation of the intersection point: an oriented frame is obtained by the chosen frame in the $(x,y)$-disk followed by an oriented frame in the $(z,w)$-disk which is contained in the push-off surface.

{\it This sign is independent of the initial choice of orientation on the $(x,y)$-disk.} Analogously, the sign of a classical crossing on a knot diagram is independent of the choice of the orientation of the knot.

Looking at the movie representation of a cross-cap nearby a branch point, we can see the push-off and compute intersection numbers. In Fig. 20, immediately before the branch point, the push-off and the curve on the surface $F$ have trivial linking number. After the branch point occurs, the linking number is $+1$. Thus the push-off disk must intersect the original disk $+1$ time.

\noindent
\includegraphics[width=5in]{connect}

Massey's theorem on normal Euler classes was originally conjectured by Whitney. One imagines that Whitney made the conjecture upon forming the connected sum of several copies of the cross-cap (Fig. 21). The {\it non-orientable genus} of a surface is the number of projective plane connect-summands. The connected sum of $n$ cross-caps has genus $n$ and Euler characteristic{\footnote {not to be confused with the normal Euler class}} $\chi(\#_{i=1}^n CC(\R P^2))= 2-n.$ If all the cross-caps are positive, the normal Euler class is $2n.$ If they are all negative, then the normal Euler class is $-2n.$ Otherwise, switching any connect-summand  from a positive cross-cap to a negative cross-cap changes the normal Euler class by $4$. Thus Whitney conjectured that the normal Euler class of a surface of non-orientable genus $n$ is found among the set of integers $\{-2n,-2n+4,\ldots,2n-4,2n\}$. This theorem was originally proven by Massey~\cite{Massey}. Kamada gave an elementary proof much later \cite{Kamada}. 

In addition to the cross-cap embeddings, we have the Roman surface, Boy's surface, and the so-called girl's surface. Each is a general position map of the projective plane into $3$-space. There is an embedding of the projective plane into $4$-space that lifts to the Roman surface, and there are immersions of the projective plane that lift to each of the Boy's and girl's surface. In this section, we will work with these lifted surfaces when indicating the movies and decker-sets. Keep in mind, {\it the associated general position maps can be obtained by removing the crossing information in the movies.} 

\subsection{The Roman  surface}\label{Roman}

\begin{center}

\includegraphics[width=5in]{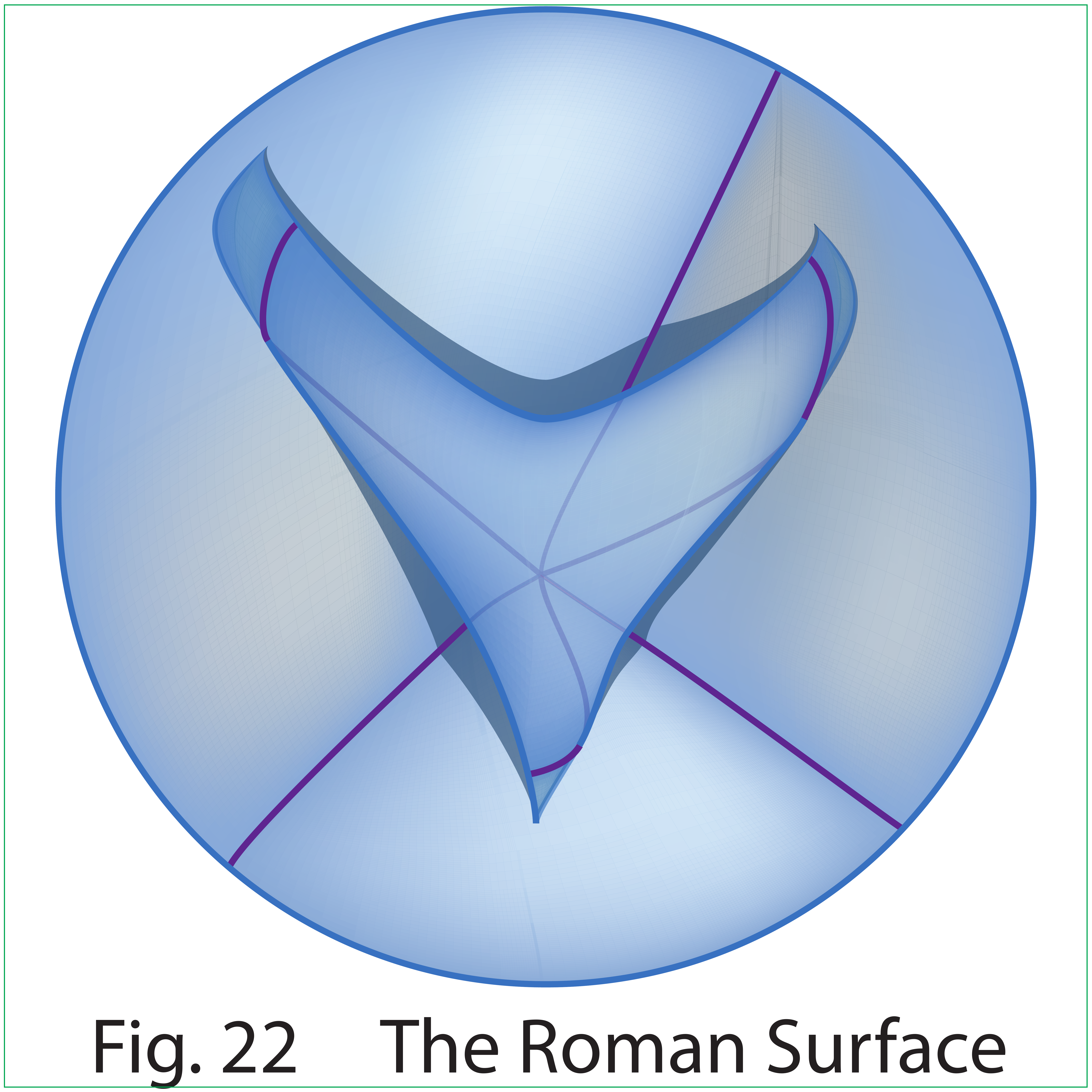} \\ 
\includegraphics[width=5in]{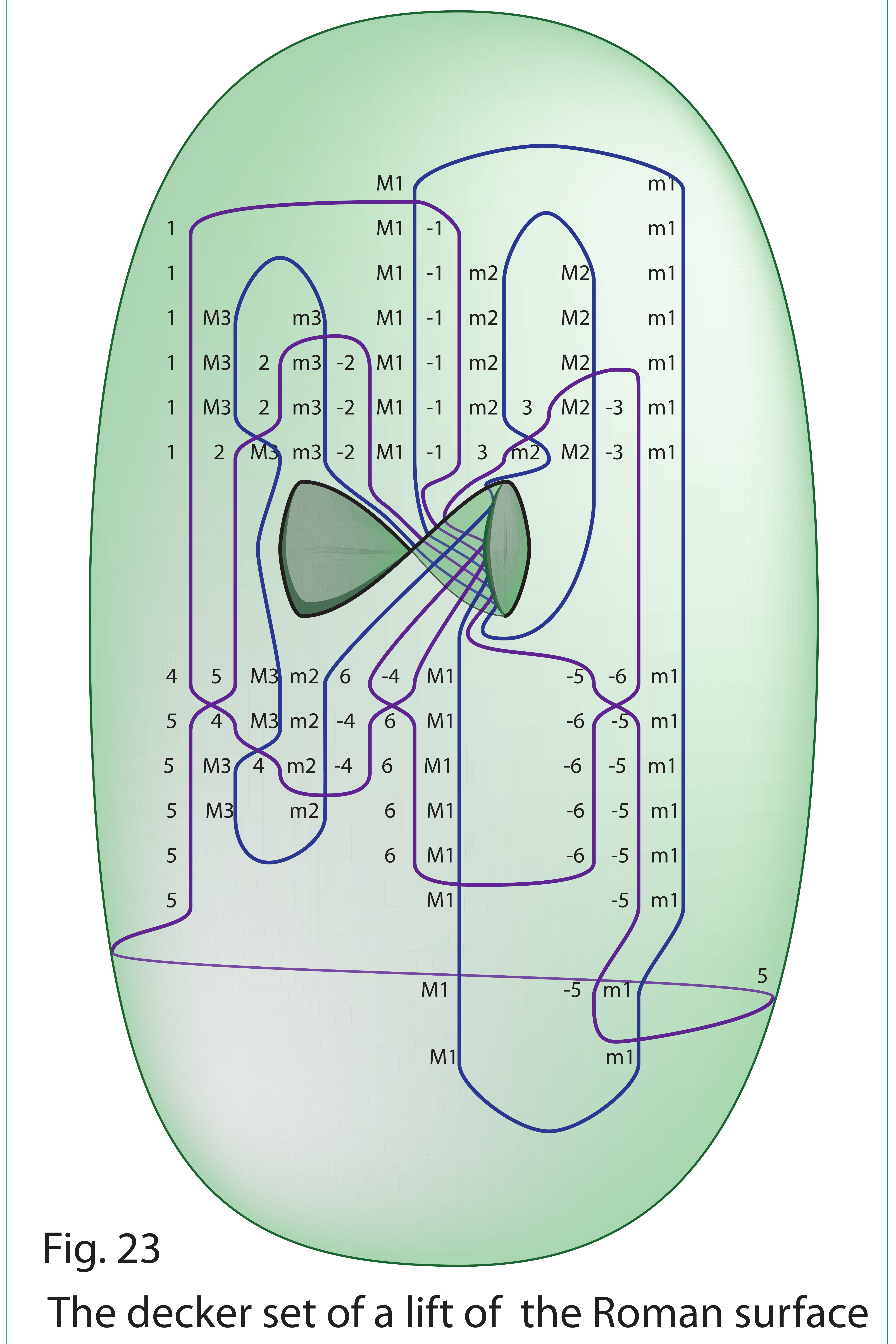}
\end{center}

\begin{wrapfigure}[17]{l}{3.02in}
\includegraphics[width=3in]{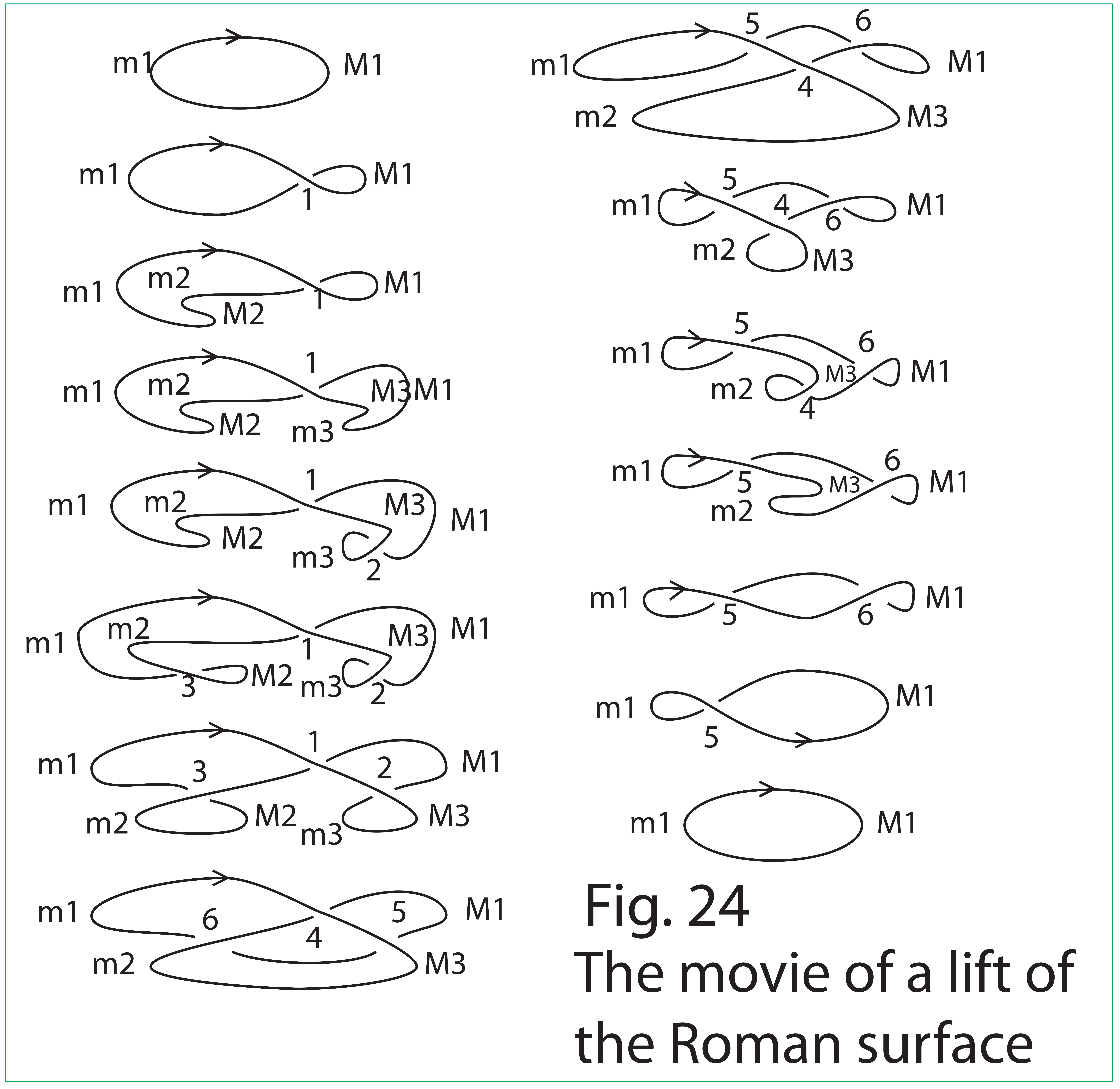}   \\ 
\end{wrapfigure}
The Roman surface  (Fig. 22) has six branch points and a single triple point. It, too, can be embedded in $4$-space as the accompanying movie indicates (Fig 24 see also Fig. 22).  Its $3$-dimensional image can be constructed by taking three coordinate disks  $\{(x,y,0): x^2 +y^2 \le 1 \}$, $\{(x,0,z); x^2+z^2 \le 1\}$, and $\{(0,y,z): y^2+z^2\le 1\}$ and gluing to these the four spherical triangles in which the product of the coordinates is positive. It has 6 branch points and one triple point. This embedding is ambiently isotopic to one of the standard cross-caps. We leave the reader two exercises: (1) which cross-cap is this embedding ambiently isotopic to? (2) construct the isotopy.

\subsection{The Boy's surface}\label{Boy}
\begin{center}
\includegraphics[width=2in]{BoyChart}\end{center} 
\begin{center}
\includegraphics[width=3.5in]{Boydeck}
\end{center}

\begin{wrapfigure}[19]{l}{3.02in}\vspace{-.5in}
\includegraphics[width=3in]{Boymovie}   
\end{wrapfigure}
The Boy's surface~\cite{HilbertCV} is an immersion of the projective plane into $3$-di\-men\-sional space (Fig. 25). An immersion has no branch points. So, by Massey's Theorem, this immersion does not lift to an embedded projective plane in $4$-space. Our movie presentation (Fig. 27) of the surface indicates a lifting of the surface to an immersed surface in $4$-space that has exactly one triple point. One can analyze the double decker set  (Fig. 26) of Boy's surface directly to determine that it does not lift to an embedding. Doing so herein would take us too far afield. However, we have attempted to make specific liftings of both Boy's and girl's surfaces. A careful examination of the stills indicates that there is a single inconsistency between the crossings in two stills in each of the movies. These double points are indicted upon the decker sets as a pair of small dots (double points in $4$-space) that appear on the double decker set and for which the parity of the labels change.

\subsection{The girl's surface}\label{girl}

\begin{center}
\includegraphics[width=3in]{girls_chart}
\end{center}

\begin{center}
\includegraphics[width=3in]{girl_deck}
\end{center}

\begin{wrapfigure}[17]{l}{3.02in}
\includegraphics[width=3in]{girls_movie}   
\end{wrapfigure}
Girl's surface (Fig. 28) appears in Ap\'{e}ry's book~\cite{Ap}. It was rediscovered and discussed in a paper by Goodman and Kossowski who labeled it $g$. The nick-name ``girls surface" was used privately, as a pun on Boy's surface and since the rediscovery was done by two women. The name has stuck. Its double point set coincides with the double point set of Boy's surface, but the double decker set (Fig. 29) differs since one of the three loops that emanates from the triple point has an extra quarter twist in it. It too, being an immersed surface, does not lift to an embedding. We have depicted an immersion in $4$-space with a single double point (Fig. 30).

\section{Standard Klein Bottles}
\label{Klein}

\begin{center}
\noindent
\includegraphics[width=2in]{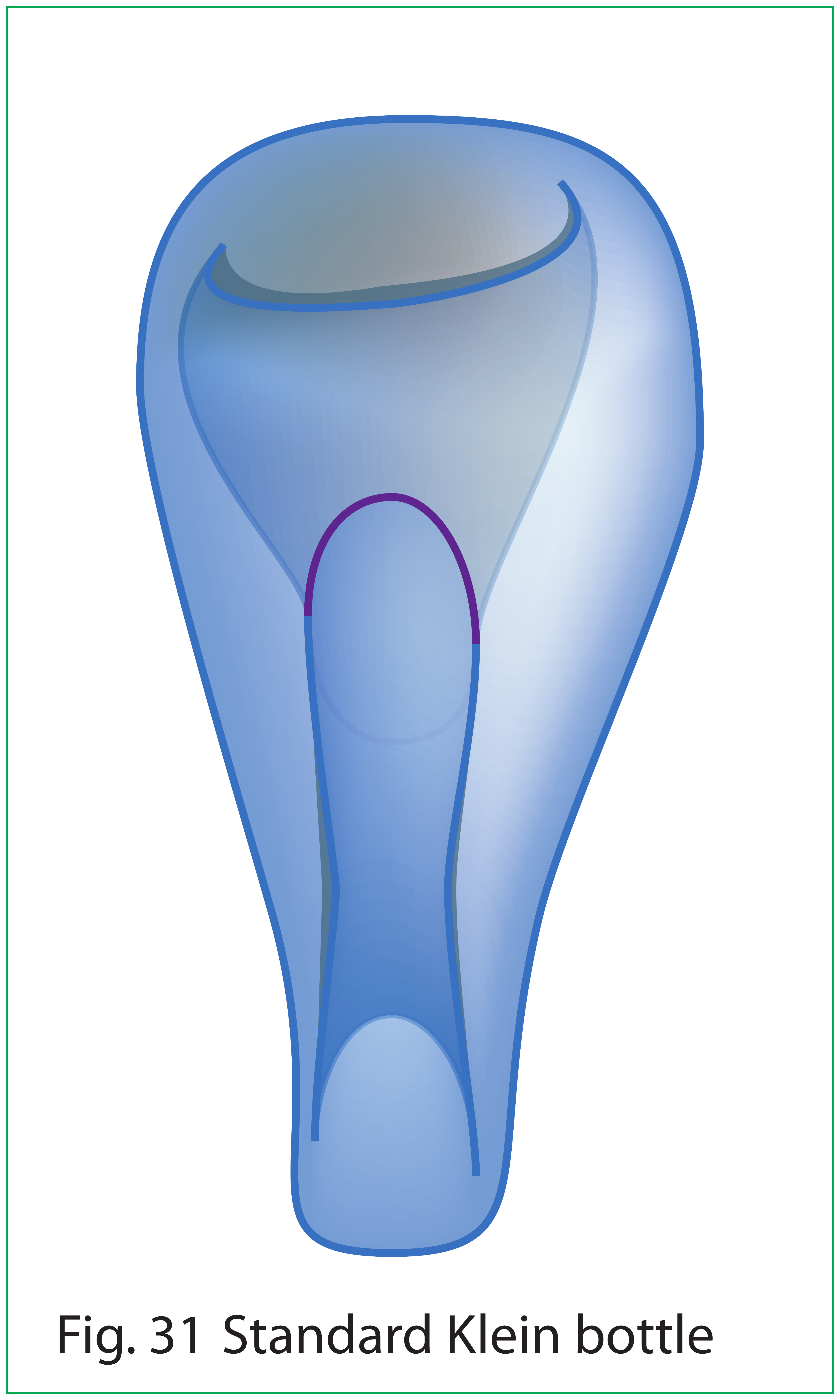}
\end{center}

The standard Klein bottle as represented, for example, by an Acme Klein bottle is indicated in Fig. 31 as it is viewed from the side with the handle protruding towards the viewer. The surface represented by an Acme Klein bottle is a surface in $3$-space, but it can be lifted to an embedding in $4$-space. Here and subsequently, we refer to these surfaces in $4$-space.

\begin{center}
\noindent
\includegraphics[width=4in]{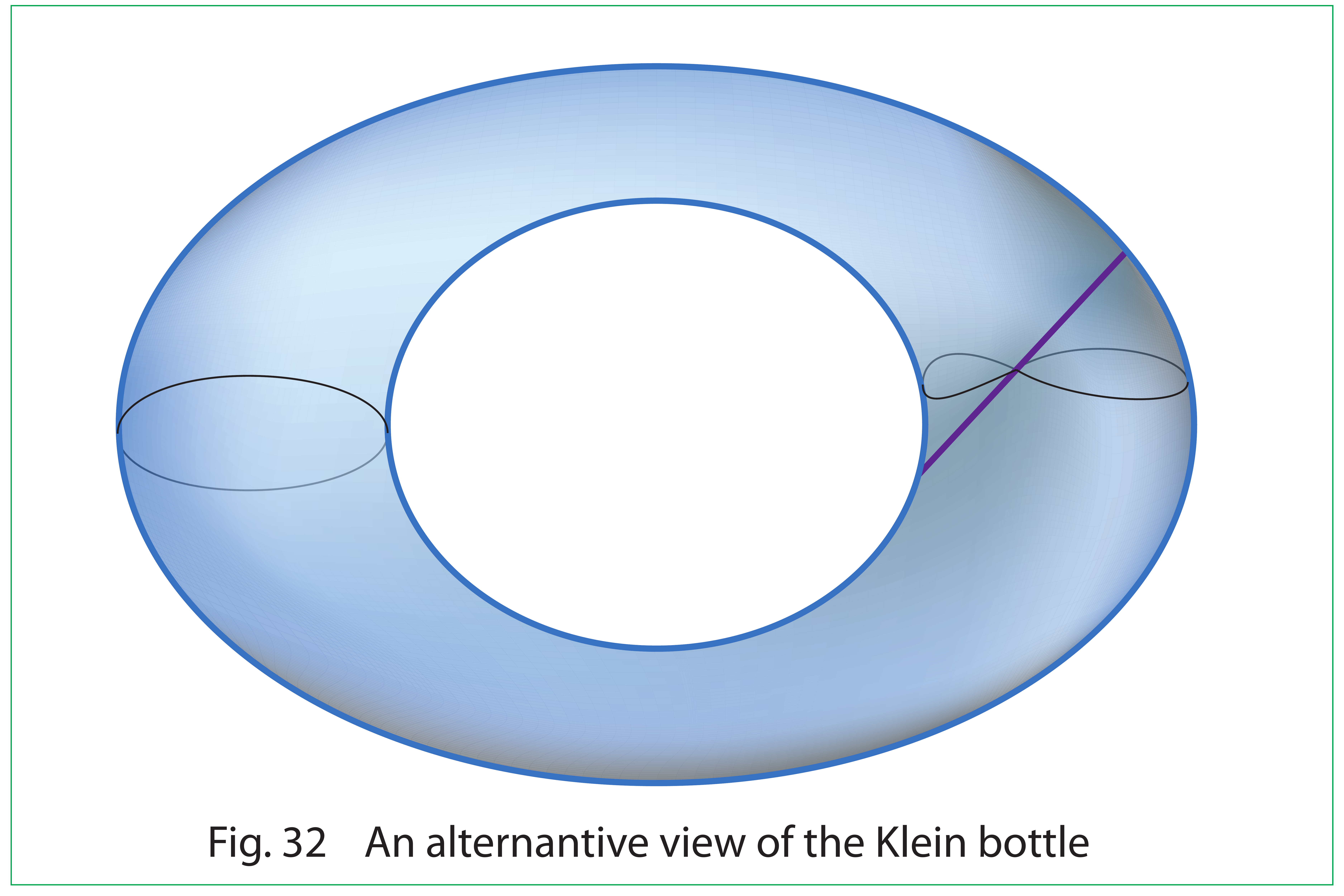}
\end{center}

An alternative view appears in Fig. 32. In this embedding, we imagine that the right-hand-side of a torus is turned over. One can imagine a circular cross-section first being represented by the edge of a coin. As it is flipped, the edge appears straight-on to an observer. Alternatively consider a hyperbolic paraboloid as represented by a popular chip (crisp) that is sold in a can. As the crisp is inverted, a small figure-8 appears on the right edge of the crisp, and gradually moves to the left edge before it disappears. Our diagram here elapses this motion. 

\begin{center}
\noindent
\includegraphics[width=4in]{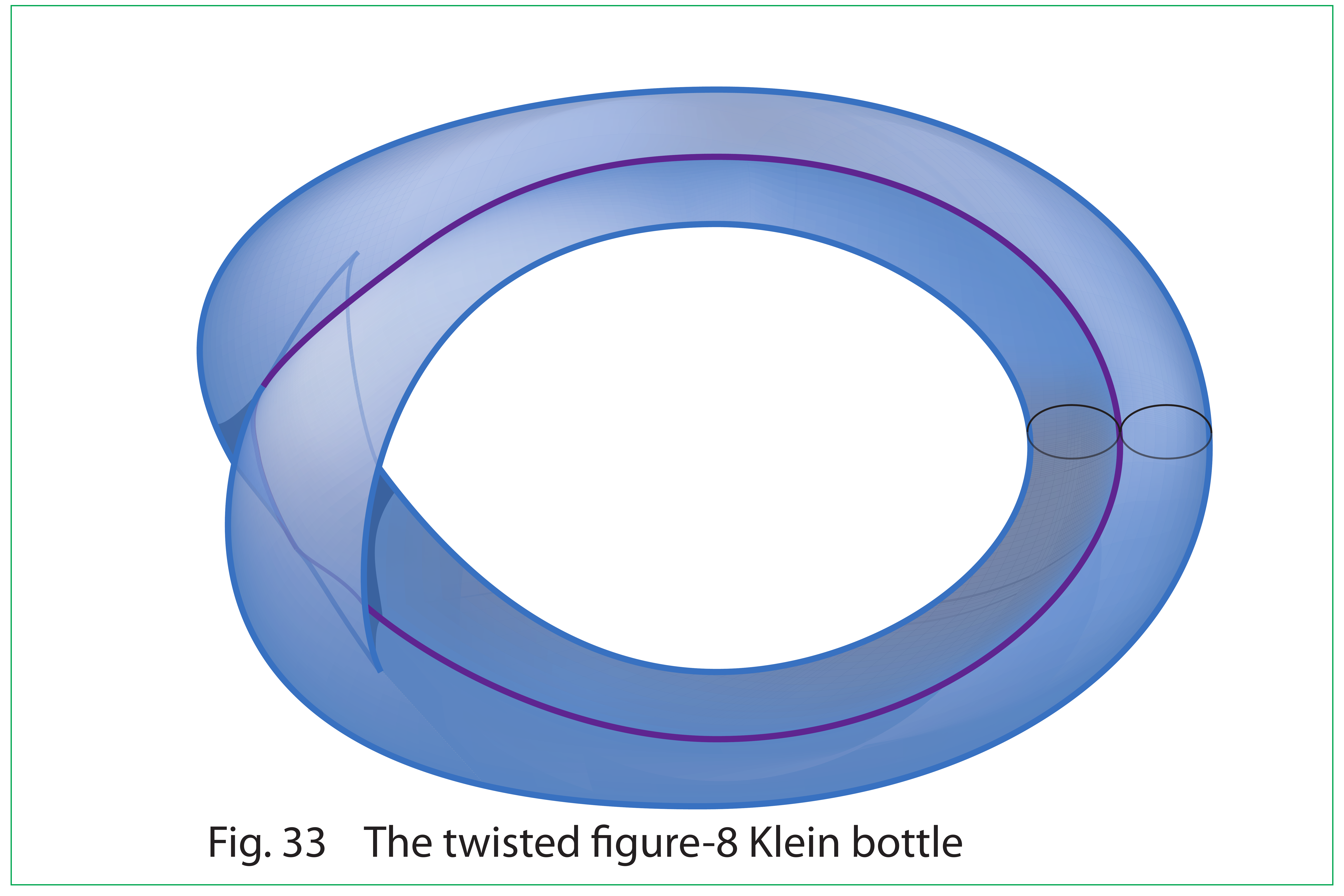}
\end{center}

Consider a standard M\"{o}bius band as represented as a paper model. Now place a second M\"{o}bius band transverse along the core circle. The cross-section of the figure is homeomorphic to a  $+$. In the first and third quadrants of any cross-section attach a quarter circle. The resulting figure-$8$ is invariant under a $180^\circ$ rotation. The resulting Klein bottle is indicated in Fig. 33. 

Each of these three maps of the Klein bottle in $3$-space is the projection of an embedded Klein bottle in $4$-space. Some of our movie representations indicate these liftings. All of these Klein bottles are unknotted: that is they are all ambiently isotopic to the connected sum of a positive and negative cross-cap. The processes of performing the isotopies are indicated in Fig. 34 through Fig. 42.

\begin{center}
\includegraphics[width=4.25in]{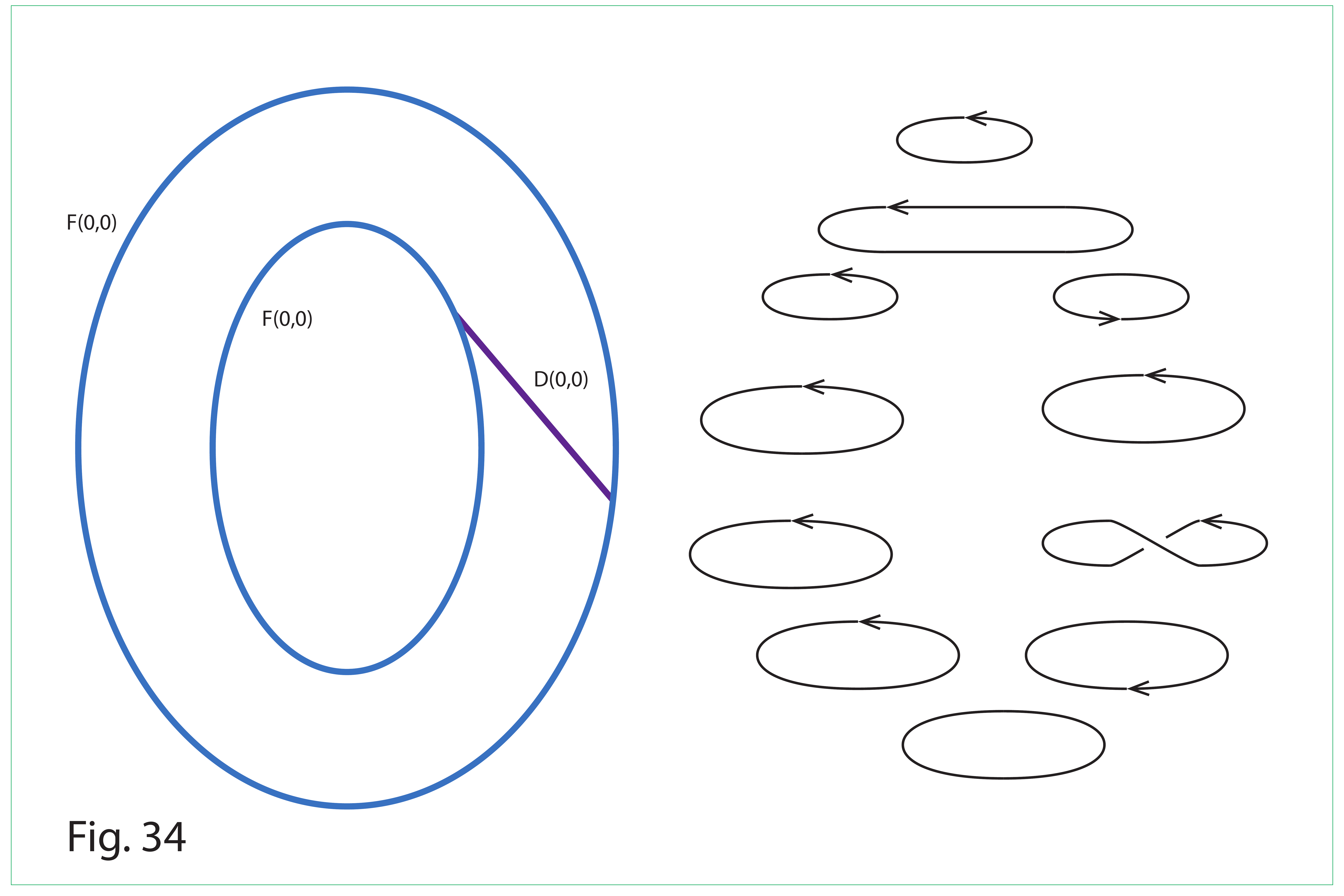}
\includegraphics[width=4.25in]{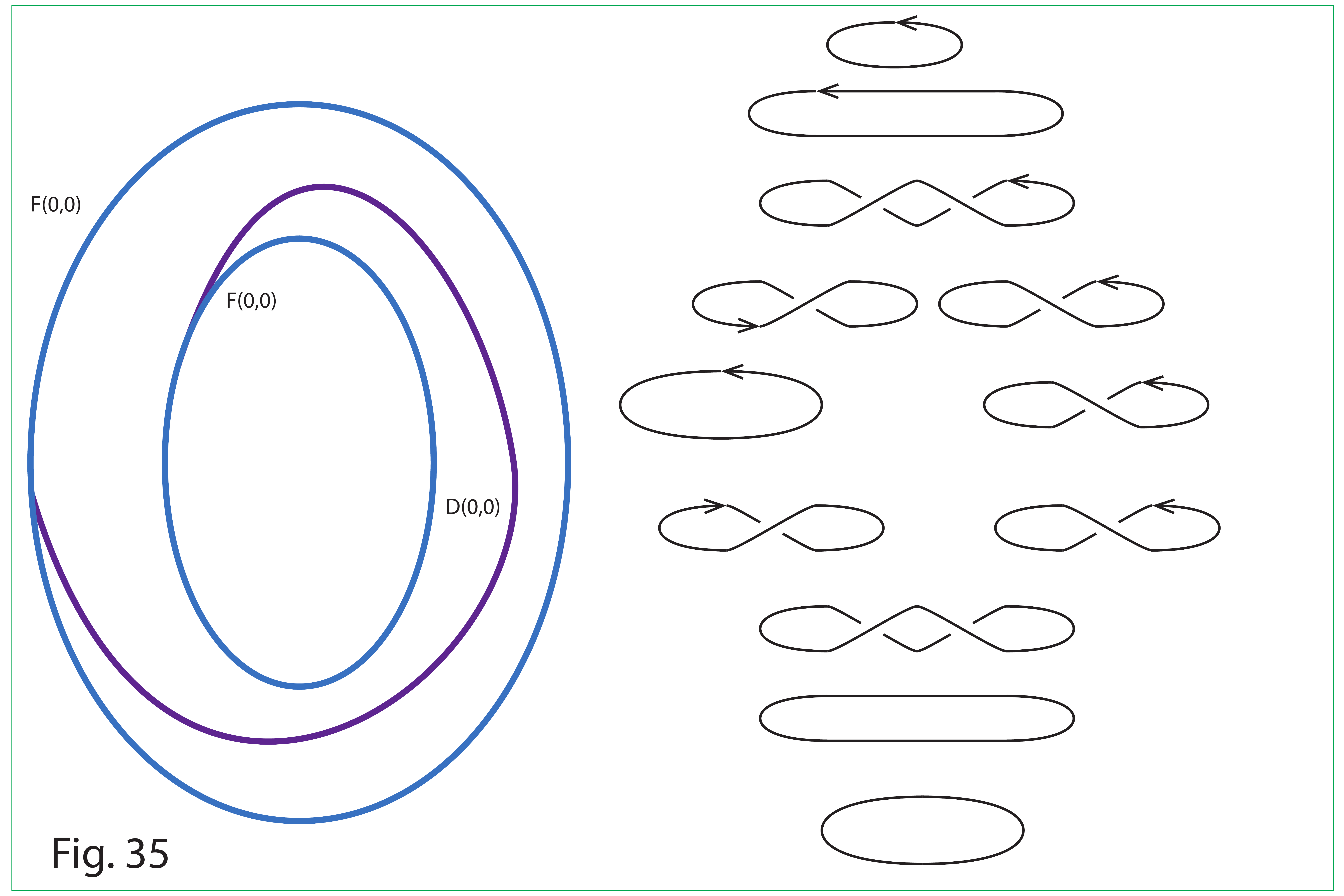}
\includegraphics[width=4.5in]{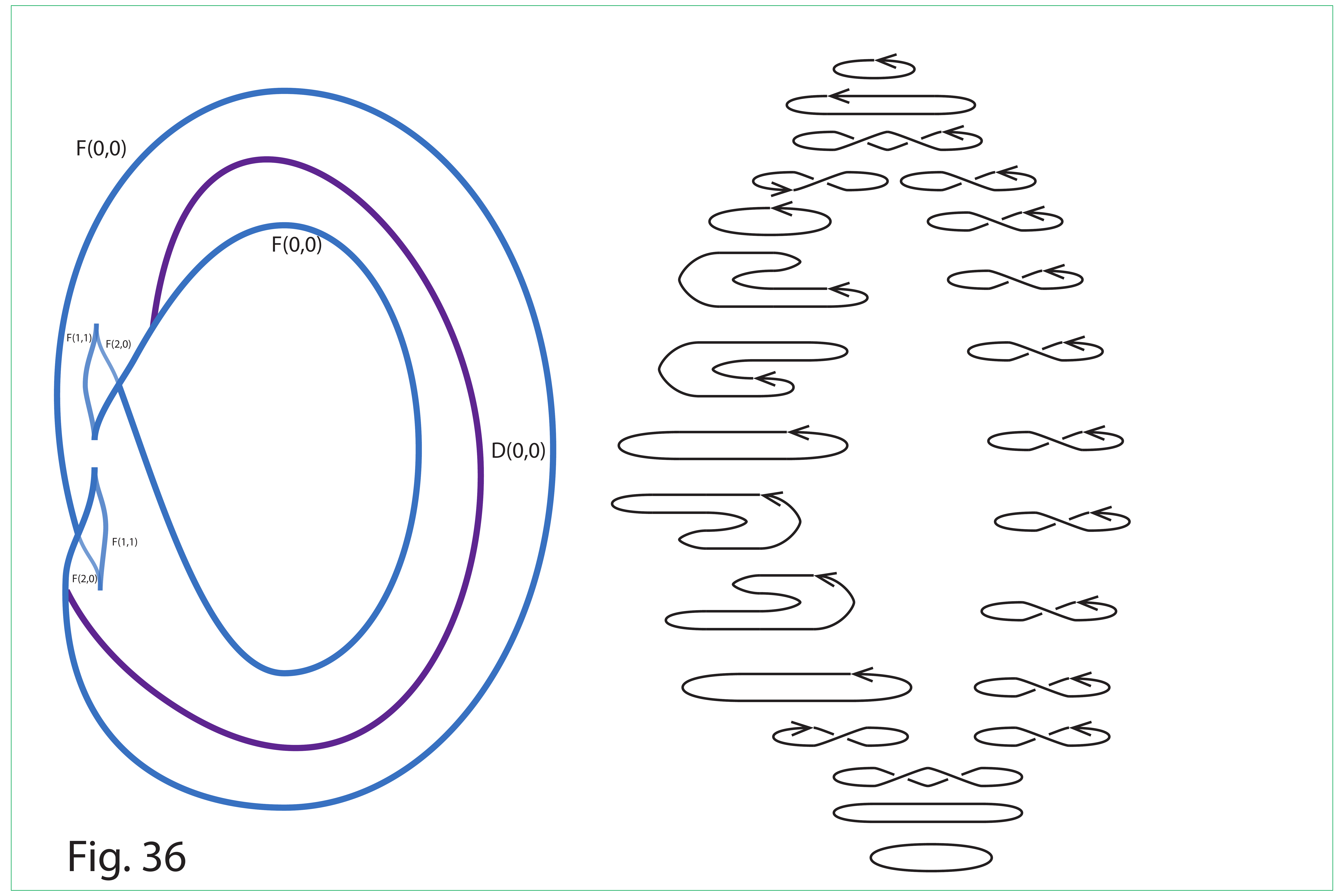}
\includegraphics[width=4.5in]{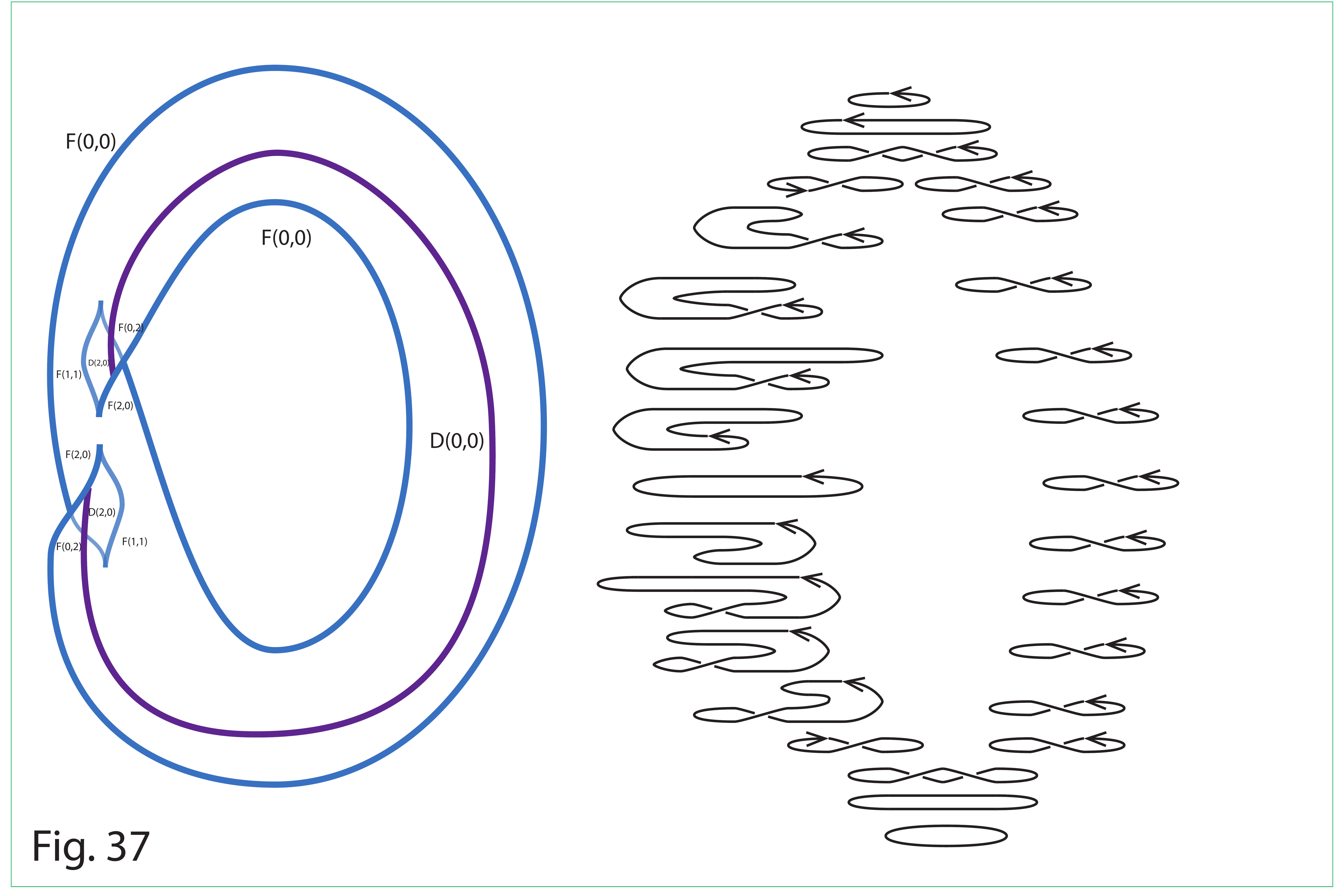}
\includegraphics[width=4.5in]{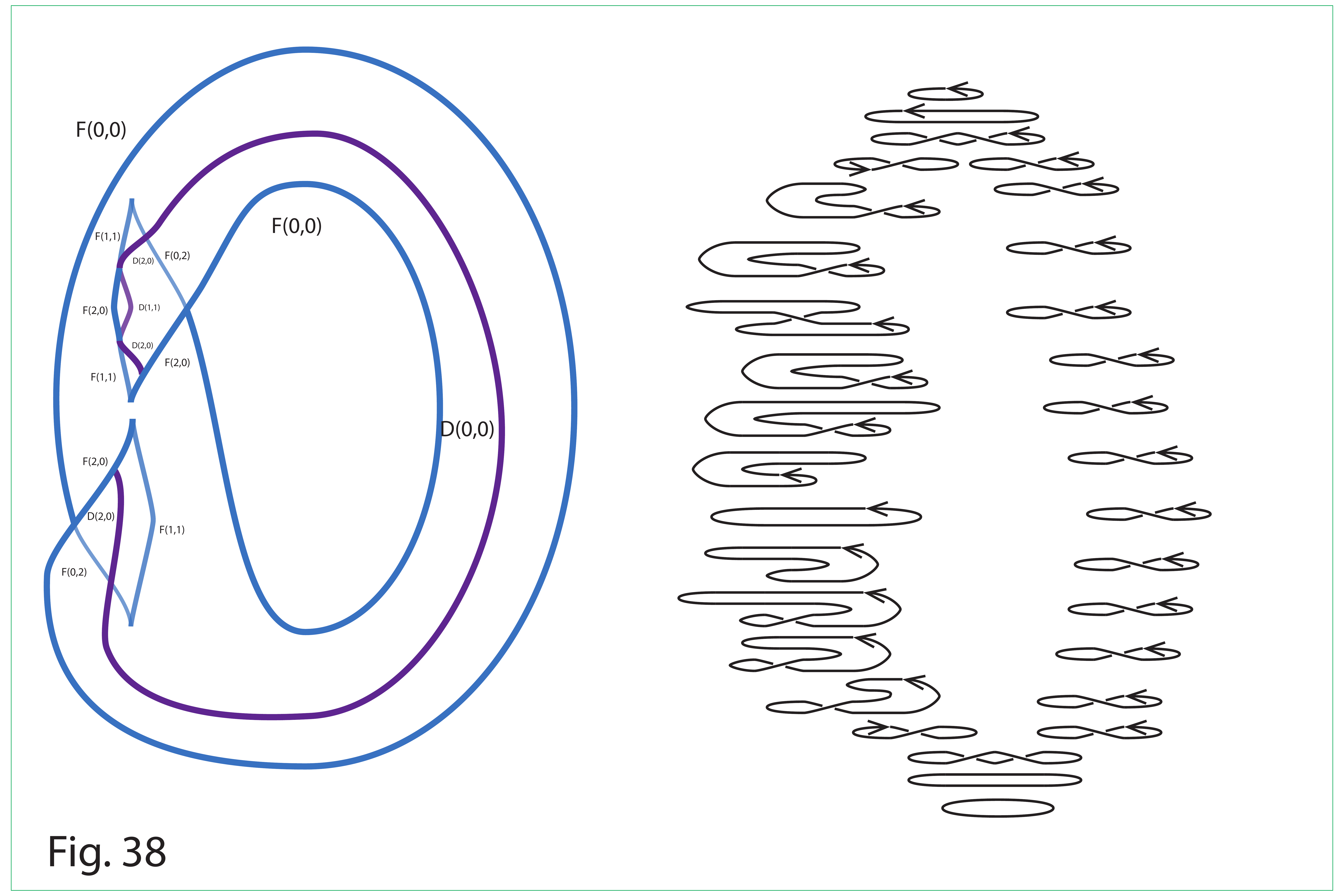}
\includegraphics[width=4.5in]{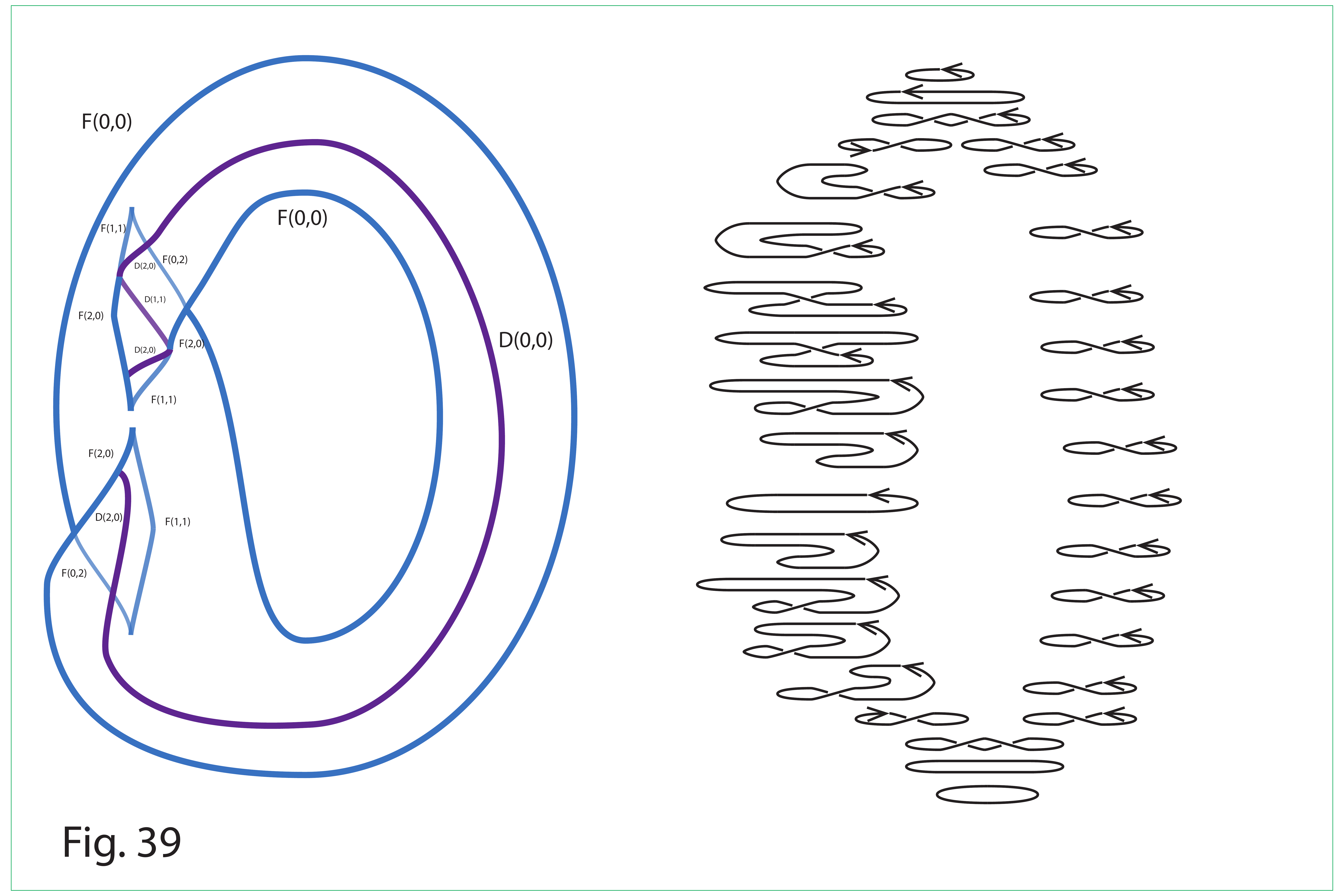}
\includegraphics[width=4.5in]{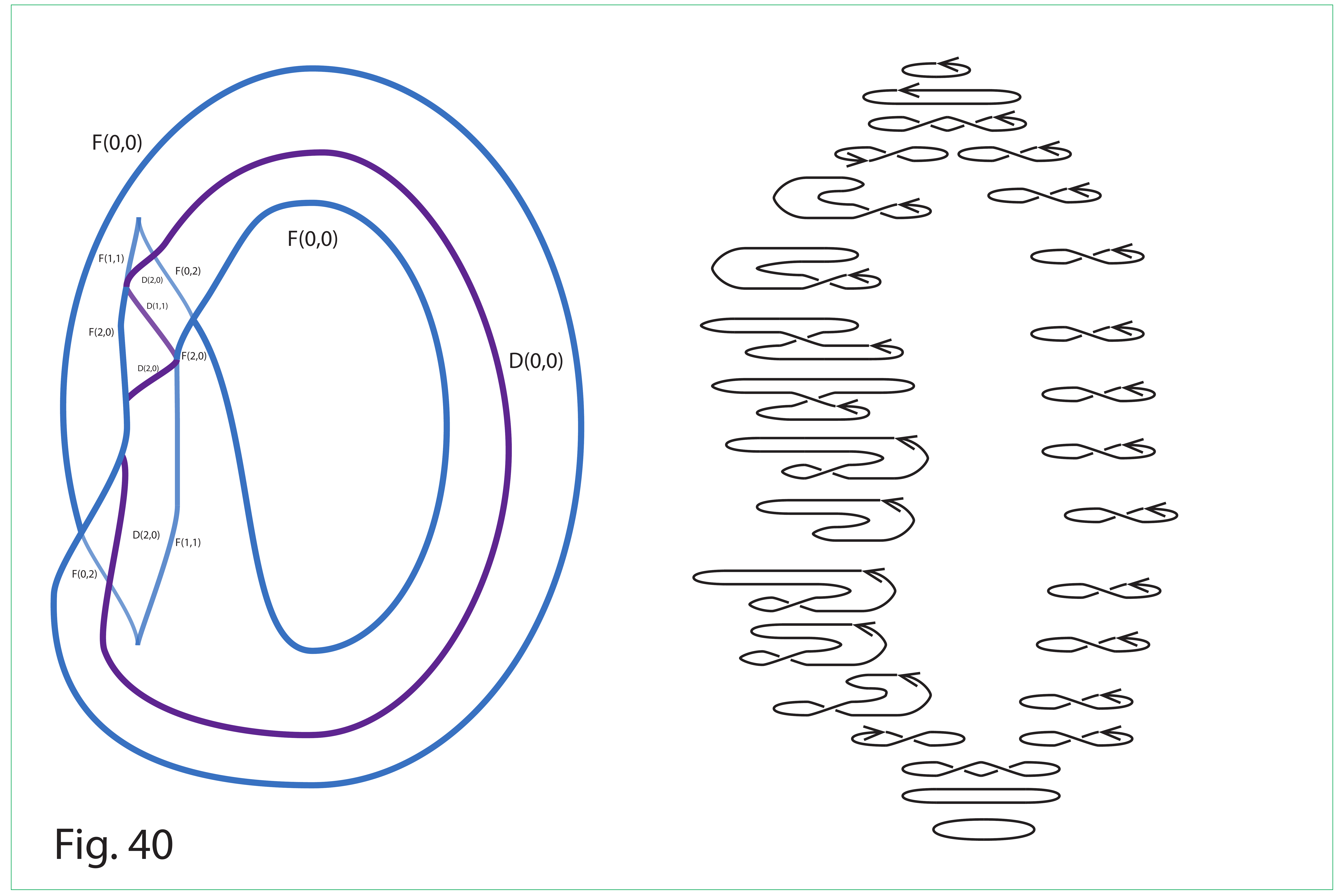}
\includegraphics[width=4.5in]{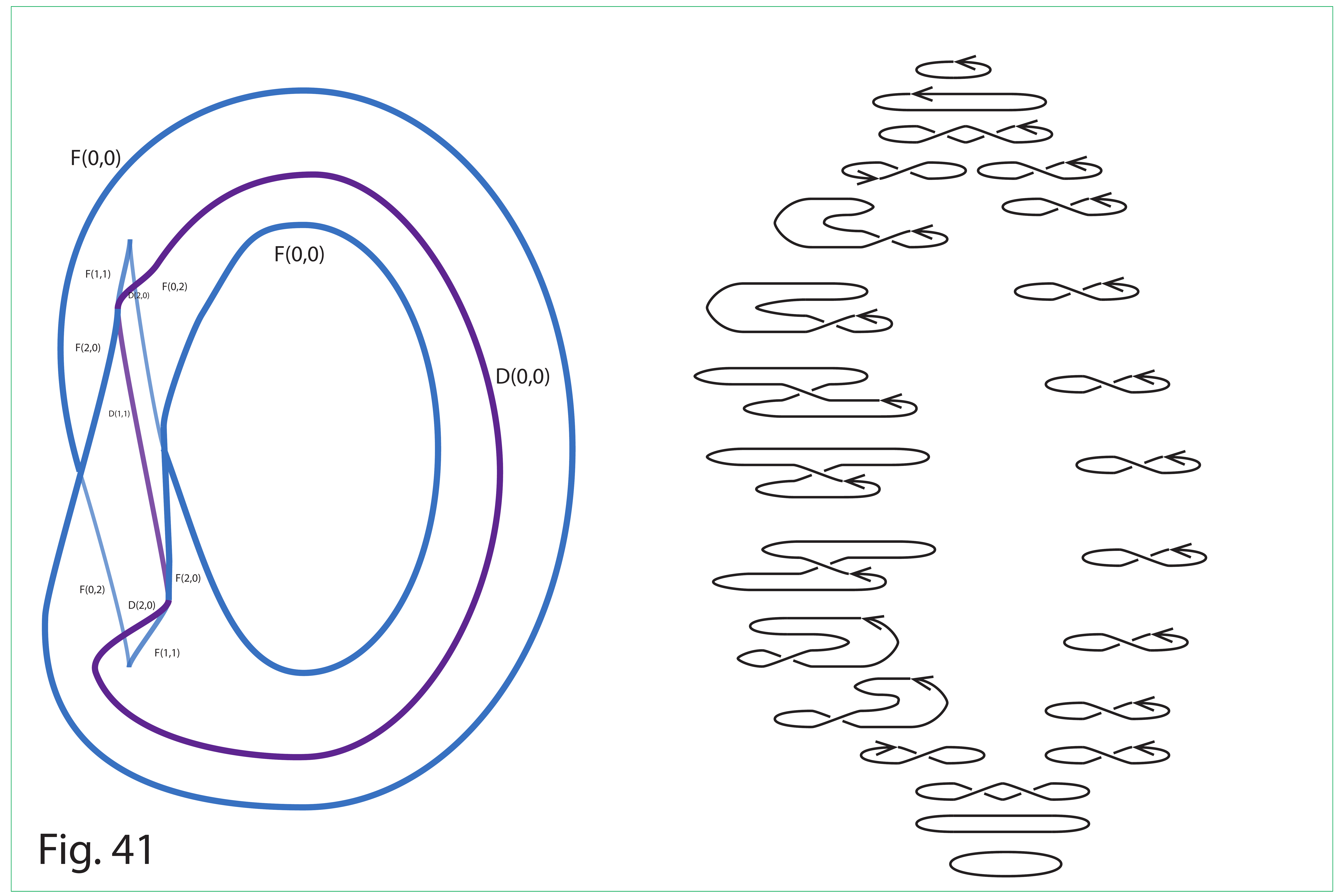}
\includegraphics[width=4.5in]{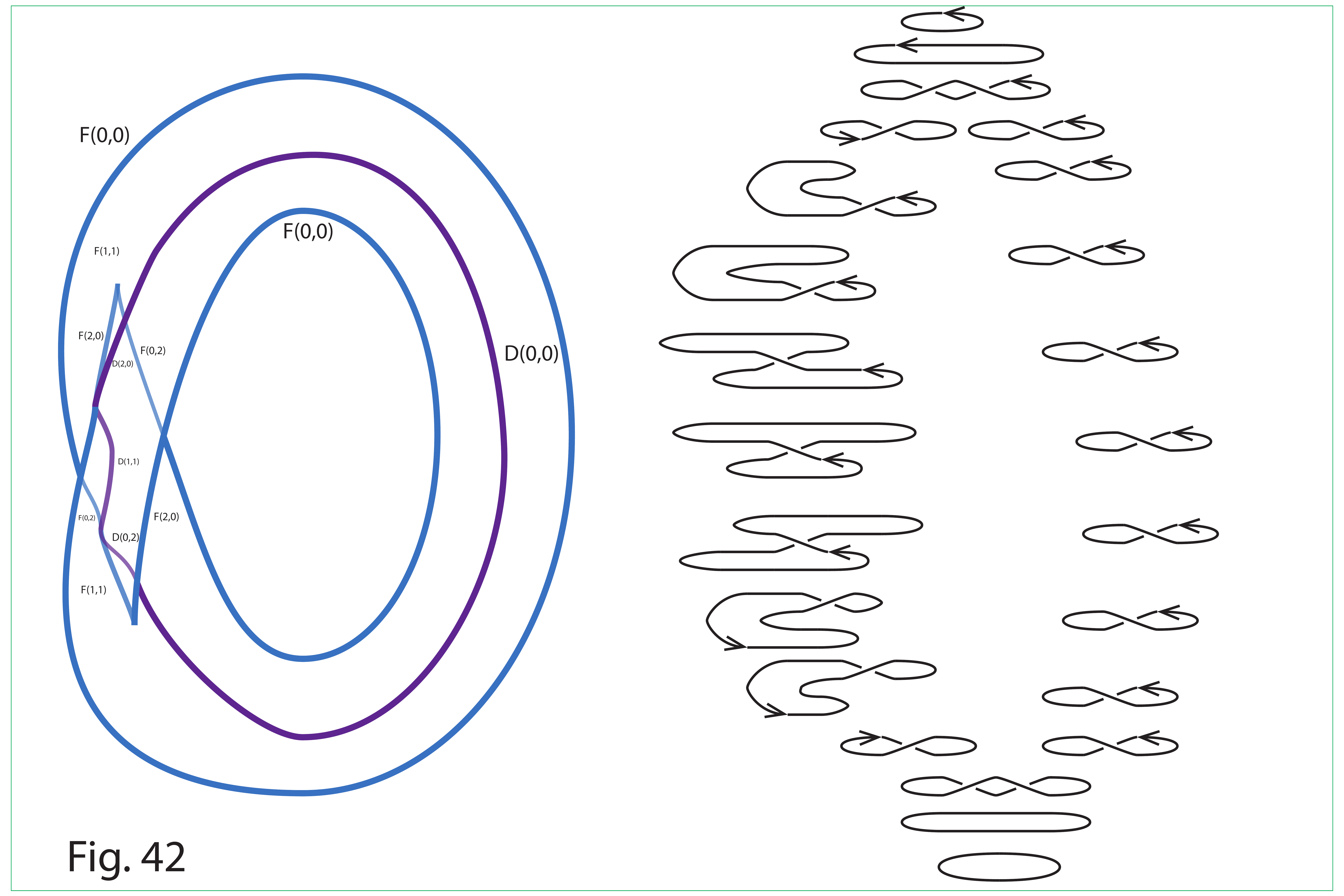}
\end{center}

The first ambient isotopy that is illustrated travels from the ``torus with a wen" to the ``twisted-figure-$8$" embeddings of the Klein bottle. In this set of illustrations movies that correspond to the given charts are also indicated. In this way, the reader can become familiar with the correspondence between movies and charts. Orientations are introduced on the cross-sectional curves for two reasons. First, the fold line upon which a branch point is introduced is ambiguous if an un-oriented move is given. The figure-$8$ on the right is born with its branch point on the left and dies with it on the right. This can be seen by comparing the orientations of the circles above and below it. Second, the orientation cannot be carried through the entire movie. The saddle point at the bottom is necessarily not orientation preserving. Let us describe the isotopy.

From Fig. 34 to Fig. 35, the branch points pass over saddles and a minimal point on the surface following the moves (Fig. 40 and 41, p. 81 \cite{CS:Book}) also optima on the double point set (``candy-cane"-moves, Fig. 50, p. 85 \cite{CS:Book}) occur. From Fig. 35 to Fig. 36, a  pair of swallow-tails (Fig. 39, p. 81 \cite{CS:Book}) are introduced. One can imagine that these swallow-tails correspond to twisting the tube that has no branch points. From Fig. 36 to Fig. 37, the branch points move past the folds that are in the back of the figure. From Fig. 37 to Fig. 38, a pair of redundant $\psi$-moves are added (Fig. 46, p. 83 \cite{CS:Book}). From Fig. 38 to Fig. 39, the top branch point moves through a cusp (Fig. 42, p. 82 \cite{CS:Book}). From Fig. 39 to Fig 40, a beak-to-beak hyperbolic confluence of cusps occurs (Fig. 39, p. 81 \cite{CS:Book}). This is followed (Fig. 40 to Fig. 41) by a hyperbolic confluence of branch points (Fig. 53, p. 86 \cite{CS:Book}). Finally, From Fig. 41 to Fig. 42, a double line passes over a fold line near a cusp (Fig. 45, p. 83 \cite{CS:Book}).

\noindent
\includegraphics[width=2.25in]{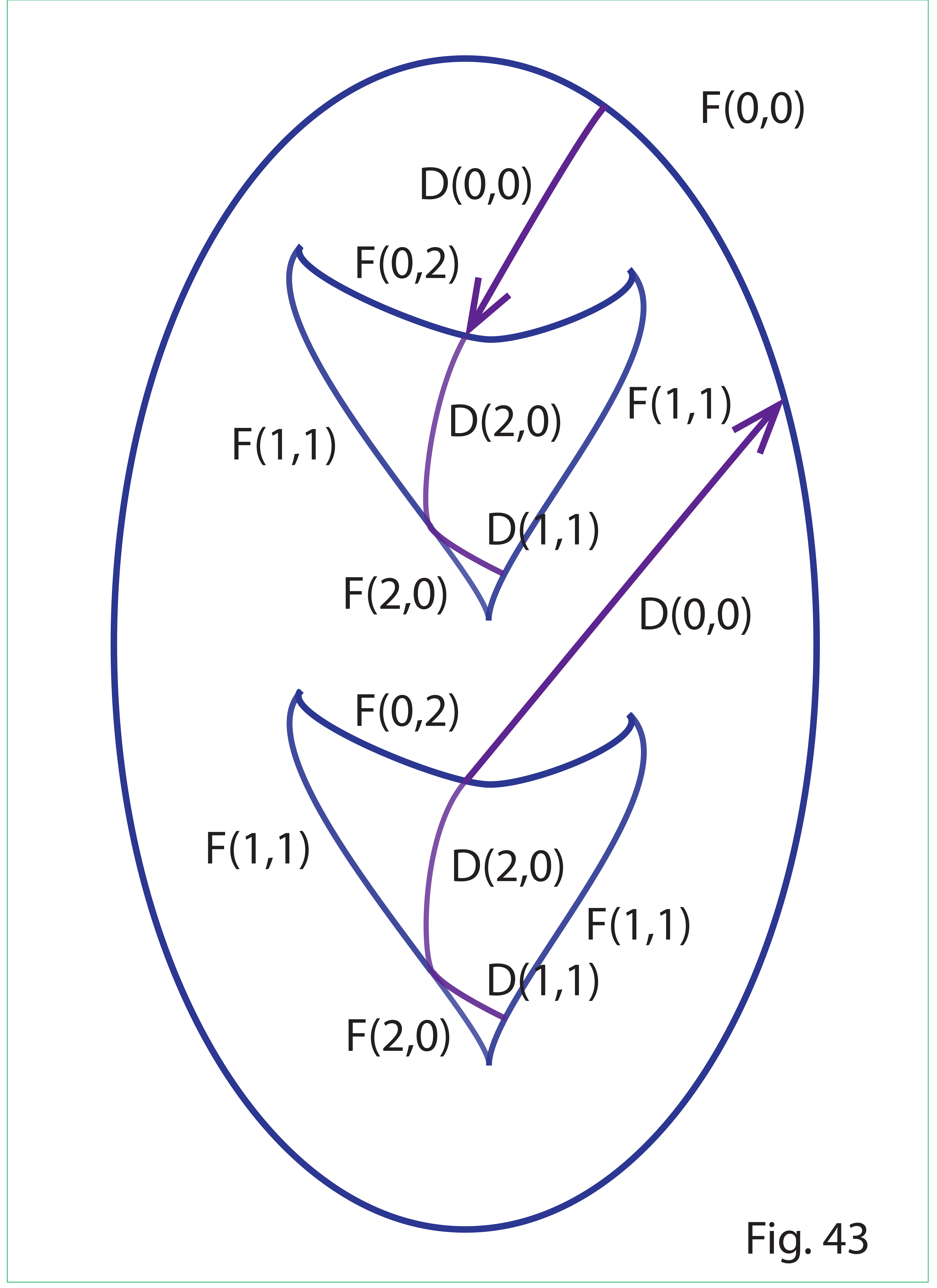}
\includegraphics[width=2.25in]{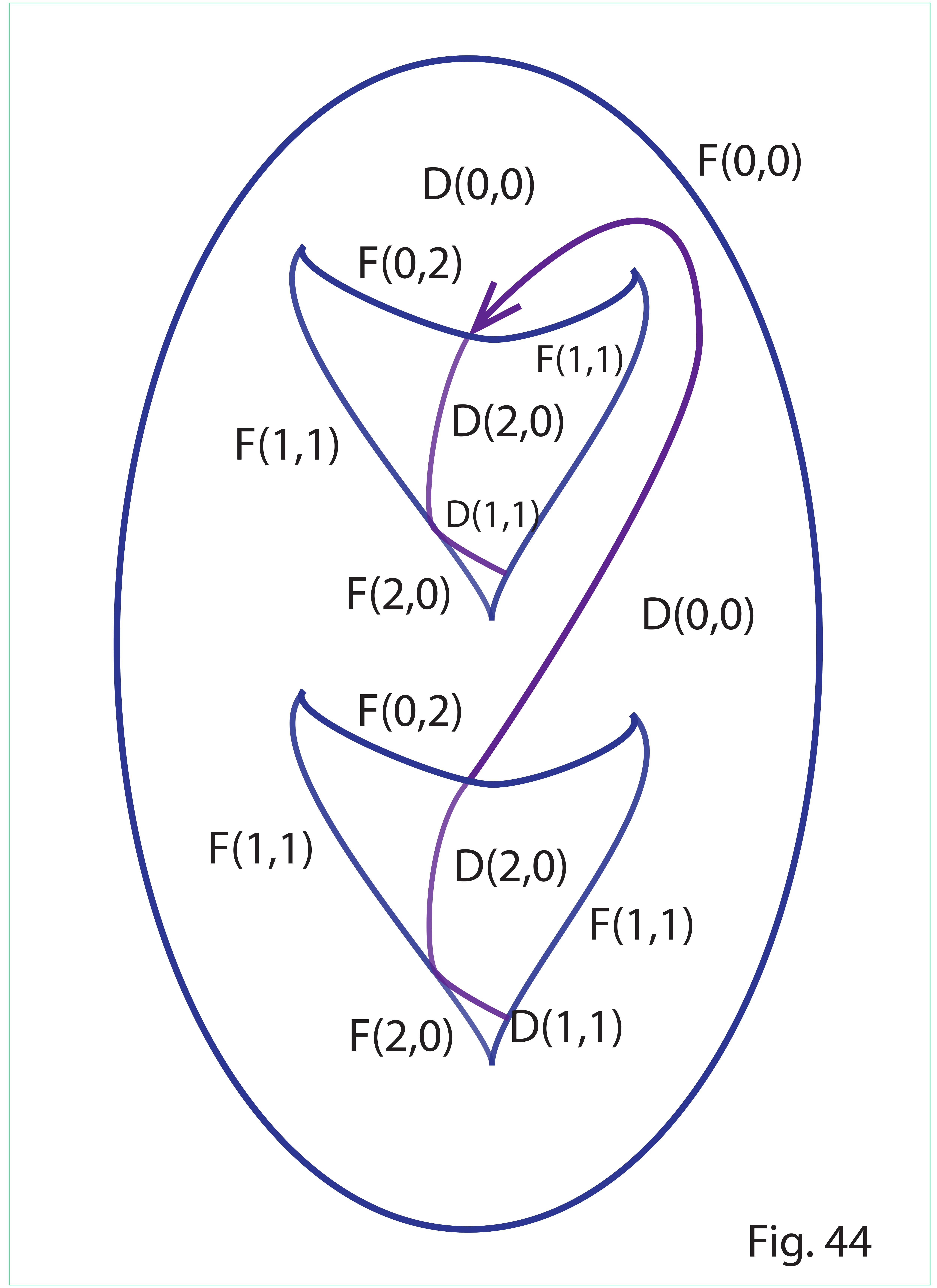}
\includegraphics[width=2.25in]{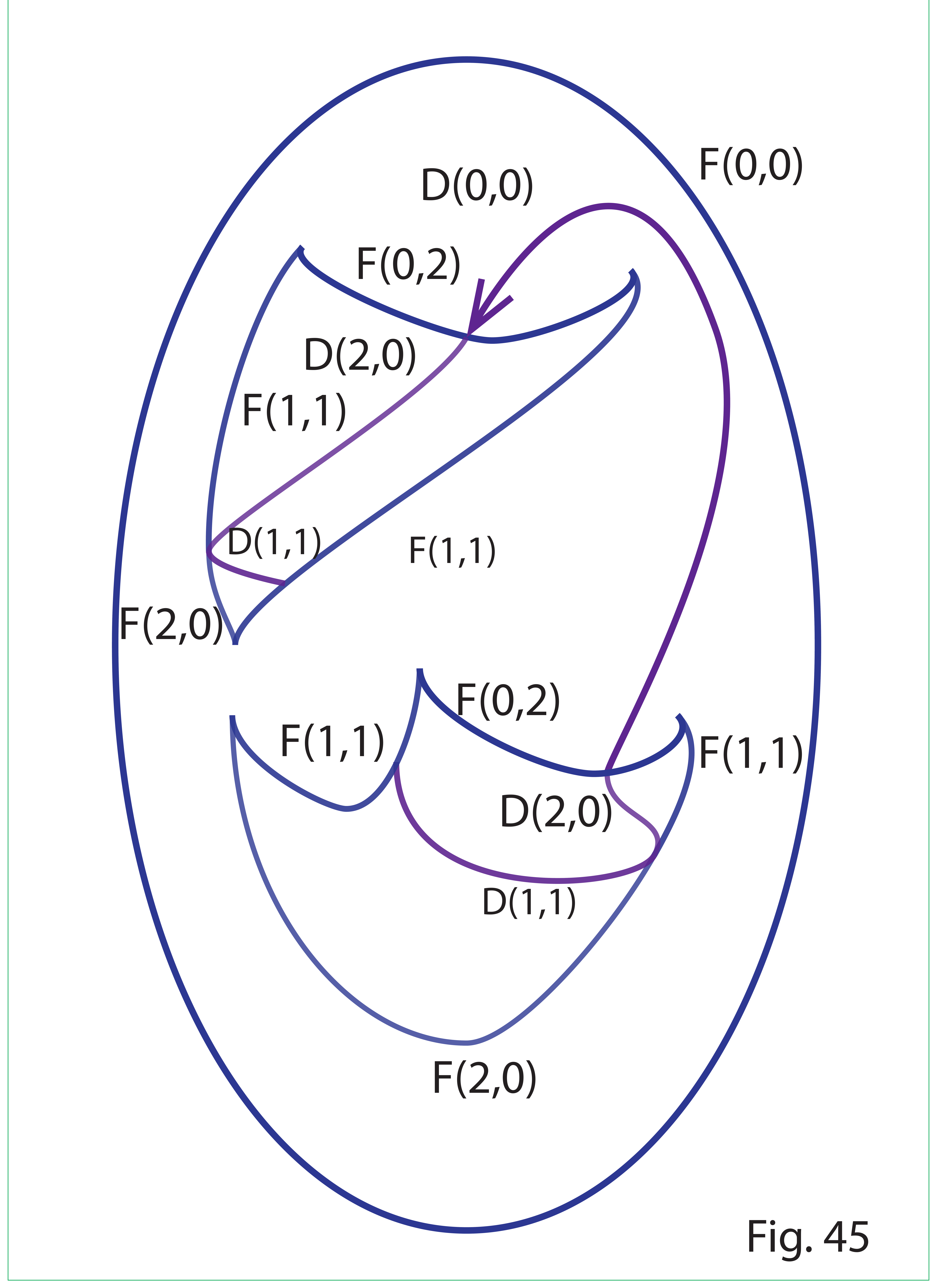}
\includegraphics[width=2.25in]{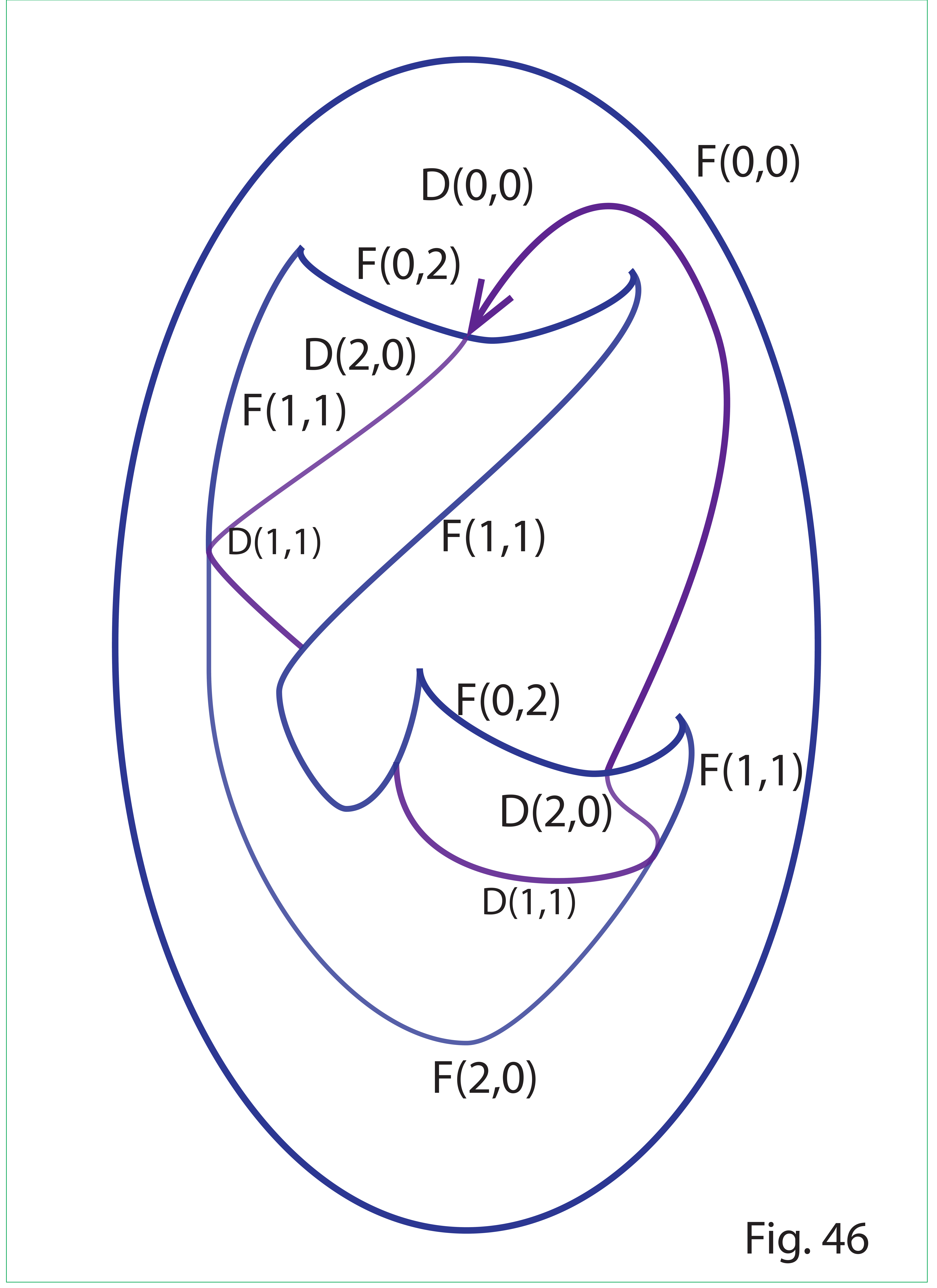}
\includegraphics[width=2.25in]{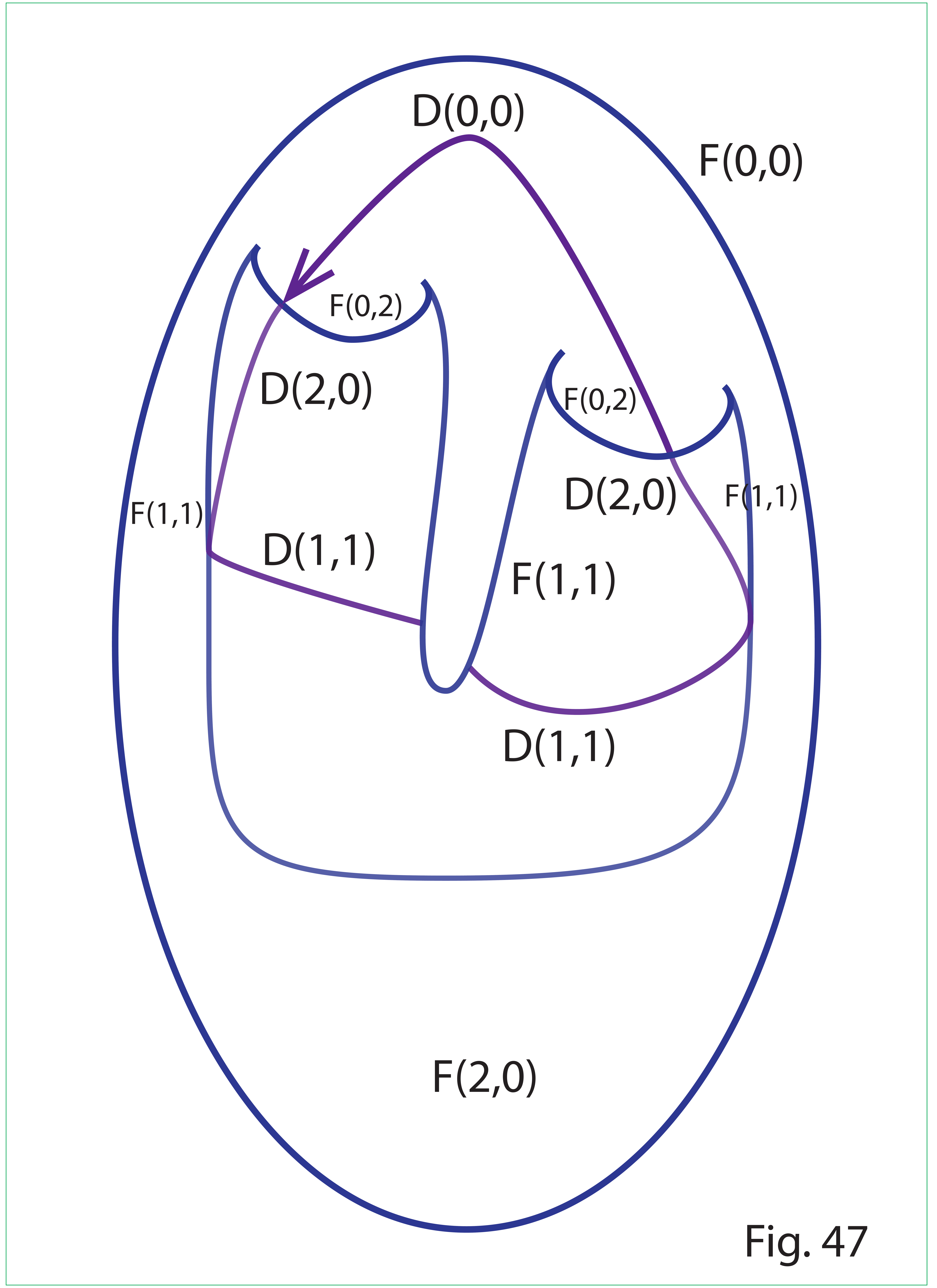} 
\includegraphics[width=2.25in]{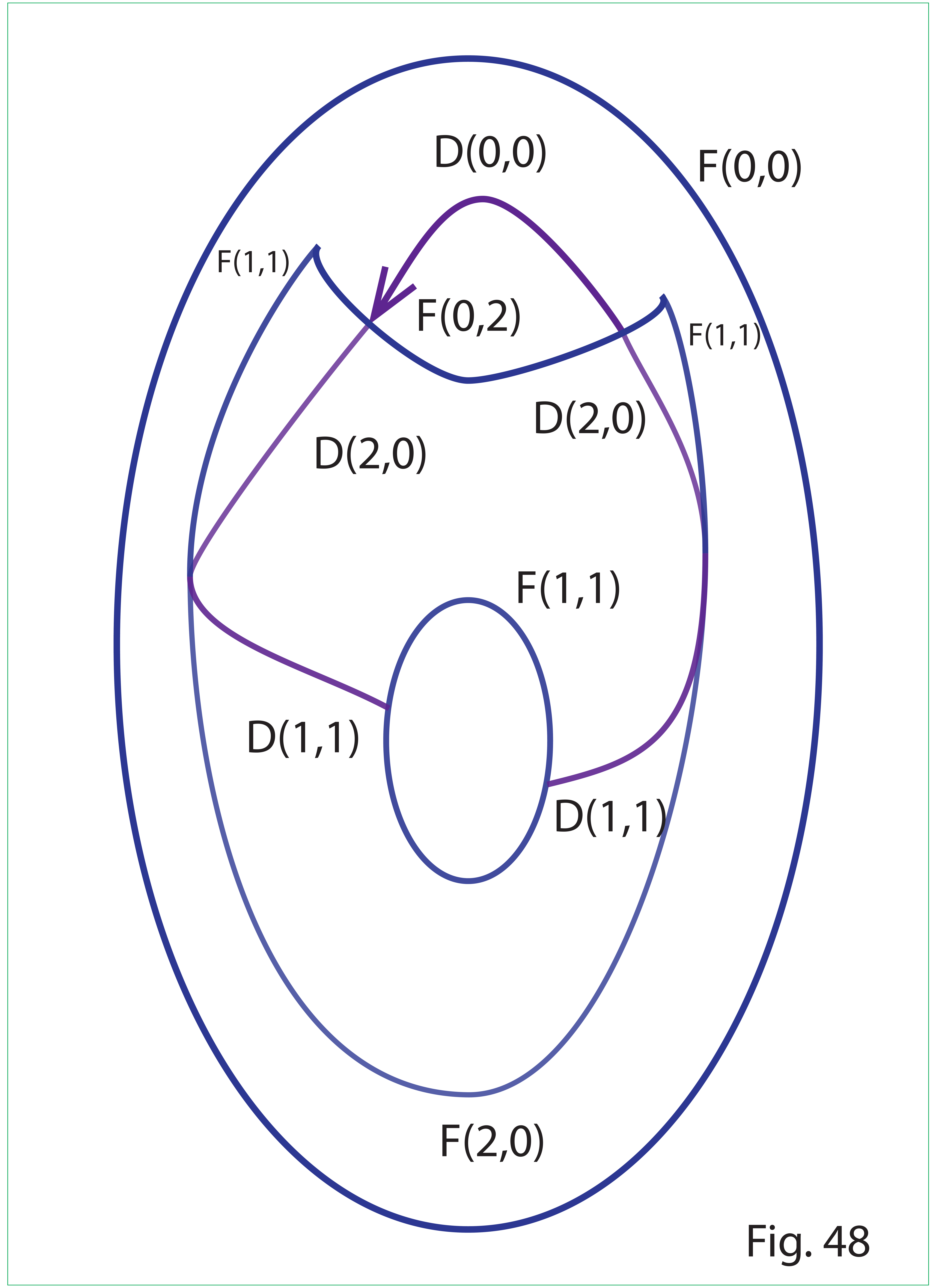} 
\includegraphics[width=2.25in]{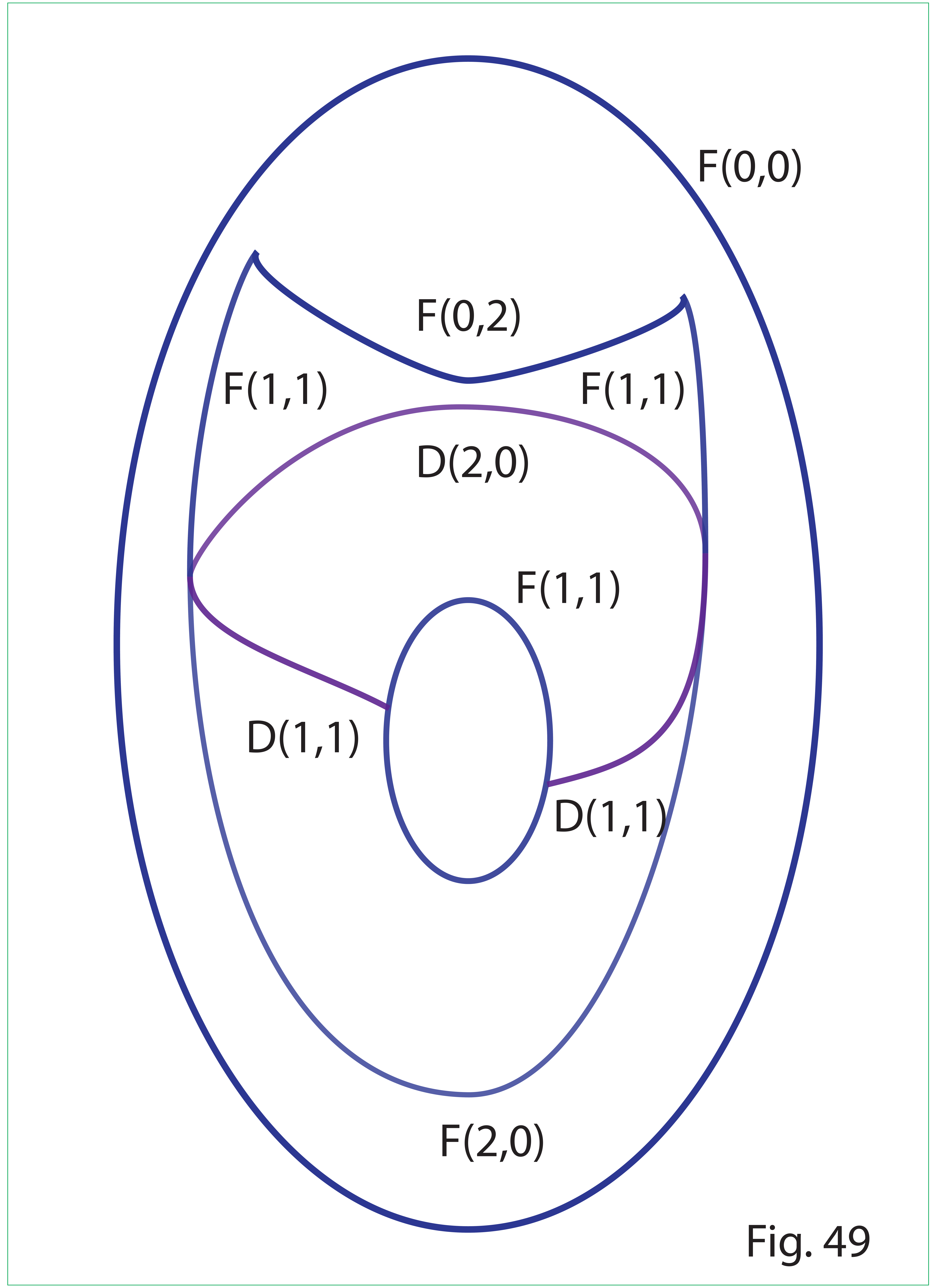}
\includegraphics[width=2.25in]{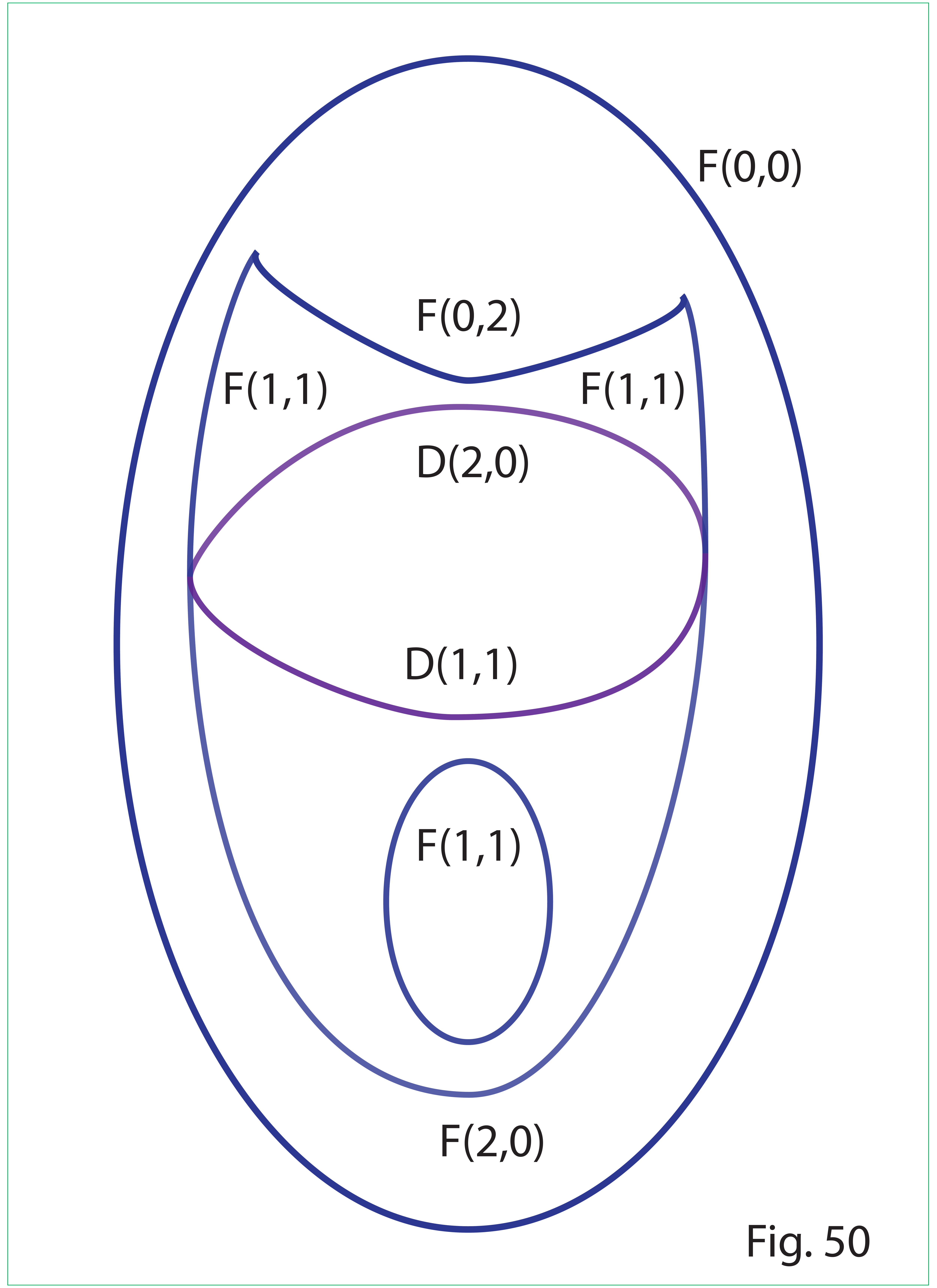}
\includegraphics[width=2.25in]{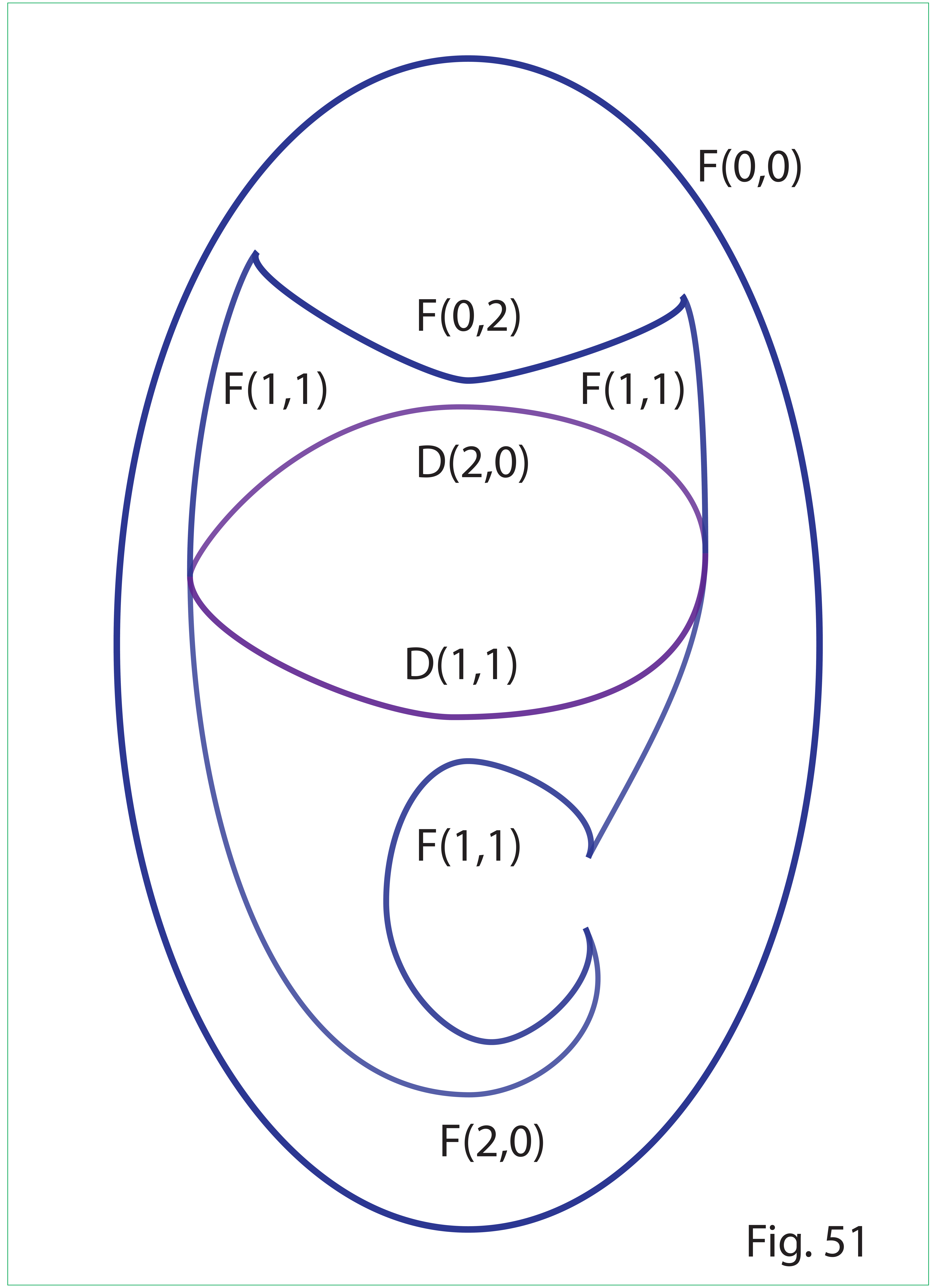} 
\includegraphics[width=2.25in]{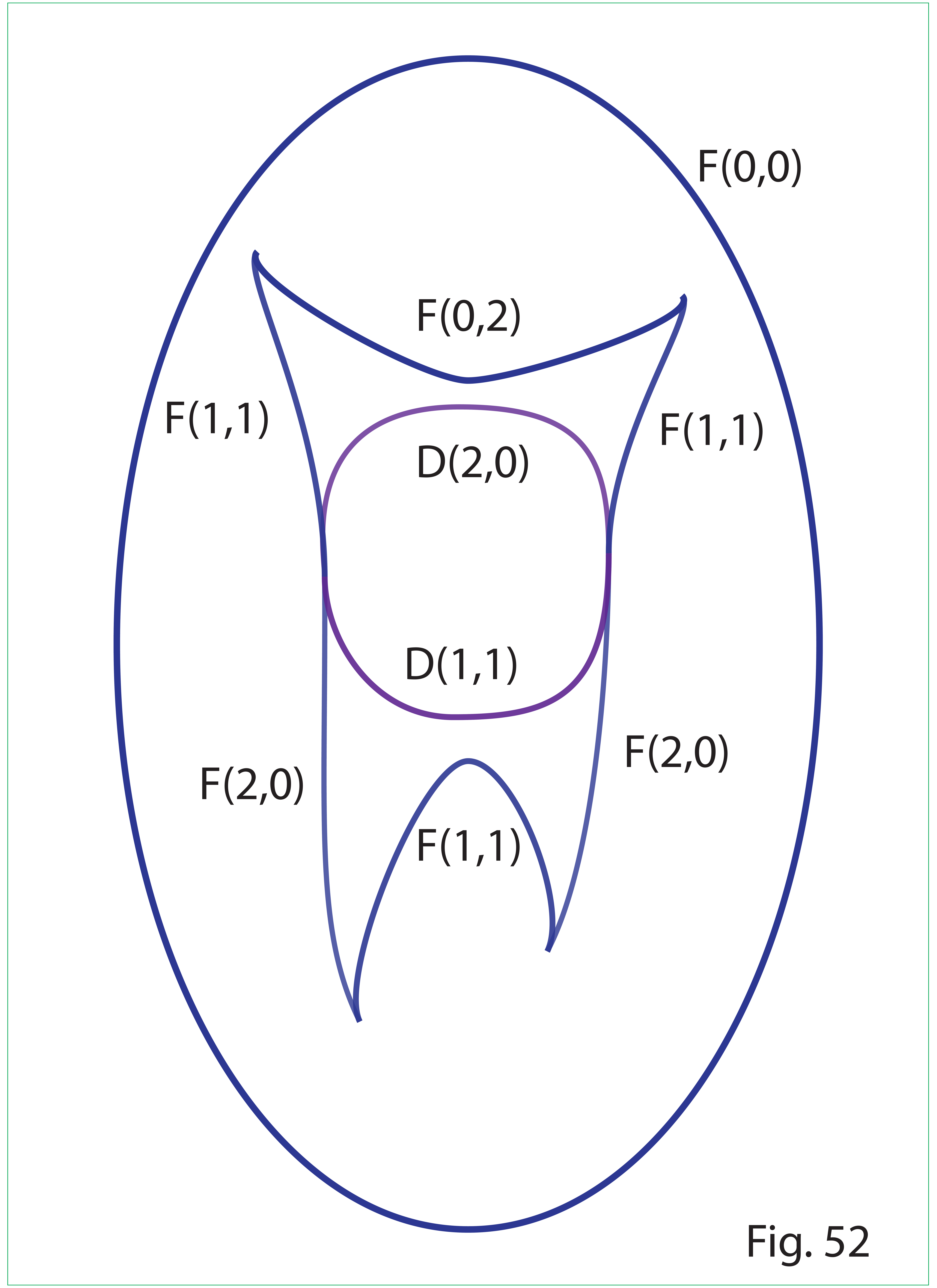} 

The connected sum of a positive and a negative cross-cap is shown in Fig. 43. The signs of the cross-caps is illustrated (somewhat artificially) as an arrow on the image of the double curve. Fold lines and double curves are indicated with their respective depths in-front and behind the folds. It is easy, in this case, to reconstruct movies of the surfaces. For example, see the movies for the cross-caps that are depicted above. Let us discuss an isotopy  to the standard Klein bottle which we represent via a sequence of closely connected charts.
Fig. 44 is related to Fig. 43 by the hyperbolic confluence of branch points (Fig. 53, page 86 \cite{CS:Book}). In addition, a maximal point of the double point set is pushed through a branch point (Fig. 50, p. 85 \cite{CS:Book}). 
To achieve the surface illustrated in Fig. 45, a horizontal cusp occurs towards the bottom of the illustration  (Fig. 38, p. 80 \cite{CS:Book}), and the bottom branch point passes through a cusp (Fig. 42, p. 80 \cite{CS:Book}).  Fig. 46 follows from Fig. 45 via a beak-to-beak  hyperbolic confluence of cusps (Fig. 36, p. 75 \cite{CS:Book}).
Fig. 47 differs from the  Fig. 46 by means of interchanging the heights of a variety of critical events. Fig. 48 differs from Fig. 47 by a sequence of moves. The third cusp from the left is inverted via a horizontal cusp (Fig. 38 p. 80 \cite{CS:Book}). Then the two central cusps are coalesced via a beak-to-beak hyperbolic confluence of cusps (Fig. 36, p. 75 \cite{CS:Book}). Fig. 49 differs from Fig. 48 by moving the double curve behind the folded sheet. The moves involved are among the chart moves illustrated in Fig. 59, p. 89 \cite{CS:Book}.  Fig. 50 differs from Fig. 49 by moving a branch point over a saddle (Fig. 41, p. 81 \cite{CS:Book}) and subsequently pushing a minimal point on the double curve through a branch point (Fig. 50, p. 85 \cite{CS:Book}). Finally, a hyperbolic confluence of branch points (Fig. 53, p. 86 \cite{CS:Book}) joins the double curves together 
Fig. 51 follows from Fig. 50 by means of a beak-to-beak hyperbolic confluence of cusps (Fig. 36, p. 75 \cite{CS:Book}). 
Fig. 52 follows by means of a horizontal cusp  (Fig. 38 p. 80 \cite{CS:Book}).

\includegraphics[width=2.25in]{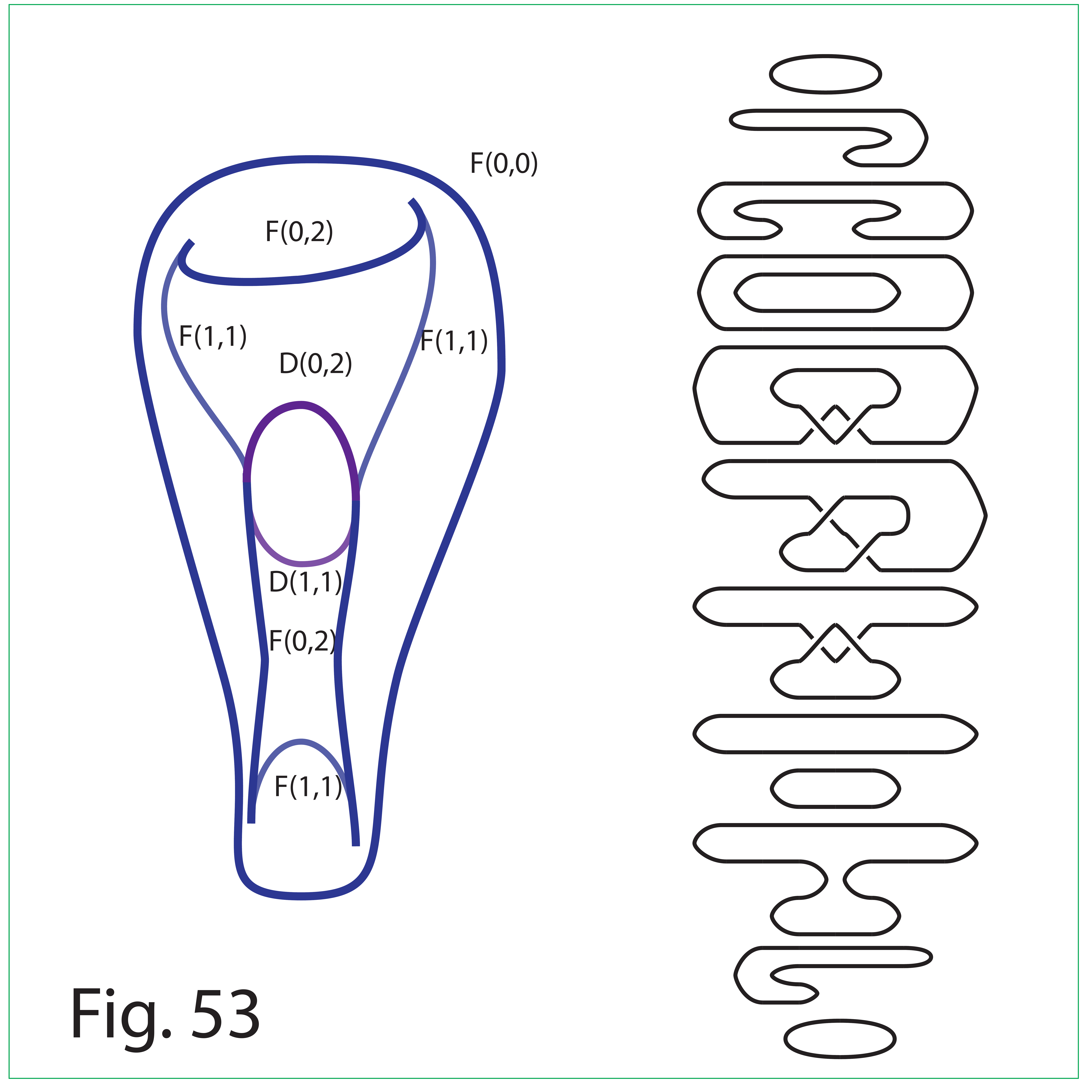}
\includegraphics[width=2.25in]{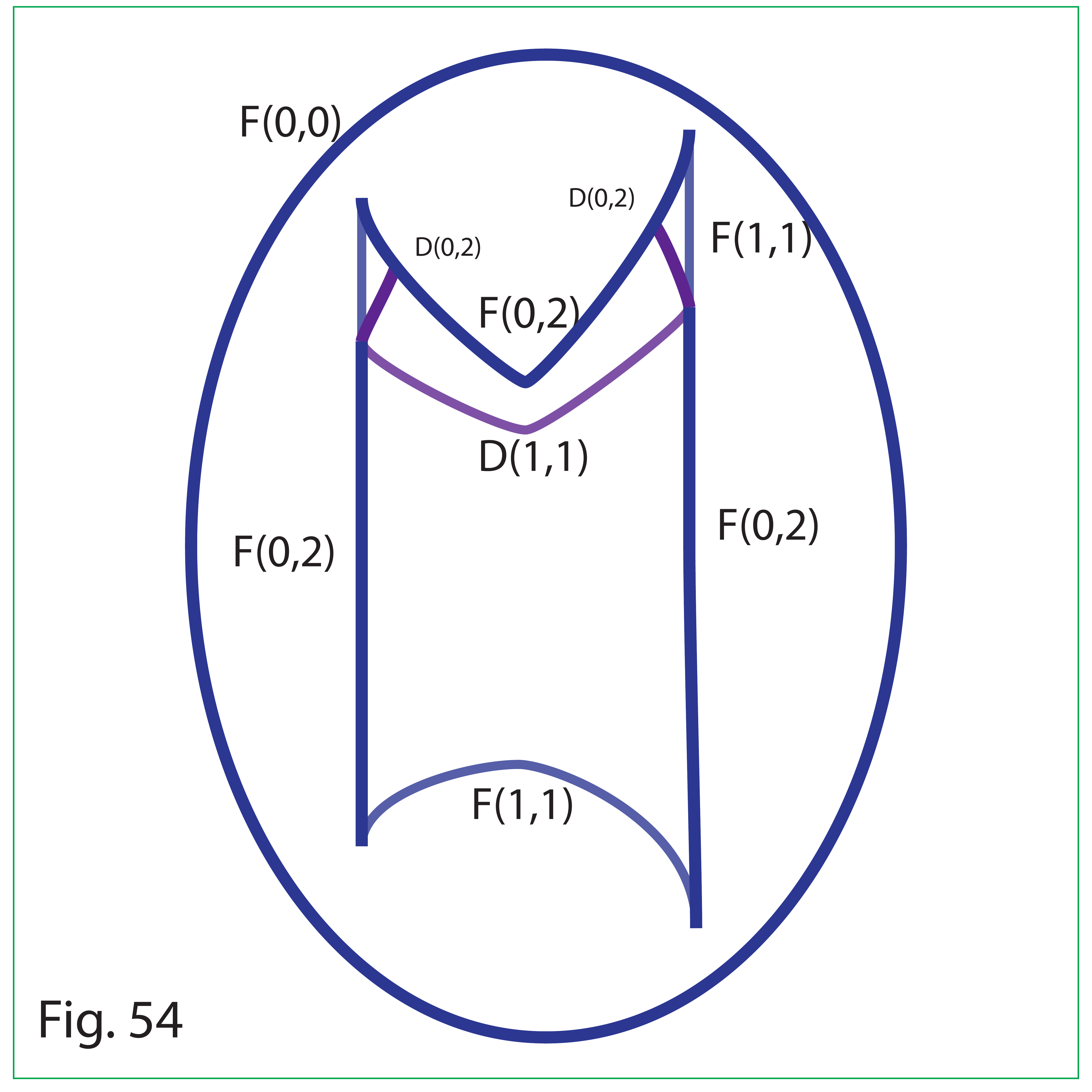}
\includegraphics[width=2.25in]{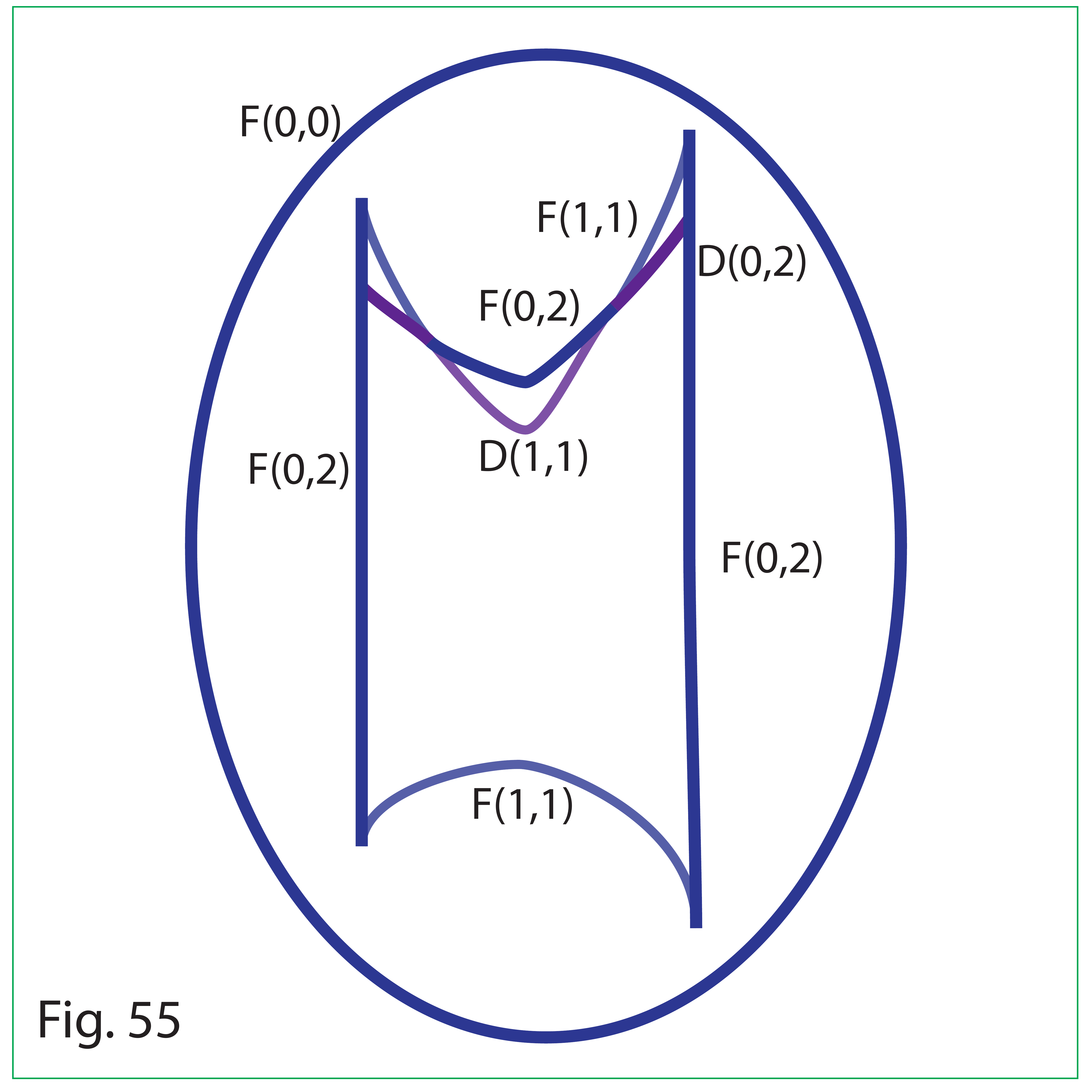}
\includegraphics[width=2.25in]{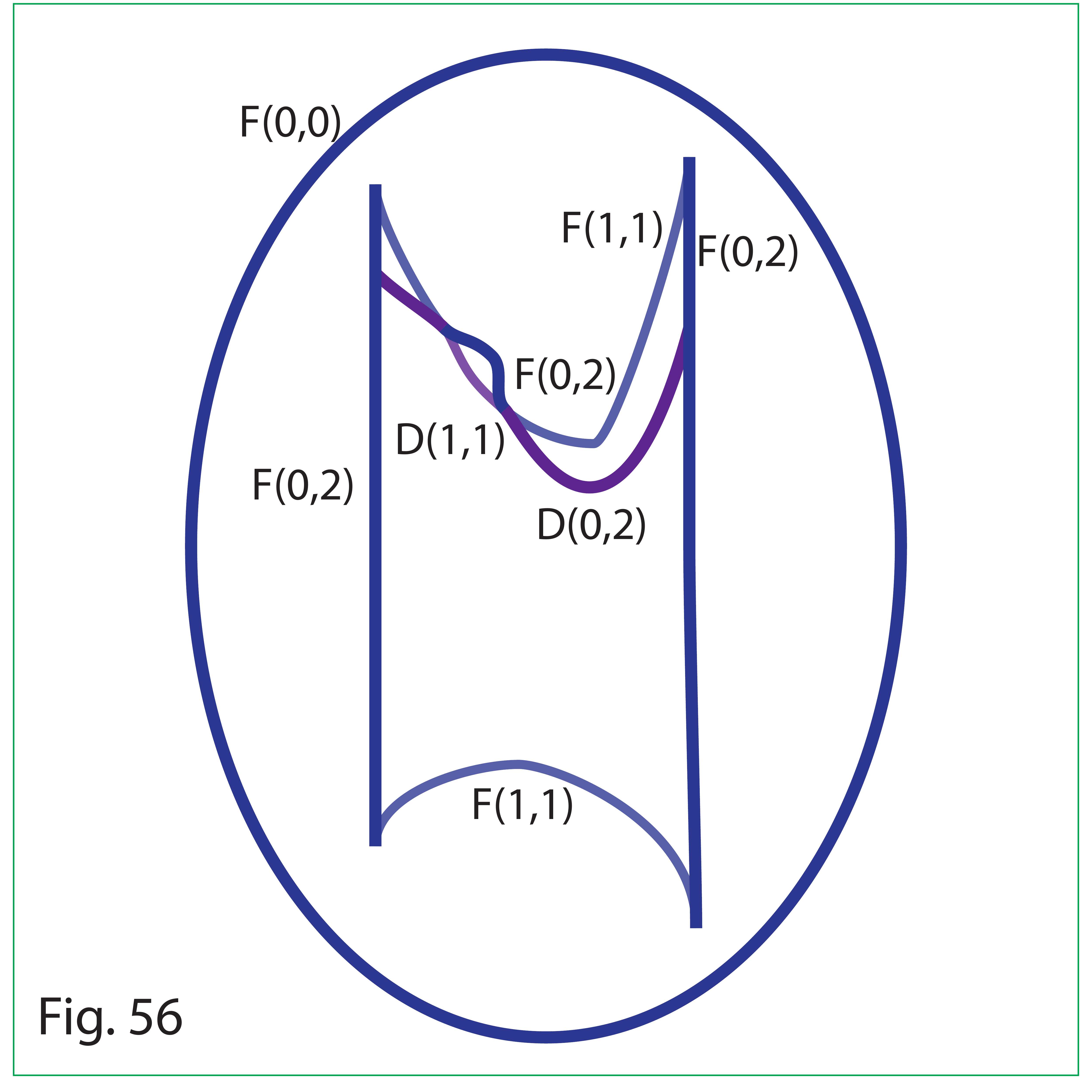}
\includegraphics[width=2.25in]{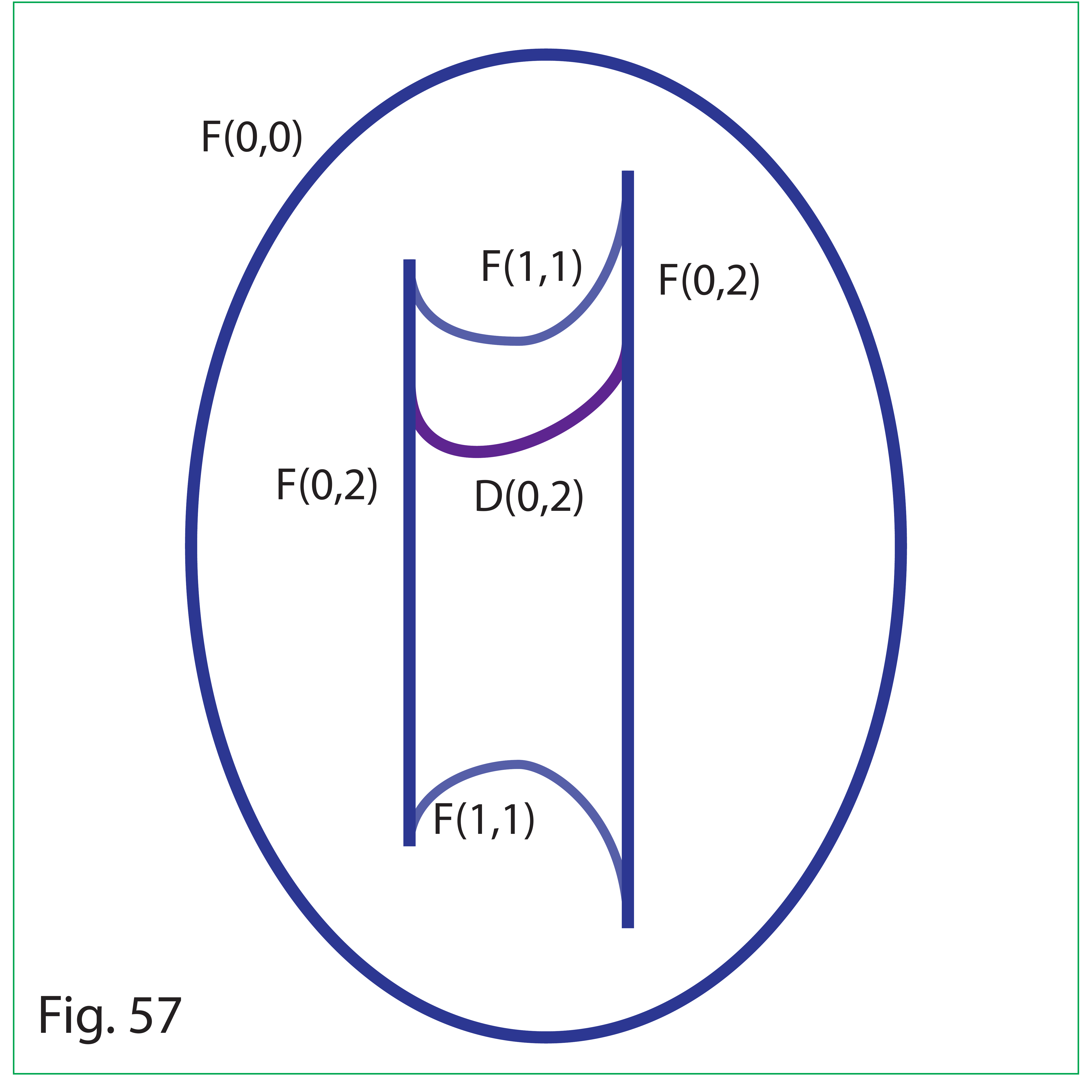}
\includegraphics[width=2.25in]{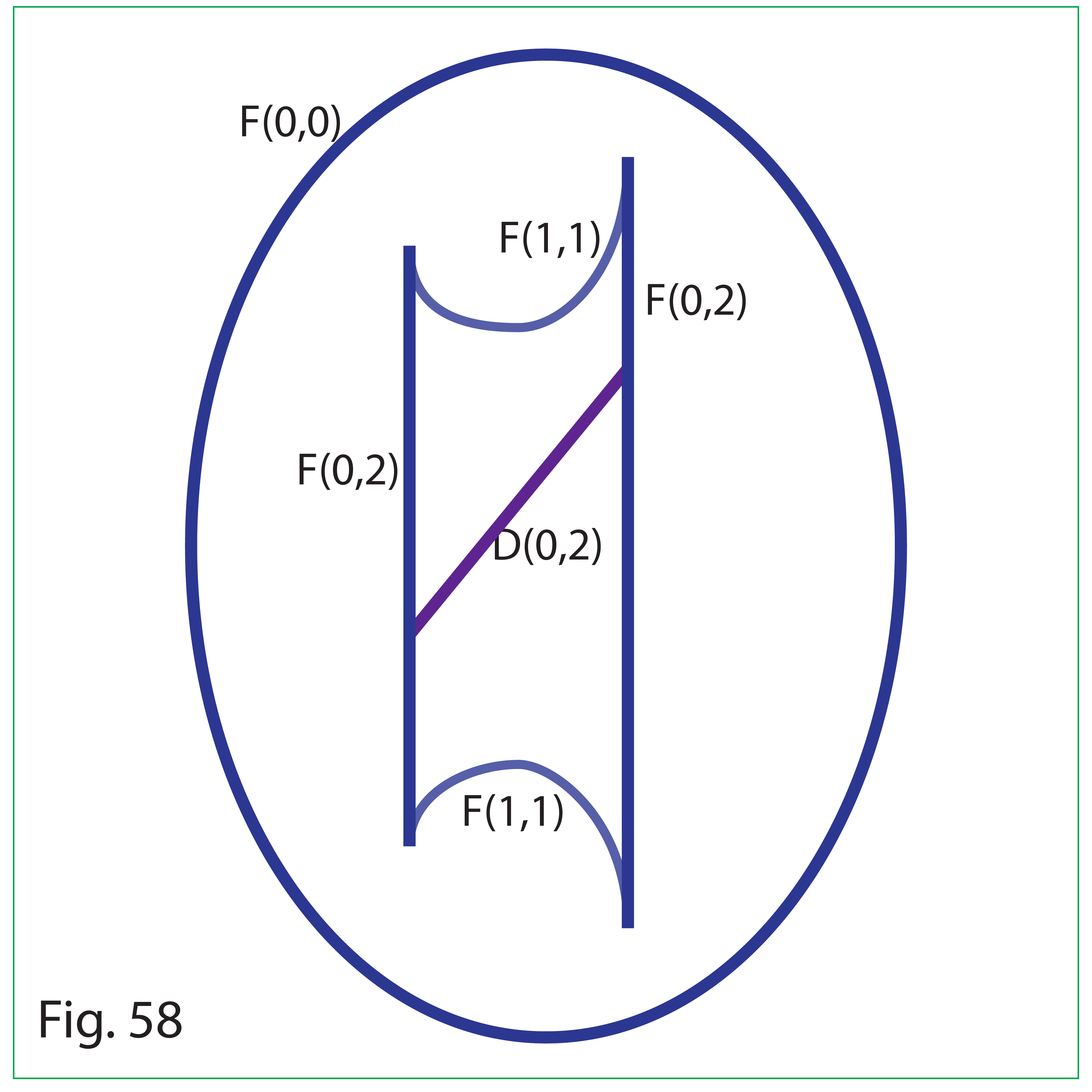}
\includegraphics[width=2.25in]{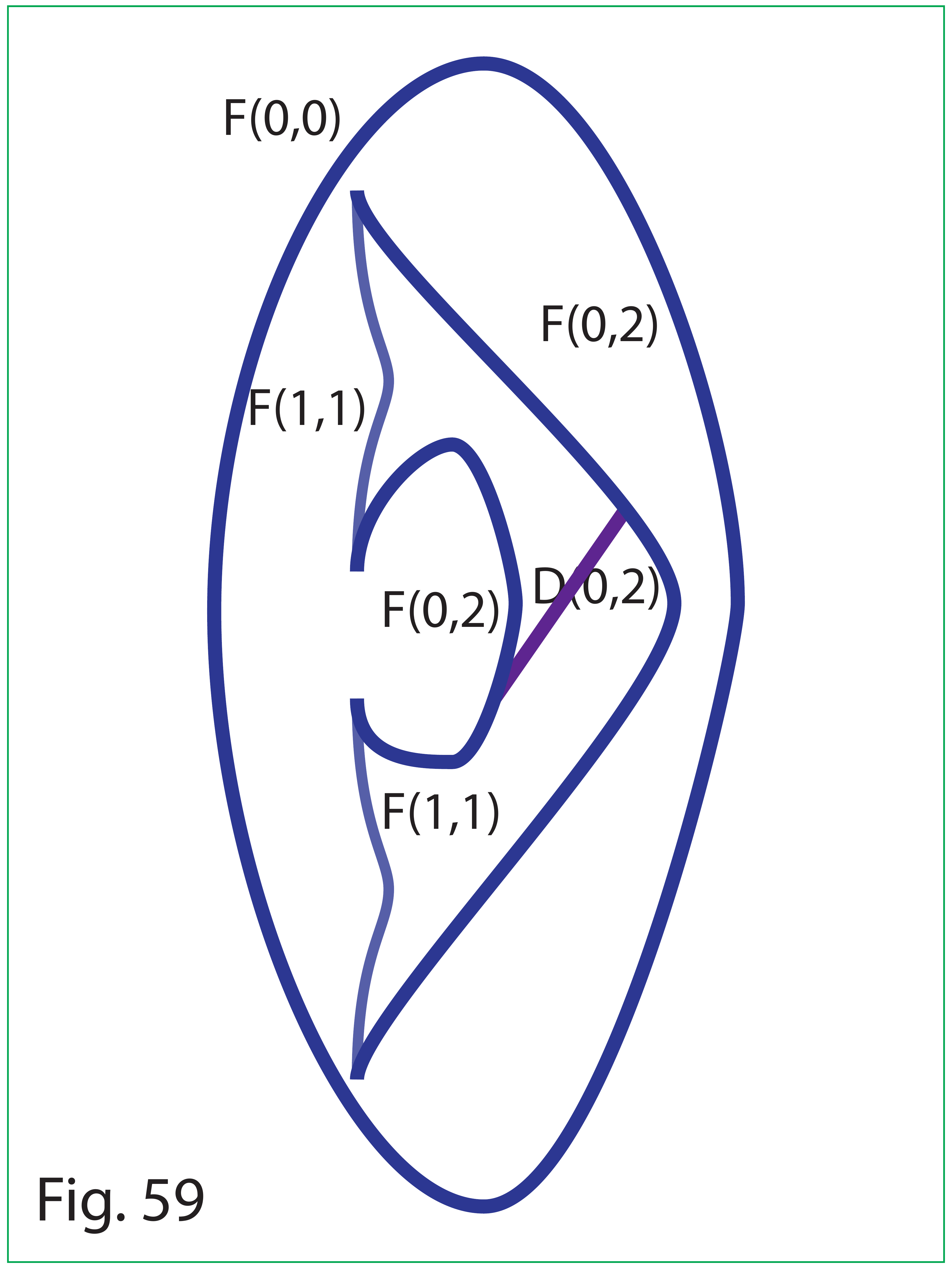}
\includegraphics[width=2.25in]{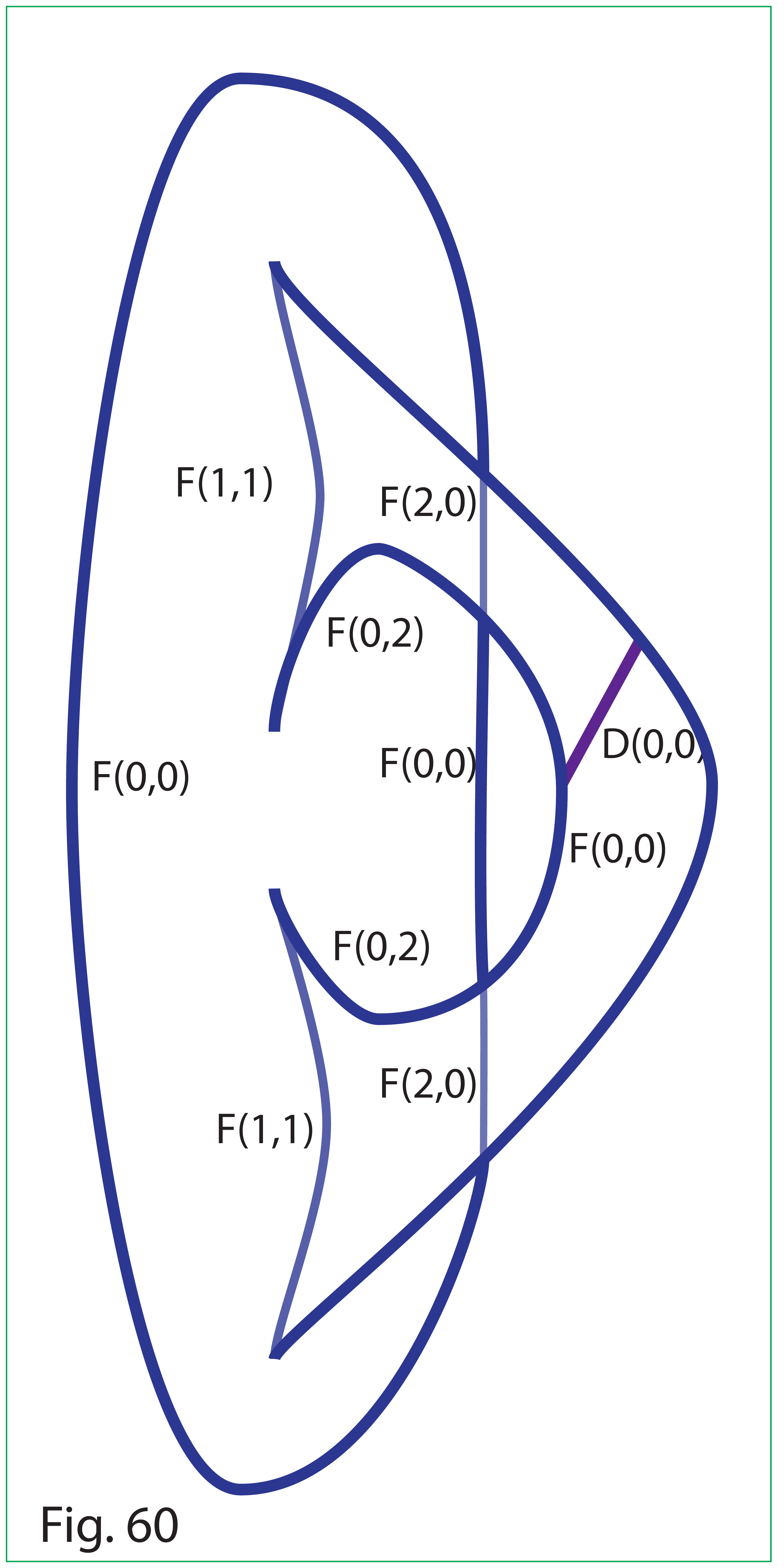}
\includegraphics[width=2.25in]{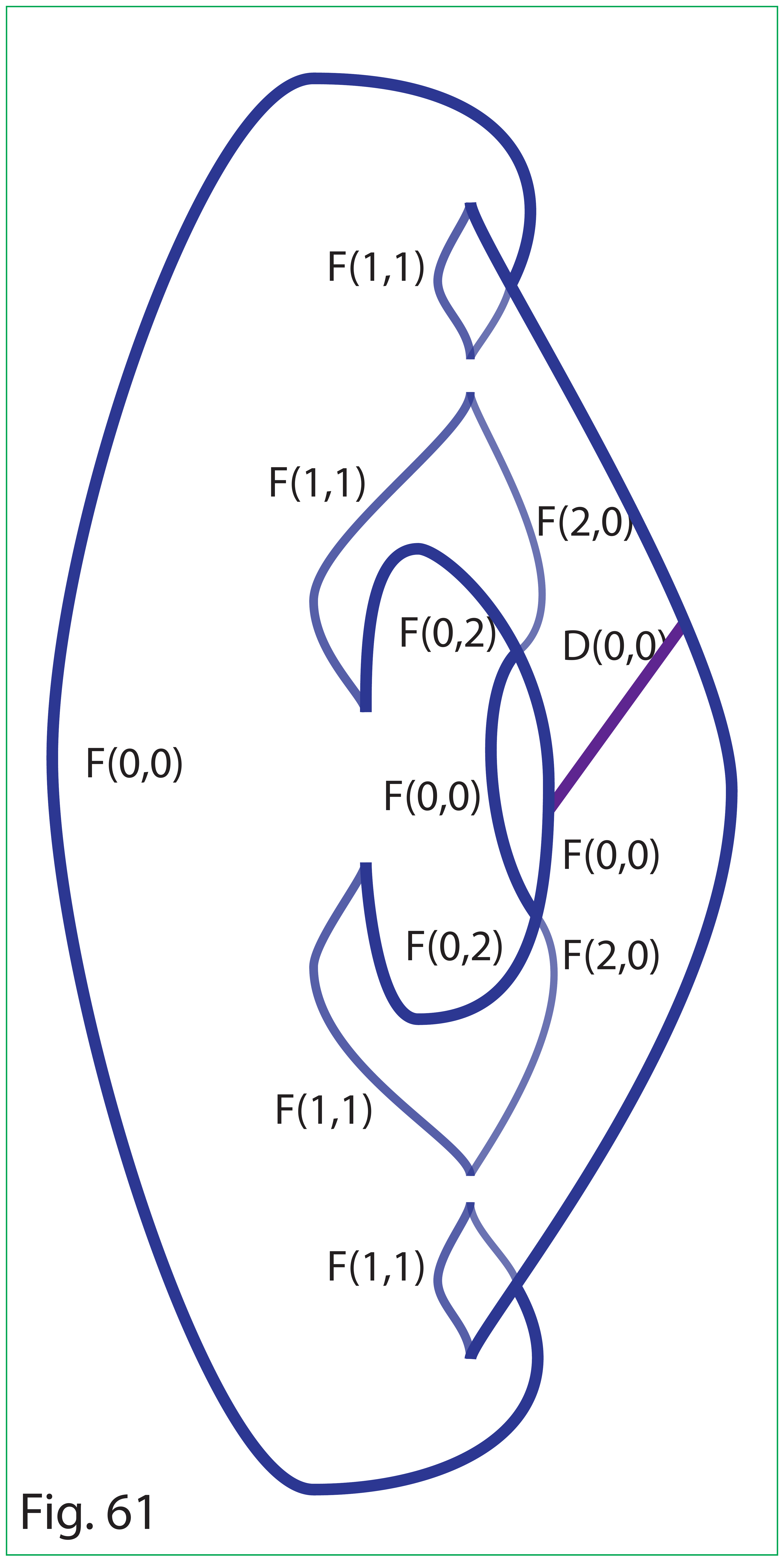}
\includegraphics[width=2.25in]{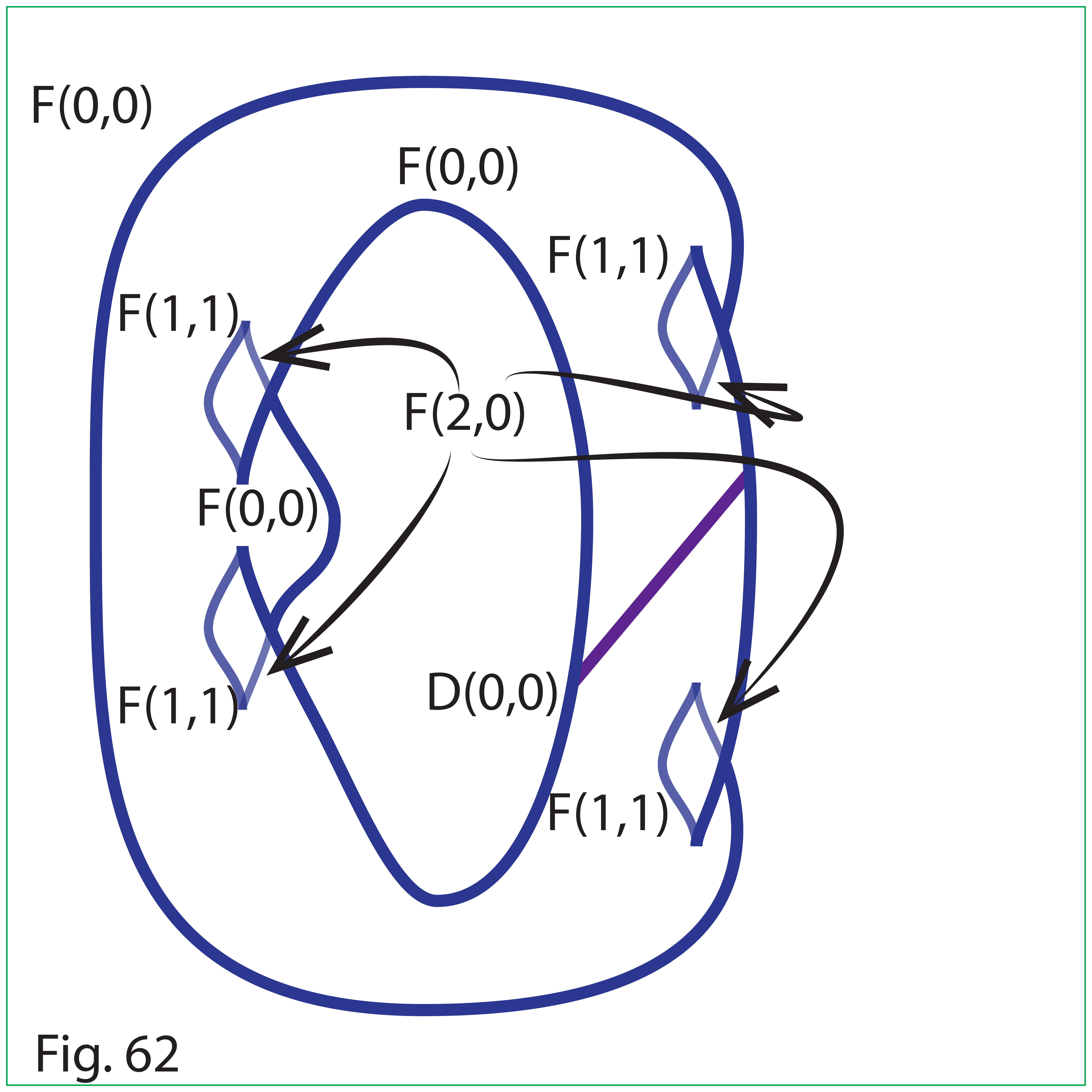}
\includegraphics[width=4.5in]{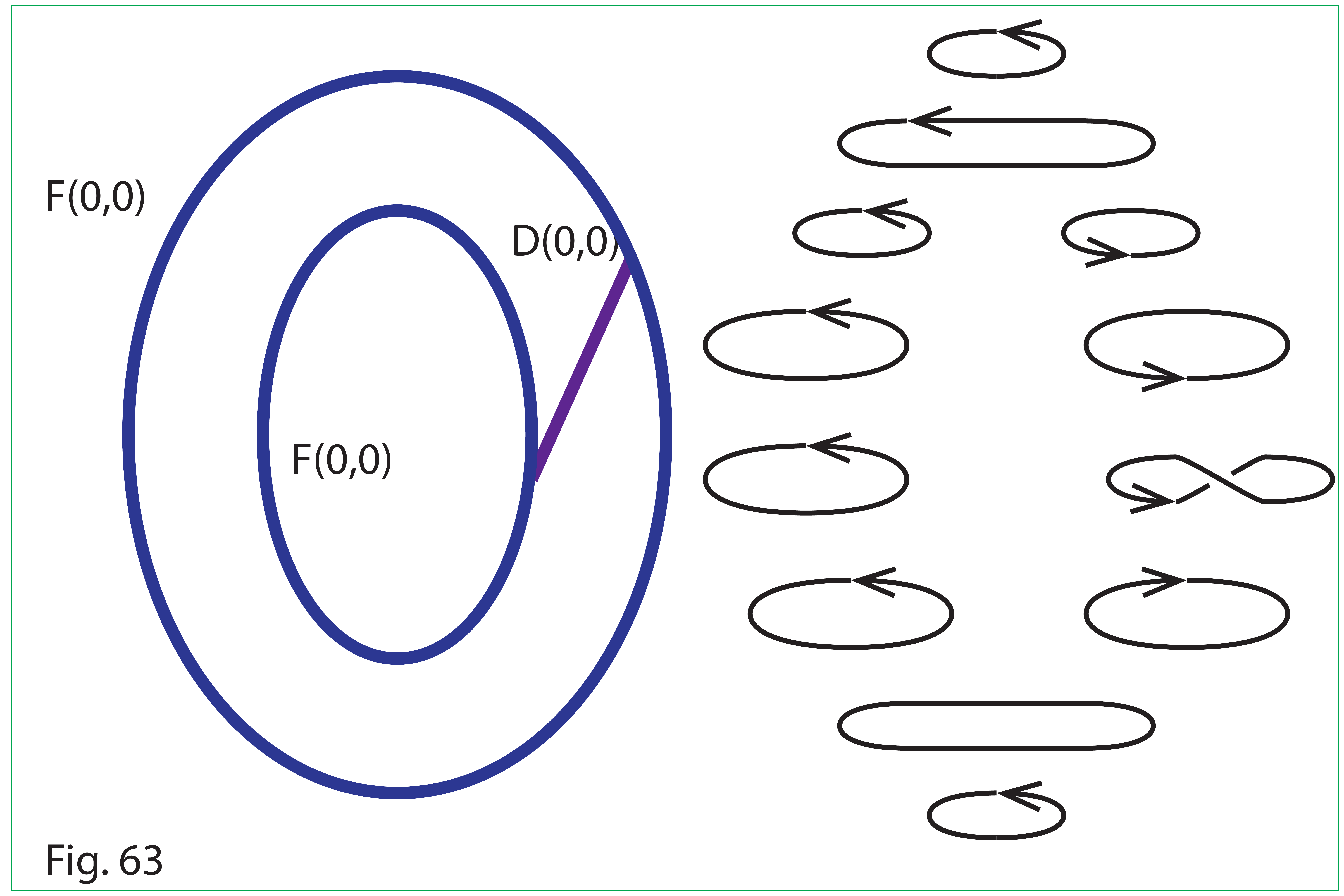}

Fig. 52 and Fig 53 have the same charts without labels. However the labels on the top arc of double points differ. The result of this difference is that the``neck" of the Klein bottle passes through the front face in Fig. 53, but it passes through the back face in Fig. 52. 
Clearly, each still in the movies that describe Fig. 52 could be flipped in a horizontal line to get a movie for which the neck passed through the front. However, the cusps at the top become reversed. We leave the reader an exercise to invert these cusps by using swallow-tail moves. Furthermore, to achieve such a reflected picture by means of movie moves can be quite complicated.

In Fig. 53 through Fig. 63, we demonstrate how to move from the standard Klein bottle to the ``torus with a wen" picture. We first illustrate the standard image via a movie to its right. The surface moves to the illustration in Fig. 54 via a horizontal confluence of branch points (Fig. 53, page 86 \cite{CS:Book}) together with a branch point passing over a saddle
 (Fig. 41, p. 81 \cite{CS:Book}). To get to the illustration in Fig. 55,  we apply two cases of branch points passing through cusps. (Fig. 42, p. 82 \cite{CS:Book}). The chart in Fig. 56 follows by passing a double point over a fold near a saddle (Fig. 44, p. 82 \cite{CS:Book}). {}From Fig. 56 to Fig. 57, we remove redundant $\psi$-bounces (Fig. 46, p. 83 \cite{CS:Book}).  A minimal point of the double point set then (Fig. 57 to Fig. 58) is pushed through a branch point (Fig.50, p. 85 \cite{CS:Book}). A  pair of horizontal cusps (Fig. 38, p. 80, \cite{CS:Book}) occurs to create the illustration in Fig. 59. At this stage it should be clear that a handle with a wen is attached to a sphere. To achieve Fig. 60, the folds are pushed to the outside right (See Fig. 59, p. 89 on the top right \cite{CS:Book}). A pair of beak-to-beak hyperbolic confluence of cusps (Fig. 36, p. 79 \cite{CS:Book}) occurs, folds move over optima on the fold set (Fig. 59, lower right, \cite{CS:Book}), and finally four swallow tails (Fig. 39, p. 81 \cite{CS:Book}) occur (Fig. 61 to Fig. 62)
 
 In Fig. 63 a movie is given. The orientations on the individual curves are followed up until the orientation reversing handle is attached. The orientations on the right-most figure-$8$ determine which fold the two branch points occur.

\section{Connected sums}
\label{knotted}

So far, all of the non-orientable surfaces that we have seen that are embedded in $4$-space are classified as unknotted since they are either cross-cap embeddings or ambiently isotopic to the connected sum of a pair of cross-caps. In the case of the standard (Acme) Klein bottle, one can in fact determine a disk bundle over the circle that is embedded in $4$-space and is bounded by the surface in question.

\subsection{The connected sum of the spun trefoil and a cross-cap}
One possible way of constructing a knotted example is to take a knotted sphere and form the connected sum of this with a standard cross-cap. The result may, or may not, be knotted. In Fig. 64, we present a movie that was created by forming the connected sum of the spun trefoil with a single cross cap. The movie was obtained by taking Kamada's \cite{Kamada:JKTR,Kamada:book} braid chart of the spun trefoil and perturbing the black vertices. At the ``bottom-most" black vertex, we changed the saddle band to one that reverses orientation.  To see that the surface is knotted, we have three-colored the movie. Specifically, such a coloring corresponds to representation of the fundamental group into the symmetric group $S_3$ on three letters $\{1,2,3\}$. The colors correspond to transpositions 
$1 \leftrightarrow (23)$, $2\leftrightarrow (13)$, and $3\leftrightarrow(12)$. In fact, one can easily compute the fundamental group of this example and demonstrate that it is precisely $S_3$.  The colors are indicated as labels from the set $\{1,2,3\}$ indicated on each of the arcs of the stills in each movie. 

\includegraphics[width=4.75in]{spunmovie}

\includegraphics[width=4in]{spunchart}

The chart of the previous movie is indicated in Fig. 65. Double lines, fold lines, and their respective depths are indicated. The vertically aligned sequences of positive integers indicate ``color vectors" for this surface. These vectors indicate the colors on the sheets of the surface in the $z$-direction. Thus a line of the form $\{(x_0,y_0,tz,): t\in[0,1]\}$ will pierce the surface many times and the sheets at which these piercings occur have colors, for example, $1,2,3,2,1,1$ from the $z=1$ level to the $z=0$ level, or as they are encountered by the observer. The color vectors of the regions of this chart are determined by the colors in the associated movie.

\subsection{The connected sum of the $2$-twist-spun trefoil and a cross-cap}
As in the previous case, we started from Kamada's braid chart for a $2$-knot, and modified one of the black vertices. In the current case the $2$-knot is the $2$-twist-spun trefoil. We illustrate the chart of this example in Fig. 66. Again, the black vertices have been perturbed to fulfill the current notion of a chart. In a neighborhoods of the branch points, we have not explicitly given indices for the folds or the double curves. Instead the movie that follows the chart indicates the perturbations explicitly. Also, line width and opacity helps indicate depth. 

The movie (Fig. 67) of this example indicates that the this knotted projective plane also can be $3$-colored. Thus it, too, is non-trivial. An open question that JSC and SK  are pursuing is whether or not this example is the same as the previous example. 

 \begin{center}
\includegraphics[width=5in]{fat}
{\Large{\sf The chart of the $2$-twist-spun trefoil connect sum a cross-cap}}
\end{center}

\begin{center}
\includegraphics[width=4.9in]{twtwshmovielong}
 {\Large{\sf The movie of the previous chart}}
\end{center}

\newpage

\subsection{The connected sum of the $3$-twist-spun trefoil and a cross-cap}

The next two illustrations depict the connected sum of the $3$-twist-spun trefoil  and a cross-cap. In the movie description that is given, Kamada's braid techniques are used. Thus there are seven times at which a single crossing is inserted or deleted. These correspond to Kamada's black vertices. However, in the chart description, they each have been perturbed to the generic situation in which three cusps, a saddle, and a branch point occur. The indices on the fold and double point sets are read from the left of the stills in the movies. Beware, these readings are not consistent with the previous examples. The inconsistency is caused because the side of the diagrams on which braid closures occur differ from those in the previous examples.

It is a standard conjecture that the $3$-twist-spun trefoil connected sum with a cross-cap is isotopic to a standard cross-cap. We believe the movie move techniques as applied to chart descriptions can be used to demonstrate explicit isotopies. However, our attempts to give such isotopies thus far have, sadly, not yielded any results. 

\includegraphics[width=5in]{thrtw4bcck}

\includegraphics[width=5in]{three4BCCK}

\subsection{General Question}

Let $K$ denote a classical knot. Let $\tau_n(K)$ denote Zeeman's $n$-twist-spinning construction applied to $K$. Let $T_n(\pm) K = \tau_n(K) \sharp CC(\R P^2(\pm))$ denote the connected sum of the $n$-twist-spin of $K$ and the $\pm$-cross-cap embedding of the projective plane. A well-known problem is to show in general that   
$$T_n(\pm)(K) \cong T_{n+2}(\pm)(K)$$ 
where $\cong$ denotes ambiently isotopic. The isotopy is believed since the two surfaces have isomorphic fundamental groups. However, it is possible that these knotted surfaces are distinct. The problem of their possible isotopy was brought to our attentions by Paul Melvin at a talk that he delivered at Knots in Washington. 

Our impression of this conjectural isotopy is that movie move and chart move techniques will be able to detect the isotopies in specific cases. The authors JSC and SK are currently preparing a sequel that addresses many of the problems that are associated to this question. At the point of this writing, we do not yet have a solution.

\medskip

\begin{flushleft}
Yongju Bae \\ 
Department of Mathematics \\ 
Kyungpook National University  \\
Daegu, Korea 
E-mail address: {\tt  ybae@knu.ac.kr}
\end{flushleft}

\begin{flushleft}
J. Scott Carter \\ 
Department of Mathematics \\ 
University of South Alabama \\ 
Mobile, AL 36688 \\
USA\\ 
E-mail address: {\tt carter@southalabama.edu} 
\end{flushleft}

\begin{flushleft}
Seonmi Choi \\ 
Department of Mathematics \\ 
Kyungpook National University  \\
Daegu, Korea 
E-mail address: {\tt  csm123c@gmail.com}
\end{flushleft}

\begin{flushleft}
Sera Kim \\ 
Department of Mathematics \\ 
University of South Alabama \\ 
Mobile, AL 36688 \\
USA\\ 
E-mail address: {\tt srkim85@gmail.com} 
\end{flushleft}

\end{document}